\documentclass[11pt]{article}
\usepackage{amssymb}
\usepackage{amsmath}
\usepackage{times}
\usepackage{graphicx}
\usepackage[all]{xy}
\usepackage{xcolor}
\usepackage{accents}

\title{Continuous colorings on compact spaces\indent}
\author{No\'e DE RANCOURT, Dominique LECOMTE$^1$, Miroslav ZELEN\'Y$^2$}
\date{\today}

\def\ufootnote#1{\let\savedthfn\thefootnote\let\thefootnote\relax
\footnote{#1}\let\thefootnote\savedthfn\addtocounter{footnote}{-1}}

\newcommand{\ca}{{\bf\Pi}^{1}_{1}}

\newcommand{\borone}{{\bf\Delta}^{0}_{1}}

\newcommand{\bormtwo}{{\bf\Pi}^{0}_{2}}

\newtheorem{thm} {Theorem} [section]
\newtheorem{defi} [thm] {Definition}
\newtheorem{cor} [thm] {Corollary}
\newtheorem{lem} [thm] {Lemma}
\newtheorem{prop} [thm] {Proposition}

\newtheorem{them} {Theorem} [subsection]

\newtheorem{coro} [them] {Corollary}
\newtheorem{lemm} [them] {Lemma}

\begin{document}

\maketitle

\centerline{to appear in Bollettino dell'Unione Matematica Italiana}\bigskip

\centerline{$\bullet$ Universit\' e de Lille, CNRS, UMR 8524, Laboratoire Paul Painlev\'e}

\centerline{F-59 000 Lille, France}

\centerline{nderancour@univ-lille.fr}\bigskip

\centerline{$\bullet^1$ Universit\' e Paris 6, Institut de Math\'ematiques de Jussieu-Paris Rive Gauche,}

\centerline{Projet Analyse Fonctionnelle, Couloir 16-26, 4\`eme \'etage, Case 247,}

\centerline{4, place Jussieu, 75 252 Paris Cedex 05, France}

\centerline{dominique.lecomte@upmc.fr}\medskip

\centerline{$\bullet^1$ Universit\'e de Picardie, I.U.T. de l'Oise, site de Creil,}

\centerline{13, all\'ee de la fa\"\i encerie, 60 107 Creil, France}\bigskip

\centerline{$\bullet^2$ Charles University, Faculty of Mathematics and Physics, Department of Mathematical Analysis}

\centerline{Sokolovsk\'a 83, 186 75 Prague, Czech Republic}

\centerline{zeleny@karlin.mff.cuni.cz}\bigskip\bigskip\bigskip\bigskip\bigskip

\ufootnote{{\it 2020 Mathematics Subject Classification.}~Primary: 03E15, Secondary: 54H05, 37B05, 37B10}

\ufootnote{{\it Keywords and phrases.}~analytic complete, antichain, basis, Cantor-Bendixson rank, compact, continuous coloring, continuous homomorphism, flip-conjugacy, graph, subshift}

\ufootnote{{\it Acknowledgements.}~No\'e de Rancourt acknowledges support from the Labex CEMPI (ANR-11-LABX-0007-01)}

\noindent {\bf Abstract.} We study several natural classes of graphs on a zero-dimensional metrizable compact space having no continuous coloring. We compare these graphs with the quasi-order $\preceq^i_c$ associated with injective continuous homomorphisms. We prove the existence of an antichain basis for these classes. We determine the size of such an antichain  basis. We provide a concrete antichain basis when there is a countable one. We also provide some related quasi-orders and equivalence relations which are analytic complete as sets.\bigskip\bigskip\bigskip\bigskip\bigskip

\vfill\eject\baselineskip=12.9pt

\section{$\!\!\!\!\!\!$ Introduction}\indent

 The present work is the continuation of the study of continuous 2-colorings initiated in [L]. All our relations will be binary. A {\bf coloring} of a relation $R$ on a set $X$ is a map $c$ from $X$ into a set $\kappa$ with the property that $c(x)\!\not=\! c(y)$ if $(x,y)\!\in\! R$. We will call this a $\kappa$-coloring. In practice, $\kappa$ will be a countable cardinal, equipped with the discrete topology. We say that $R$ (or $(X,R)$) is a {\bf graph} if $R$ is symmetric and does not meet the {\bf diagonal} $\Delta (X)\! :=\!\{ (x,x)\mid x\!\in\! X\}$ of $X$. We set $R^{-1}\! :=\!\{ (x,y)\!\in\! X^2\mid (y,x)\!\in\! R\}$, and 
$s(R)\! :=\! R\cup R^{-1}$ is the {\bf symmetrization} of $R$. We compare our relations with the following quasi-order:\medskip

\centerline{$(X,R)\preceq^i(Y,S)\Leftrightarrow\exists h\! :\! X\!\rightarrow\! Y\mbox{ injective with }R\!\subseteq\! (h\!\times\! h)^{-1}(S).$}\medskip

\noindent If this holds, then we say that $h$ is an injective {\bf homomorphism} from 
$(X,R)$ into $(Y,S)$. In the present article, we work with the quasi-order $\preceq^i_c$ associated with injective continuous homomorphisms. All our topological spaces will be zero-dimensional, except where indicated, to ensure the existence of enough continuous functions between them. We write $(X,R)\prec^i_c(Y,S)$ when $(X,R)\preceq^i_c(Y,S)$ and $(Y,S)\not\preceq^i_c(X,R)$. The material in [L] shows that the structure of 
$\preceq^i_c$ is complex on a number of classes of graphs. Recall that a {\bf basis} for a quasi-order $(\mathcal{Q},\leq )$ is a subclass $\mathcal{B}$ of $\mathcal{Q}$ such that any element of $\mathcal{Q}$ is $\leq$-above an element of $\mathcal{B}$. We are interested in basis as small as possible for the inclusion, which means that their elements are pairwise $\leq$-incomparable (if this last property is satisfied, then we say that we have a $\leq$-{\bf antichain}). Note that an antichain basis is always made of minimal elements of the considered class. Conversely, let 
$\equiv^i_c\ :=\ \preceq^i_c\cap\ (\preceq^i_c)^{-1}$ be the equivalence relation associated with $\preceq^i_c$. Note that we can derive an antichain basis from a basis made of minimal elements by choosing an element in each $\equiv^i_c$-equivalence class, using the axiom of choice if necessary.\medskip

\noindent $\bullet$ Theorem 1.10 in [L] shows that there is no antichain basis for the class of graphs on a zero-dimensional metrizable compact space ({\bf 0DMC} for short; we will also use similar abbreviations like MC or 0DM) having no continuous 2-coloring. This theorem in fact gives the same result for graphs $(X,R)$ with $R$ countable. The situation is completely different for {\bf closed graphs}, which leads to the first class we study. A compactness argument shows that any closed graph on a 0DMC space has a continuous $\aleph_0$-coloring.

\begin{thm} \label{mainclosed} Let $\kappa\! <\!\aleph_0$ be a cardinal.\smallskip

\noindent (a) There is a $\preceq^i_c$-basis made of minimal elements for the class of closed graphs on a 0DMC space having no continuous $\kappa$-coloring.\smallskip

\noindent (b) Such a basis can be $\{ (1,\emptyset )\}$ if $\kappa\! =\! 0$, $\{ (2,\{ (0,1),(1,0)\} )\}$ if $\kappa\! =\! 1$, and has size $2^{\aleph_0}$ if $\kappa\!\geq\! 2$.\end{thm}

\noindent $\bullet$ The case of graphs induced by a function has been considered since the very beginning of the study of definable colorings in [K-S-T], and also in [Co-M], [L], [P] and [T-V] for instance. If\medskip
 
\centerline{$f\! :\!\mbox{Domain}(f)\!\subseteq\! X\!\rightarrow\!
\mbox{Range}(f)\!\subseteq\! X$}\medskip

\noindent is a partial function, then the graph {\bf induced} by $f$ is 
$G_f\! :=\! s\big(\mbox{Graph}(f)\big)\!\setminus\!\Delta (X)$. The end of Section 9 in [L] shows that there is no antichain basis for the class of graphs induced by a partial homeomorphism on a 0DMC space with countable domain having no continuous 2-coloring. So we will focus on {\bf graphs induced by a total homeomorphism}.\medskip

 The following example was essentially introduced in [L-Z]. We consider a converging sequence with its limit in the Cantor space $2^\omega$, for instance $\mathbb{X}_1\! :=\!\{ 0^n1^\infty\mid n\!\in\!\omega\}\cup\{ 0^\infty\}$, which is a countable MC space. We define a homeomorphism $f_1$ of $\mathbb{X}_1$ by $f_1(0^\infty )\! :=\! 0^\infty$ and 
${f_1(0^{2n+\varepsilon}1^\infty )\! :=\! 0^{2n+1-\varepsilon}1^\infty}$ if $\varepsilon\!\in\! 2$. We will see that 
$(\mathbb{X}_1,G_{f_1})$ has no continuous $\aleph_0$-coloring.

\vfill\eject\baselineskip=12.9pt

\begin{thm} \label{mainhomeo} Let $\kappa\!\leq\!\aleph_0$ be a cardinal.\smallskip

\noindent (a) There is a $\preceq^i_c$-basis made of minimal elements for the class of graphs, induced by a homeomorphism of a 0DMC space, having no continuous $\kappa$-coloring.\smallskip

\noindent (b) ([L], Theorems 1.17 (b) and 1.13 (d)) Such a basis can be $\{ (1,G_{0\mapsto 0})\}$ if $\kappa\! =\! 0$, 
$\{ (2,G_{\varepsilon\mapsto 1-\varepsilon})\}$ if $\kappa\! =\! 1$, has size $2^{\aleph_0}$ if $\kappa\! =\! 2$, and can be 
$\{ (\mathbb{X}_1,G_{f_1})\}$ if $\kappa\!\geq\! 3$.\end{thm}

 This result can be refined when $\kappa\! =\! 2$ if we consider the Cantor-Bendixson rank of the considered spaces. Recall from 6.C in [K1] that if $X$ is a topological space, then the {\bf Cantor-Bendixson derivative} of $X$ is 
$X'\! :=\!\{ x\!\in\! X\mid x\mbox{ is a limit point of }X\}$. The {\bf iterated Cantor-Bendixson derivatives} are defined by 
$X^0\! :=\! X$, $X^{\alpha +1}\! :=\! (X^\alpha )'$, and, when $\lambda$ is a limit ordinal, 
$X^\lambda\! :=\!\bigcap_{\alpha <\lambda}~X^\alpha$. Note that if $f$ is a homeomorphism of $X$, then all the derivatives are $f$-{\bf invariant}, i.e., $f[X^\alpha ]\! =\! X^\alpha$ if $\alpha$ is an ordinal. If $X$ is a countable MC space, then the 
{\bf Cantor-Bendixson rank} of $X$ is the least countable ordinal $\alpha_0$ such that 
$X^{\alpha_0}\! =\!\emptyset$. If moreover $X$ is nonempty, then this rank is a successor ordinal, by compactness. More generally, the {\bf Cantor-Bendixson rank} of a Polish space is the least countable ordinal $\alpha_0$ such that $X^\alpha\! =\! X^{\alpha_0}$ for each $\alpha\!\geq\!\alpha_0$, so that the Cantor space $2^\omega$ has Cantor-Bendixson rank zero.\medskip

 The following examples are of particular interest here.\medskip

\noindent - The odd cycles $(2q\! +\! 3,C_{2q+3})$, for $q\!\in\!\omega$. In this case, the formula 
$f(i)\! :=\! (i\! +\! 1)\mbox{ mod }(2q\! +\! 3)$ defines a homeomorphism of the discrete MC space $2q\! +\! 3$, whose Cantor-Bendixson rank is one. We set $C_{2q+3}\! :=\! G_f$, and the fact that $(2q\! +\! 3,C_{2q+3})$ has no (continuous) 2-coloring is classical.\medskip

\noindent - $\mathbb{X}_1$ has Cantor-Bendixson rank two.\medskip

\noindent - We also consider subshifts, which are particular dynamical systems widely studied in symbolic dynamics. We refer to the book [Ku] for basic notions and definitions.

\begin{defi} Let $A$ be a finite set of cardinality at least two.\smallskip

\noindent (a) The {\bf shift map} $\sigma\! :\! A^\mathbb{Z}\!\rightarrow\! A^\mathbb{Z}$ is defined by the formula 
$\sigma (\alpha )(k)\! :=\!\alpha (k\! +\! 1)$.\smallskip

\noindent (b) A {\bf two-sided subshift} is a closed $\sigma$-invariant subset $\Sigma$ of $A^\mathbb{Z}$.\end{defi}

 The restriction of the homeomorphism $\sigma$ to a two-sided subshift $\Sigma$ induces a graph 
$(\Sigma ,G_{\sigma_{\vert\Sigma}})$ that we will denote by $(\Sigma ,G_\sigma )$. If $f$ is a bijection of the set $X$ and 
$x\!\in\! X$, then the $f$-{\bf orbit} of $x$ is\medskip

\centerline{$\mbox{Orb}_f(x)\! :=\!\{ f^k(x)\mid k\!\in\!\mathbb{Z}\}$}\medskip

\noindent (also denoted by $\mbox{Orb}(x)$ when the context is clear). If 
$x\!\in\! A^{-\omega}$ and $y\!\in\! A^\omega$, then 
$z\! :=\! x\!\cdot\! y\!\in\! A^\mathbb{Z}$ is defined by $z(i)\! :=\! y(i)$ and 
$z(-i\! -\! 1)\! :=\! x(-i)$ when $i\!\in\!\omega$. If 
$w\!\in\! A^{<\omega}\!\setminus\!\{\emptyset\}$, then $w^{-\infty}\! :=\!\cdots\! ww$ is in $A^{-\omega}$, 
$w^\infty\! :=\! ww\!\cdots$ is in $A^\omega$ and $w^\mathbb{Z}\! :=\! w^{-\infty}\!\cdot\! w^\infty$. Note that 
$(2q\! +\! 3,C_{2q+3})$ can be seen as a two-sided subshift by putting 
${}_{2q+3}\Sigma\! :=\!\mbox{Orb}_\sigma\Big(\big( 0\!\cdots\! (2q\! +\! 2)\big)^\mathbb{Z}\Big)$. Recall that if $X,Y$ are metrizable compact spaces and $f,g$ are homeomorphisms of $X,Y$ respectively, then $(X,f),(Y,g)$ (or $f,g$) are 
{\bf conjugate} (resp., {\bf flip}-{\bf conjugate}) if there is a homeomorphism $\varphi\! :\! X\!\rightarrow\! Y$ such that 
$\varphi\!\circ\! f\! =\! g\!\circ\!\varphi$ (resp., $\varphi\!\circ\! f\! =\! g\!\circ\!\varphi$ or 
$\varphi\!\circ\! f\! =\! g^{-1}\!\circ\!\varphi$). We will see that 
$(\mathbb{X}_1,G_{f_1})$ is not conjugate to the shift of a two-sided subshift. We set\medskip

\centerline{$\mathcal{P}\! =\!\big\{ {\bf p}\! :=\! 
(l,\lambda_0,\cdots ,\lambda_l,m,\varepsilon_0,\cdots ,\varepsilon_{l-1})\!\in\!\omega^{l+3}\!\times\! 2^l\mid
\forall i\!\leq\! l~\lambda_i\! >\! 0\mbox{ is even and }m\! <\!\lambda_0\mbox{ is odd}\big\} .$}\medskip

\noindent We associate to each ${\bf p}\!\in\!\mathcal{P}$ a two-sided subshift as follows. We fix disjoint injective families of symbols $(a^i_j)_{i,j\in\omega}$ and $(b_i)_{i\in\omega}$, and set 
$A_{\bf p}\! :=\!\{ a^i_j\mid i\!\leq\! l\wedge j\! <\!\lambda_i\}\cup\{ b_i\mid i\! <\! m\}$, which is finite of cardinality at least two (in fact at least three).

\vfill\eject\baselineskip=12.9pt

 We then set, for $i\!\leq\! l$, $w_i\! :=\! 
a^i_0\!\cdots\! a^i_{\lambda_i-1}\!\in\! A_{\bf p}^{<\omega}\!\setminus\!\{\emptyset\}$, and define\medskip

\centerline{$\Sigma_{\bf p}\! :=\!\bigcup_{i\leq l}~~\mbox{Orb}_\sigma(w_i^\mathbb{Z})\cup
\bigcup_{i<l}~~
\mbox{Orb}_\sigma (w_{i+\varepsilon_i}^{-\infty}\!\cdot\! w_{i+1-\varepsilon_i}^{\infty})\cup
\mbox{Orb}_\sigma\big( w_l^{-\infty}\!\cdot\! b_0\!\cdots\! b_{m-1}(w_0^\infty )\big)\mbox{,}$}\medskip

\noindent a countable MC space with Cantor-Bendixson rank two. We will see that 
$(\Sigma_{\bf p},G_\sigma )$ has no continuous $2$-coloring. We set, for each 
${\bf p}\!\in\!\mathcal{P}$, $\mathcal{F}_{\bf p}\! :=\!\{ {\bf p}'\!\in\!\mathcal{P}\mid (\Sigma_{{\bf p}'},G_\sigma )\equiv^i_c(\Sigma_{\bf p},G_\sigma )\}$.\medskip
 
 For $\kappa\! =\! 2$, we prove the following.
 
\begin{thm} \label{main} (a) The family $\{ ({}_{2q+3}\Sigma ,G_\sigma )\mid q\!\in\!\omega\}\cup
\{ (\mathbb{X}_1,G_{f_1})\}\cup\{ (\Sigma_{\bf p},G_\sigma )\mid {\bf p}\!\in\!\mathcal{P}\}$ is a concrete $\preceq^i_c$-basis of size $\aleph_0$ for the class of graphs, induced by a homeomorphism of a countable MC space with Cantor-Bendixson rank at most two, having no continuous 2-coloring. Moreover, for each ${\bf p}\!\in\!\mathcal{P}$, the set 
$\mathcal{F}_{\bf p}$ is finite, and choosing $\mbox{min}_{\text{lex}}~\mathcal{F}_{\bf p}$ provides a $\preceq^i_c$-antichain basis of size $\aleph_0$.\smallskip

\noindent (b) (see [L], Theorem 1.17 (b)) If $\xi\!\geq\! 3$ is a countable ordinal, then there is a $\preceq^i_c$-basis made of minimal elements for the class of graphs, induced by a homeomorphism of a countable MC space with Cantor-Bendixson rank at most $\xi$, having no continuous 2-coloring, and any such basis must have size $2^{\aleph_0}$.\smallskip

\noindent (c) ([L], Theorem 1.15) If $\xi$ is a countable ordinal, then any $\preceq^i_c$-basis made of minimal elements for the class of graphs, induced by a homeomorphism of a 0DMC space with Cantor-Bendixson rank at most $\xi$, having no continuous 2-coloring, must have size $2^{\aleph_0}$.\end{thm}

 Note that the proof of Theorem \ref{mainclosed} (b) will show that it has no such refinement when $\kappa\!\geq\! 2$.\medskip

\noindent $\bullet$ The class of {\bf graphs induced by the shift of a two-sided subshift} is a natural subclass of the previous one. 

\begin{thm} \label{mainsubshift} Let $\kappa\!\leq\!\aleph_0$ be a cardinal.\smallskip

\noindent (a) There is a $\preceq^i_c$-basis made of minimal elements for the class of graphs, induced by the shift of a two-sided subshift, having no continuous $\kappa$-coloring.\smallskip

\noindent (b) Such a basis can be $\{ (\mbox{Orb}_\sigma (0^\mathbb{Z}),G_\sigma )\}$ if $\kappa\! =\! 0$, 
$\big\{\big(\mbox{Orb}_\sigma\big( (01)^\mathbb{Z}\big),G_\sigma\big)\big\}$ if $\kappa\! =\! 1$, and has size $2^{\aleph_0}$ if 
$\kappa\!\geq\! 2$.\end{thm}

 Here again, this result can be refined when $\kappa\! =\! 2$ if we consider the Cantor-Bendixson rank of the considered spaces. As $(\mathbb{X}_1,G_{f_1})$ is not conjugate to the shift of a two-sided subshift, we have to introduce some other examples.\medskip
 
 We set ${}_0\Sigma\! :=\!\mbox{Orb}_\sigma (0^\mathbb{Z})\cup\mbox{Orb}_\sigma (0^{-\infty}\!\cdot\! 10^\infty )\mbox{, }
{}_1\Sigma\! :=\!\mbox{Orb}_\sigma (0^\mathbb{Z})\cup\mbox{Orb}_\sigma (1^\mathbb{Z})\cup
\mbox{Orb}_\sigma (0^{-\infty}\!\cdot\! 1^\infty )$ and, for $q\!\in\!\omega$, 
${}_{2q+2}\Sigma\! :=\!\mbox{Orb}_\sigma (0^\mathbb{Z})\cup
 \mbox{Orb}_\sigma\big( (1,\cdots ,2q\! +\! 2)^\mathbb{Z}\big)\cup
 \mbox{Orb}_\sigma (0^{-\infty}\!\cdot\! (1,\cdots ,2q\! +\! 2)^\infty )$.
 
\begin{thm} \label{mainshift} (a) The family $\{ ({}_n\Sigma ,G_\sigma )\mid n\!\in\!\omega\}\cup
\{ (\Sigma_{\bf p},G_\sigma )\mid {\bf p}\!\in\!\mathcal{P}\}$ is a concrete $\preceq^i_c$-basis of size $\aleph_0$ for the class of graphs, induced by the shift of a countable two-sided subshift with Cantor-Bendixson rank at most two, having no continuous 2-coloring. Moreover, choosing $\mbox{min}_{\text{lex}}~\mathcal{F}_{\bf p}$ for each 
${\bf p}\!\in\!\mathcal{P}$ provides a $\preceq^i_c$-antichain basis of size $\aleph_0$.\smallskip

\noindent (b) ([L], Theorem 1.17 (b) and Corollary 10.12) If $\xi\!\geq\! 3$ is a countable ordinal, then there is a $\preceq^i_c$-basis made of minimal elements for the class of graphs, induced by the shift of a countable two-sided subshift with Cantor-Bendixson rank at most $\xi$, having no continuous 2-coloring, and any such basis must have size $2^{\aleph_0}$.\smallskip

\noindent (c) ([L], Corollary 10.12) If $\xi$ is a countable ordinal, then any $\preceq^i_c$-basis made of minimal elements for the class of graphs, induced by the shift of a two-sided subshift with Cantor-Bendixson rank at most $\xi$, having no continuous 2-coloring, must have size $2^{\aleph_0}$.\end{thm}
 
\vfill\eject\baselineskip=13.2pt
 
 A {\bf dynamical system} $(X,f)$ is given by a homeomorphism $f$ of a metrizable compact space $X$. If $X$ is homeomorphic to $2^\omega$, then we say that $(X,f)$ is a {\bf Cantor dynamical system}. A dynamical system (or 
$f$) is {\bf minimal} if $\mbox{Orb}_f(x)$ is dense in $X$ for each $x\!\in\! X$. The set of homeomorphisms of $2^\omega$ is denoted by $\mathcal{H}(2^\omega )$. It is a Polish group when equipped with the topology whose basic neighbourhoods of the identity are of the form $\{ h\!\in\!\mathcal{H}(2^\omega )\mid\forall i\! <\! n~~h[O_i]\! =\! O_i\}$, where 
$(O_i)_{i<n}$ ranges over all finite families of clopen subsets of $2^\omega$. By Lemma 4.1 in [Me], the space $\mathbb{M}$ of minimal homeomorphisms of $2^\omega$ is a Polish space. The equivalence relation of flip-conjugacy on $\mathbb{M}$ is denoted by $FCO$. The standard way to compare analytic equivalence relations on standard Borel spaces is the Borel reducibility quasi-order $\leq_B$ (see, for instance, [G]). Recall that if $X,Y$ are standard Borel spaces and $E,F$ are analytic equivalence relations on $X,Y$ respectively, then 
$$(X,E)\leq_B(Y,F)\Leftrightarrow\exists\varphi\! :\! X\!\rightarrow\! Y\mbox{ Borel with }
E\! =\! (\varphi\!\times\!\varphi )^{-1}(F).$$
Theorem 13.2 in [L] essentially shows that $FCO$ is Borel reducible to the (analytic) restriction of $\equiv_c^i$ to the set of irreflexive relations $G$ on a fixed countable dense  subset of $2^\omega$ such that $(2^\omega ,G)$ has no continuous 2-coloring.  A very recent result, in [De-GR-Ka-Kun-Kw], asserts that $FCO$ is analytic complete as a set. As a consequence, this restriction is also analytic complete. We make such statement more systematic, and partly in relation with the classes $\mathcal{C}_\kappa$ (resp., $\mathcal{H}_\kappa$) of closed graphs (resp., of graphs induced by a homeomorphism) introduced in Section 3.\medskip

 Let $\kappa\! <\!\aleph_0$, $\mathfrak{C}_\kappa$ be the set of closed graphs on $2^\omega$ having no continuous 
$\kappa$-coloring, and 
$$E^\mathfrak{C}_\kappa\! :=\!\{ (G,H)\!\in\!\mathfrak{C}_\kappa^2\mid 
(2^\omega ,G)\equiv^i_c(2^\omega ,H)\} .$$ 
Now let $\kappa\!\leq\! 3$, $\mathfrak{H}_\kappa$ be the set of homeomorphisms of 
$2^\omega$ whose induced graph has no continuous $\kappa$-coloring, and 
$E^\mathfrak{H}_\kappa\! :=\!\{ (f,g)\!\in\!\mathfrak{H}_\kappa^2\mid 
(2^\omega ,G_f)\equiv^i_c(2^\omega ,G_g)\}$. Now let $\kappa\!\leq\!\aleph_0$, $D$ be a countable dense subset of $2^\omega$, $\mathfrak{D}_\kappa$ be the set of graphs $G$ on $D$ such that $(2^\omega ,G)$ has no continuous $\kappa$-coloring, and 
$E^\mathfrak{D}_\kappa\! :=\!\{ (G,H)\!\in\!\mathfrak{D}_\kappa^2\mid (2^\omega ,G)\equiv^i_c(2^\omega ,H)\}$.

\begin{thm} \label{Cor1} The spaces $\mathfrak{C}_\kappa$, $\mathfrak{H}_\kappa$ and $\mathfrak{D}_\kappa$ are Polish, and $FCO$ is Borel reducible to the analytic equivalence relations $E^\mathfrak{C}_\kappa$, $E^\mathfrak{H}_\kappa$ and 
$E^\mathfrak{D}_\kappa$. In particular, these relations are analytic complete as sets.\end{thm}

\section{$\!\!\!\!\!\!$ Fixed points}\indent

 The set $F_f\! :=\!\{ x\!\in\!\mbox{Domain}(f)\mid f(x)\! =\! x\}$ of fixed points of $f$ is very much related to the continuous colorings of $G_f$. The next two results are essentially Proposition 7.2 and Corollary 7.3 in [L]. We recall them for the convenience of the reader.

\begin{prop} \label{fp} Let $X$ be a first countable space, and $f\! :\! X\!\rightarrow\! X$ be a partial continuous function. If $F_f$ is not an open subset of $\mbox{Domain}(f)$, then there is no continuous $\aleph_0$-coloring of $G_f$.\end{prop}

\noindent\bf Proof.\rm\ We argue by contradiction, which gives $c\! :\! X\!\rightarrow\!\aleph_0$. Let $(C_i)_{i\in\aleph_0}$ be the  partition of $X$ into clopen sets given by $C_i\! :=\! c^{-1}(\{ i\})$. As $F_f$ is not open in $\mbox{Domain}(f)$, we can find 
$x\!\in\! F_f$ and $(x_n)_{n\in\omega}\!\in\! (\mbox{Domain}(f)\!\setminus\! F_f)^\omega$ converging to $x$. Note that $f(x_n)$ is different from $x_n$, and $\big( f(x_n)\big)_{n\in\omega}$ converges to $f(x)\! =\! x$. Let $i$ with $x\!\in\! C_i$. Then we may assume that $x_n,f(x_n)\!\in\! C_i$. This implies that $\big( x_n,f(x_n)\big)\!\in\! G_f\cap C_i^2$, which is the desired contradiction.\hfill{$\square$}

\vfill\eject
 
\begin{cor} \label{corfp} Let $X$ be a 0DM space, and $f\! :\! X\!\rightarrow\! X$ be a partial continuous function with closed domain.\smallskip

\noindent (a) Exactly one of the following holds:\smallskip

(1) $F_f$ is an open subset of $\mbox{Domain}(f)$,\smallskip

(2) there is no continuous $\aleph_0$-coloring of $G_f$.\smallskip

\noindent (b) If $\mbox{Domain}(f)$ is clopen in $X$, $F_f$ is an open subset of $\mbox{Domain}(f)$, $f$ is injective and 
$1\!\leq\!\kappa\!\leq\!\aleph_0$, then $(X,G_f)$ has a continuous $\kappa$-coloring if and only if 
$(X\!\setminus\! F_f,G_f\cap (X\!\setminus\! F_f)^2)$ has a continuous $\kappa$-coloring.\end{cor}

\noindent\bf Proof.\rm\ (a) Assume that (1) holds. Note that 
$s\big(\mbox{Graph}(f_{\vert\text{Domain}(f)\setminus F_f})\big)$ and $\Delta (X)$ are disjoint and closed in $X^2$. 22.16 in [K1] gives $C\!\subseteq\! X^2$ clopen with 
${\Delta (X)\!\subseteq\! C\!\subseteq\! X^2\!\setminus\! s\big(\mbox{Graph}(f_{\vert\text{Domain}(f)\setminus F_f})\big)}$. The relation $C$ gives a continuous $\aleph_0$-coloring of $G_f$ since $X$ is zero-dimensional and second countable. So (2) does not hold.\medskip

 If $F_f$ is not an open subset of $\mbox{Domain}(f)$, then we apply Proposition \ref{fp}.\medskip
   
\noindent (b) Let $c\! :\! X\!\setminus\! F_f\!\rightarrow\!\kappa$ be a continuous $\kappa$-coloring of 
$(X\!\setminus\! F_f,G_f\cap (X\!\setminus\! F_f)^2)$. We extend $c$ to $X$ by setting $c(x)\! :=\! 0$ if $x\!\in\! F_f$. As 
$F_f$ is a clopen subset of the clopen set $\mbox{Domain}(f)$, this extension is continuous. If $\big( x,f(x)\big)\!\in\! G_f$, then $x\!\notin\! F_f$. As $f$ is injective and $F_f$ is $f$-invariant, $f[\mbox{Domain}(f)\!\setminus\! F_f]\cap F_f\! =\!\emptyset$, so that $f(x)\!\notin\! F_f$. Thus $c(x)\!\not=\! c\big( f(x)\big)$, showing that $c$ is a $\kappa$-coloring of $(X,G_f)$. Conversely, any continuous $\kappa$-coloring of $(X,G_f)$ defines a continuous $\kappa$-coloring of 
$(X\!\setminus\! F_f,G_f\cap (X\!\setminus\! F_f)^2)$ by restriction.\hfill{$\square$}\medskip

 In particular, $({}_n\Sigma ,G_\sigma )$ has no continuous $\aleph_0$-coloring if $n$ is 1 or even.

\begin{cor} \label{exfp} Let $X$ be a 0DM space, and $f$ be a homeomorphism of $X$. Then exactly one of the following holds:\smallskip

(1) $F_f$ is an open subset of $X$,\smallskip

(2) $(\mathbb{X}_1,G_{f_1})\preceq^i_c(X,G_f)$.\end{cor}

\noindent\bf Proof.\rm\ If (1) holds, then Corollary \ref{corfp} provides a continuous $\aleph_0$-coloring of $G_f$. If (2) also holds, then $G_{f_1}$ also has such a coloring $c$. As $F_{f_1}\! =\!\{ 0^\infty\}$ is not an open subset of $\mathbb{X}_1$, this contradicts Proposition \ref{fp}. So assume that (1) does not hold. Then we can find an injective sequence $(x_n)_{n\in\omega}$ of points of $X\!\setminus\! F_f$ converging to a point $x$ of $F_f$. Moreover, we may assume that $\{ x_m,f(x_m)\}\cap\{ x_n,f(x_n)\}\! =\!\emptyset$ if $m\!\not=\! n$. We then set 
$\varphi (0^\infty )\! :=\! x$, $\varphi (0^{2n}1^\infty )\! :=\! x_n$ and 
$\varphi (0^{2n+1}1^\infty )\! :=\! f(x_n)$, so that $\varphi$ is a witness for the fact that 
$(\mathbb{X}_1,G_{f_1})\preceq^i_c(X,G_f)$.\hfill{$\square$}\medskip

\noindent\bf Remark.\rm\ In the introduction, we mentioned the fact that $(\mathbb{X}_1,G_{f_1})$ is not conjugate to the shift of a two-sided subshift. Here's the argument. We argue by contradiction. By Proposition 3.68 in [Ku], 
$$\exists n\!\in\!\omega ~~\forall x\!\not=\! y\!\in\!\mathbb{X}_1~~\exists k\!\in\!\mathbb{Z}~~
f_1^k(x)\vert n\!\not=\! f_1^k(y)\vert n.$$
Fix $n\!\in\!\omega$, and choose $(x,y)\! :=\! (0^{2n}1^\infty ,0^{2n+1}1^\infty )$. Then 
$f_1^k(x)\vert n\! =\! f_1^k(y)\vert n\! =\! 0^n$ for each $k\!\in\!\mathbb{Z}$, which is the desired contradiction.\medskip

 We will use the following fact, which is part of Proposition 7.6 in [L].
 
\begin{prop} \label{LZ1+} $(\mathbb{X}_1,G_{f_1})$ is $\preceq^i_c$-minimal in the class of graphs on a 0DMC space having no continuous 2-coloring.\end{prop}

\section{$\!\!\!\!\!\!$ Basis made of minimal elements}\indent

 The part (a) of Theorems \ref{mainclosed}, \ref{mainhomeo} and \ref{mainsubshift} is based on compactness. The first key fact is that it is possible to keep a big chromatic number when we take infinite decreasing sequences of graphs or spaces.

\begin{lem} \label{compa} Assume that $2\!\leq\!\kappa\! <\!\aleph_0$ is a cardinal, $X$ is a 0DMC space, $(K_p)_{p\in\omega}$ is a decreasing sequence of closed subsets of $X$, 
$(G_p)_{p\in\omega}$ is a decreasing sequence of closed graphs on $X$ such that, for each $p\!\in\!\omega$, $G_p\!\subseteq\! K_p^2$ and $(K_p,G_p)$ has no continuous 
$\kappa$-coloring, $K\! :=\!\bigcap_{p\in\omega}~K_p$ and 
$G\! :=\!\bigcap_{p\in\omega}~G_p$. Then $(K,G)$ has no continuous $\kappa$-coloring.\end{lem}

\noindent\bf Proof.\rm\ We argue by contradiction, which gives a continuous coloring 
$c\! :\! K\!\rightarrow\!\kappa$ of $(K,G)$. By 7.8 and 2.8 in [K1], there is a continuous extension $\overline{c}\! :\! X\!\rightarrow\!\kappa$ of $c$. We set, for $p\!\in\!\omega$ and $\varepsilon\!\in\!\kappa$, 
$K^p_\varepsilon\! :=\! K_p\cap\overline{c}^{-1}(\{\varepsilon\} )$, so that 
$(K^p_\varepsilon )_{\varepsilon\in\kappa}$ is a partition of $K_p$ into clopen subsets. By assumption, $\overline{c}_{\vert K_p}$ is not a $\kappa$-coloring of 
$(K_p,G_p)$, which gives $\varepsilon_p\!\in\!\kappa$ and 
$(x_p,y_p)\!\in\! G_p\cap (K^p_{\varepsilon_p})^2$. We may assume that 
$\varepsilon\! :=\!\varepsilon_p$ does not depend on $p$. By compactness of $X^2$, we may assume that $\big( (x_p,y_p)\big)_{p\in\omega}$ converges to some $(x,y)\!\in\! K^2$. It remains to note that $(x,y)\!\in\! G\cap\big( K\cap c^{-1}(\varepsilon )\big)^2$, which is the desired  contradiction.\hfill{$\square$}\medskip

 We are now ready to prove Theorem \ref{mainclosed} (a). Let $\mathcal{C}_\kappa$ be the class of closed graphs on a 0DMC space having no continuous $\kappa$-coloring. When we write $\mbox{sup}_{p\in\omega}~\lambda_p$, we always assume that 
$(\lambda_p)_{p\in\omega}$ is a strictly increasing sequence of ordinals.\medskip

\noindent\bf Proof of Theorem \ref{mainclosed} (a).\rm\ Note that $(X,G)\!\in\!\mathcal{C}_0$ exactly when 
$X\!\not=\!\emptyset$, so that the singleton in (b) is convenient. Note then that $(X,G)\!\in\!\mathcal{C}_1$ exactly when 
$G\!\not=\!\emptyset$, so that the singleton in (b) is convenient.\medskip

 So we may assume that $\kappa\!\geq\! 2$. We argue by contradiction, which gives $(X,G)\!\in\!\mathcal{C}_\kappa$ such that, for each $(X',G')\!\in\!\mathcal{C}_\kappa$ with 
$(X',G')\preceq^i_c(X,G)$, there is $(X'',G'')\!\in\!\mathcal{C}_\kappa$ with the property that $(X'',G'')\prec^i_c(X',G')$.\medskip

\noindent\bf Claim.\it\ For each $(X',G')\!\in\!\mathcal{C}_\kappa$ with $(X',G')\preceq^i_c(X,G)$, there is 
$(X'',G'')\!\in\!\mathcal{C}_\kappa$ such that $X''\!\subseteq\! X'$ and $G''\!\subsetneqq\! G'$.\rm\medskip

 Indeed, let $\tilde X$ be the projection $\mbox{proj}[G']$ of $G'$, and $\tilde G\! :=\! G'$. Note that 
$\tilde X\! =\!\mbox{proj}[\tilde G]$, $(\tilde X,\tilde G)\!\in\!\mathcal{C}_\kappa$ and 
$(\tilde X,\tilde G)\preceq^i_c(X',G')\preceq^i_c(X,G)$. This gives $(\tilde {\tilde X},\tilde {\tilde G})\!\in\!\mathcal{C}_\kappa$ with the property that $(\tilde {\tilde X},\tilde {\tilde G})\prec^i_c(\tilde X,\tilde G)$. Let $h$ be a witness for the fact that 
$(\tilde {\tilde X},\tilde {\tilde G})\preceq^i_c(\tilde X,\tilde G)$. We put $X''\! :=\! h[\tilde {\tilde X}]$ and 
$G''\! :=\! (h\!\times\! h)[\tilde {\tilde G}]$. Note that 
$(\tilde {\tilde X},\tilde {\tilde G})\preceq^i_c(X'',G'')\preceq^i_c(\tilde {\tilde X},\tilde {\tilde G})$ with witnesses $h$, $h^{-1}$ respectively. In particular, $(X'',G'')$ is in $\mathcal{C}_\kappa$, $X''\!\subseteq\!\tilde X\! =\!\mbox{proj}[G']\!\subseteq\! X'$, and 
$G''\!\subseteq\!\tilde G\! =\! G'$. If $G''\! =\! G'$, then $X''\! =\!\tilde X$ and 
$$(\tilde X ,\tilde G)\! =\! (X'',G'')\preceq^i_c(\tilde {\tilde X},\tilde {\tilde G})\prec^i_c(\tilde X ,\tilde G)\mbox{,}$$ 
which cannot be.\hfill{$\diamond$}\medskip

 We inductively construct a $\subseteq$-decreasing sequence $(X_\xi )_{\xi <\aleph_1}$ and a strictly $\subseteq$-decreasing sequence $(G_\xi )_{\xi <\aleph_1}$ such that 
$(X_0,G_0)\! =\! (X,G)$ and $(X_\xi ,G_\xi )\!\in\!\mathcal{C}_\kappa$, which will contradict the fact that $G$ is a 0DMC space. If $(X_\xi ,G_\xi )$ is constructed, then the claim gives $(X_{\xi +1},G_{\xi +1})\!\in\!\mathcal{C}_\kappa$ such that 
$X_{\xi +1}\!\subseteq\! X_\xi$ and $G_{\xi +1}\!\subsetneqq\! G_\xi$. If 
$\lambda\! =\!\mbox{sup}_{p\in\omega}~\lambda_p$ is a limit ordinal, then Lemma \ref{compa} applied to $X$, $(X_{\lambda_p})_{p\in\omega}$ and 
$(G_{\lambda_p})_{p\in\omega}$ implies, setting 
$X_\lambda\! :=\!\bigcap_{p\in\omega}~X_{\lambda_p}$ and 
$G_\lambda\! :=\!\bigcap_{p\in\omega}~G_{\lambda_p}$, that 
$(X_\lambda ,G_\lambda )\!\in\!\mathcal{C}_\kappa$. As 
$G_\lambda\!\subsetneqq\! G_{\lambda_p}$ for each $p\!\in\!\omega$, we are done.
\hfill{$\square$}

\vfill\eject\medskip

 We now study the graphs induced by a homeomorphism. Things become more complex since\medskip
 
\noindent - fixed poins can exist; when they cannot be avoided, the induced graph is not closed,\smallskip

\noindent - the intersection of such a graph $G_f$ with a closed square $C^2$ is not necessarily of the form $G_g$; it is of this form if $C$ is $f$-invariant.\medskip

  The next lemma is a first step towards invariance. It is about the preservation of the size of orbits with at least three points under $\preceq^i_c$.

\begin{lem} \label{orbthree} Let $X$ be a topological space, $f$ be a homeomorphism of $X$, $Y,g$ having the corresponding properties and satisfying 
$(X,G_f)\preceq^i_c(Y,G_g)$ with $h$ as a witness, and $x\!\in\! X$ with 
$\vert\mbox{Orb}_f(x)\vert\!\geq\! 3$. Then $h[\mbox{Orb}_f(x)]\! =\!\mbox{Orb}_g\big( h(x)\big)$, and either 
$h\!\circ\! f\! =\! g\!\circ\! h$ on $\mbox{Orb}_f(x)$, or $h\!\circ\! f\! =\! g^{-1}\!\circ\! h$ on $\mbox{Orb}_f(x)$.\end{lem}

\noindent\bf Proof.\rm\ Let $O\! :=\!\mbox{Orb}_f(x)$. As $f_{\vert O}$ is fixed point free, $\big( x,f(x)\big)\!\in\! G_f$. Thus 
$\Big( h(x),h\big( f(x)\big)\Big)\!\in\! G_g$, showing that $h\big( f(x)\big)\! =\! g^{\pm 1}\big( h(x)\big)$. In particular, 
$h[O]\!\subseteq\!\mbox{Orb}_g\big( h(x)\big)$. We set 
$$P\! :=\!\{ z\!\in\! O\mid h\big( f(z)\big)\! =\! g\big( h(z)\big)\}$$ 
and $M\! :=\!\{ z\!\in\! O\mid h\big( f(z)\big)\! =\! g^{-1}\big( h(z)\big)\}$. As $\vert O\vert\!\geq\! 3$, 
$\vert\mbox{Orb}_g\big( h(x)\big)\vert\!\geq\! 3$ by injectivity of $h$, and $P$ and $M$ are disjoint closed subsets of 
$O\! =\! P\cup M$. If $z\!\in\! P$, then $f(z)\!\in\! P$ since otherwise $f(z)\!\in\! M$, 
$h\big( f^2(z)\big)\! =\! g^{-1}\Big( h\big( f(z)\big)\Big)\! =\! h(z)$, $f^2(z)\! =\! z$ by injectivity of $h$, which contradicts the fact that $\vert O\vert\!\geq\! 3$. Thus $O\! =\! P$ or $O\! =\! M$. In particular, either 
$h\big( f^i(x)\big)\! =\! g^i\big( h(x)\big)$ for each $i\!\in\!\mathbb{Z}$, or 
$h\big( f^i(x)\big)\! =\! g^{-i}\big( h(x)\big)$ for each $i\!\in\!\mathbb{Z}$. In both cases, we get 
$h[O]\! =\!\mbox{Orb}_g\big( h(x)\big)$.\hfill{$\square$}\medskip

 The next lemma is about the preservation of the size of orbits of size two under $\preceq^i_c$ (an orbit of size two could be sent into a bigger orbit since we consider symmetrizations). Let $\mathcal{H}_\kappa$ be the class of graphs, induced by a homeomorphism of a 0DMC space, having no continuous $\kappa$-coloring. We set, for $(X,G_f)\!\in\!\mathcal{H}_\kappa$, 
$F^{X,f}_2\! :=\!\{ x\!\in\! X\mid f^2(x)\! =\! x\}$.

\begin{lem} \label{orbtwo} Let $(X,G_f),(Y,G_g)\!\in\!\mathcal{H}_\kappa$ such that  
$(X,G_f)\preceq^i_c(Y,G_g)$ with $h$ as a witness and $F_2^{X,f}$ is nowhere dense in $X$, and $x\!\in\! X$ with 
$\vert\mbox{Orb}_f(x)\vert\! =\! 2$. Then $h[\mbox{Orb}_f(x)]$ is a $g$-orbit of size two.\end{lem}

\noindent\bf Proof.\rm\ As $F_2^{X,f}$ is nowhere dense in $X$, we can find a sequence $(x_n)_{n\in\omega}$ of points of 
$X\!\setminus\! F_2^{X,f}$ converging to $x$. Note that $\vert\mbox{Orb}_f(x_n)\vert\!\geq\! 3$ for each $n\!\in\!\omega$.  We set $y_n\! :=\! h(x_n)$, so that $h[\mbox{Orb}_f(x_n)]$ is $\mbox{Orb}_g(y_n)$ by Lemma \ref{orbthree}. We set $z\! :=\! f(x)$, so that $(x,z)\!\in\! G_f$ since $\vert\mbox{Orb}_f(x)\vert\! =\! 2$. Thus $\big( h(x),h(z)\big)\!\in\! G_g$, which gives 
$\theta\!\in\!\{ -1,1\}$ with $h(z)\! =\! g^\theta\big( h(x)\big)$. Similarly, there is, for each $n\!\in\!\omega$, $\theta_n\!\in\!\{ -1,1\}$ with $h\big( f(x_n)\big)\! =\! g^{\theta_n}(y_n)$, and we may assume that $\theta_n\! =\!\theta_0$ for each $n\!\in\!\omega$. Thus $h(z)\! =\! g^{\theta_0}\big( h(x)\big)$. So we are done if $\theta\!\not=\!\theta_0$ since 
$h(z)\!\in\!\mbox{Orb}_g\big( h(x)\big)\!\setminus\!\{ h(x)\}$. So we may assume that $\theta\! =\!\theta_0$. As 
$\vert\mbox{Orb}_f(x_n)\vert\!\geq\! 3$, $f^{-1}(x_n)\!\not=\! f(x_n)$ and 
$h\big( f^{-1}(x_n)\big)\!\not=\! h\big( f(x_n)\big)$. This implies that $h\big( f^{-1}(x_n)\big)\! =\! g^{-\theta_0}(y_n)$. Thus $h\big( f^{-1}(x)\big)\! =\! g^{-\theta_0}\big( h(x)\big)$. As $\vert\mbox{Orb}_f(x)\vert\! =\! 2$, $f^{-1}(x)\! =\! f(x)\! =\! z$, so that $g^{\theta_0}\big( h(x)\big)\! =\! g^\theta\big( h(x)\big)\! =\! h(z)\! =\! g^{-\theta_0}\big( h(x)\big)$.\hfill{$\square$}\medskip

 In the next proof, we also have to deal with orbits of size one. We are now ready to prove Theorem \ref{mainhomeo} (a).\medskip

\noindent\bf Proof of Theorem \ref{mainhomeo} (a).\rm\ As in the proof of Theorem \ref{mainclosed} (a), we may assume that 
$\kappa\!\geq\! 2$. We argue by contradiction, which gives $(X,G_f)\!\in\!\mathcal{H}_\kappa$ such that, for each 
$(X',G_{f'})\!\in\!\mathcal{H}_\kappa$ with the property that $(X',G_{f'})\preceq^i_c(X,G_f)$, there is 
$(X'',G_{f''})\!\in\!\mathcal{H}_\kappa$ with $(X'',G_{f''})\prec^i_c(X',G_{f'})$.\medskip

 If $F_f$ is not an open subset of $X$, then $(\mathbb{X}_1,G_{f_1})\preceq^i_c(X,G_f)$ by Corollary \ref{exfp}. By Corollaries \ref{exfp} and \ref{corfp}, $(\mathbb{X}_1,G_{f_1})\!\in\!\mathcal{H}_\kappa$. Our assumption gives 
$(X'',G_{f''})\!\in\!\mathcal{H}_\kappa$ with the property that $(X'',G_{f''})$ is strictly $\preceq^i_c$-below 
$(\mathbb{X}_1,G_{f_1})$, which contradicts Proposition \ref{LZ1+}. This shows that $F_f$ is an open subset of $X$. Corollary \ref{corfp} then shows that we may assume that $f$ is fixed point free. In particular, there is a $\aleph_0$-coloring of $(X,G_f)$, by Corollary \ref{corfp}, so that $\kappa\! <\!\aleph_0$.\medskip

\noindent\bf Claim.\it\ For each $(X',G_{f'})\!\in\!\mathcal{H}_\kappa$ with $(X',G_{f'})\preceq^i_c(X,G_f)$, there is 
$(X'',G_{f''})\!\in\!\mathcal{H}_\kappa$ such that $f''$ is fixed point free, $(X'',G_{f''})\preceq^i_c(X',G_{f'})$ and 
$F^{X'',f''}_2$ is nowhere dense in $X''$.\rm\medskip

 Indeed, we argue by contradiction, which gives $(X',G_{f'})\!\in\!\mathcal{H}_\kappa$. As $(X',G_{f'})\preceq^i_c(X,G_f)$, there is also a $\aleph_0$-coloring of $(X',G_{f'})$. Corollary \ref{corfp} then shows that we may assume that $f'$ is fixed point free. We inductively construct a strictly $\subseteq$-decreasing sequence $(X_\xi )_{\xi <\aleph_1}$ such that 
$X_0\! =\! X'$, $X_\xi$ is $f'$-invariant and $(X_\xi ,G_{f'}\cap X_\xi^2)\!\in\!\mathcal{H}_\kappa$, which will contradict the fact that $X'$ is a 0DMC space. Assume that $X_\xi$ is constructed. Note that $F^{X_\xi ,f'_{\vert X_\xi}}_2$ is closed and not nowhere dense in $X_\xi$. This gives a nonempty clopen subset $C$ of $X_\xi$ with the property that 
$C\!\subseteq\! F^{X_\xi,f'_{\vert X_\xi}}_2$. Note that the set $U\! :=\! C\cup f'[C]$ is a nonempty clopen $f'$-invariant subset of $X_\xi$ contained in $F^{X_\xi ,f'_{\vert X_\xi}}_2$. In particular, $U$ is a ODM separable space and 
$f'_{\vert U}$ is a fixed point free continuous involution. Proposition 7.5 in [L] provides a continuous 2-coloring of 
$(U,G_{f'_{\vert U}})$. All this implies that $X_{\xi +1}\! :=\! X_\xi\!\setminus\! U\!\subsetneqq\! X_\xi$, $X_{\xi +1}$ is $f'$-invariant and $(X_{\xi +1},G_{f'}\cap X_{\xi +1}^2)\!\in\!\mathcal{H}_\kappa$. If 
$\lambda\! =\!\mbox{sup}_{p\in\omega}~\lambda_p$ is a limit ordinal, then Lemma \ref{compa} applied to $X'$, $G_{\lambda_p}\! :=\! G_{f'}\cap X_{\lambda_p}^2$ and 
$(X_{\lambda_p})_{p\in\omega}$ implies, setting 
$X_\lambda\! :=\!\bigcap_{p\in\omega}~X_{\lambda_p}$, that $X_\lambda$ is $f'$-invariant and $(X_\lambda ,G_{f'}\cap X_\lambda^2)\!\in\!\mathcal{H}_\kappa$. As 
$X_\lambda\!\subsetneqq\! X_{\lambda_p}$ for each $p\!\in\!\omega$, we are done.
\hfill{$\diamond$}\medskip
 
 We inductively construct a strictly $\subseteq$-decreasing sequence $(X_\xi )_{\xi <\aleph_1}$ such that 
$X_0\! =\! X$, $X_\xi$ is $f$-invariant and $(X_\xi ,G_f\cap X_\xi^2)\!\in\!\mathcal{H}_\kappa$, which will contradict the fact that $X$ is a 0DMC space. Assume that $X_\xi$ is constructed. Our assumption gives 
$(X',G_{f'})\!\in\!\mathcal{H}_\kappa$ with the property that $(X',G_{f'})$ is stricly $\preceq^i_c$-below 
$(X_\xi ,G_f\cap X_\xi^2)$. The claim gives $(X'',G_{f''})\!\in\!\mathcal{H}_\kappa$ such that $f''$ is fixed point free, 
$(X'',G_{f''})\preceq^i_c(X',G_{f'})$ and $F^{X'',f''}_2$ is nowhere dense in $X''$. In particular, $(X'',G_{f''})$ is strictly 
$\preceq^i_c$-below $(X_\xi ,G_f\cap X_\xi^2)$. Let $h$ be a witness for the fact that 
$(X'',G_{f''})\preceq^i_c(X_\xi ,G_f\cap X_\xi^2)$. The fact that $f''$ is fixed point free and Lemmas \ref{orbthree}, \ref{orbtwo}  imply that $X_{\xi +1}\! :=\! h[X'']\!\subseteq\! X_\xi$ is $f$-invariant. Moreover, 
$(h\!\times\! h)[G_{f''}]\!\subseteq\! G_f\cap X_{\xi +1}^2$ and 
$(X'',G_{f''})\preceq^i_c(X_{\xi +1},(h\!\times\! h)[G_{f''}])$ with $h$ as a witness, so that 
$(X_{\xi +1},G_f\cap X_{\xi +1}^2)\!\in\!\mathcal{H}_\kappa$. If $(y_0,y_1)\!\in\! G_f\cap X_{\xi +1}^2$, then let 
$x_0,x_1\!\in\! X''$ with $y_\varepsilon\! =\! h(x_\varepsilon )$. Note that 
$\big( x_0,{f''}^\theta (x_0)\big)\!\in\! G_{f''}$ for each $\theta\!\in\!\{ -1,1\}$ since $f''$ is fixed point free, which implies that 
${\Big( h(x_0),h\big( {f''}^\theta (x_0)\big)\Big)\!\in\! G_f}$. If 
$\vert\mbox{Orb}_{f''}(x_0)\vert\!\geq\! 3$, then 
$h\big( {f''}(x_0)\big)\!\not=\! h\big( {f''}^{-1}(x_0)\big)$ is of the form 
$f^\eta\big( h(x_0)\big)$ for some $\eta\!\in\!\{ -1,1\}$. This gives 
$\eta_0,\theta_0\!\in\!\{ -1,1\}$ with  
${y_1\! =\! f^{\eta_0}\big( h(x_0)\big)\! =\! h\big( {f''}^{\theta_0}(x_0)\big)}$. Thus 
$x_1\! =\! {f''}^{\theta_0}(x_0)$, $(x_0,x_1)\!\in\! G_{f''}$, and 
$(y_0,y_1)\!\in\! (h\!\times\! h)[G_{f''}]$. If 
$\vert\mbox{Orb}_{f''}(x_0)\vert\! <\! 3$, then $\vert\mbox{Orb}_{f''}(x_0)\vert\! =\! 2$ since $f''$ is fixed point free, and 
$h[\mbox{Orb}_{f''}(x_0)]$ is an $f$-orbit of size two by Lemma \ref{orbtwo}. The conclusion is the same, with 
$\eta_0\! =\!\theta_0\! =\! 1$. So we proved that $(h\!\times\! h)[G_{f''}]\! =\! G_f\cap X_{\xi +1}^2$ in any case. Now note that $(X_{\xi +1},G_f\cap X_{\xi +1}^2)\! =\! (X_{\xi +1},(h\!\times\! h)[G_{f''}])\preceq^i_c(X'',G_{f''})$ with $h^{-1}$ as a witness, so that 
$(X_{\xi +1},G_f\cap X_{\xi +1}^2)\prec^i_c(X_\xi ,G_f\cap X_\xi^2)$, proving that 
$X_{\xi +1}\!\subsetneqq\! X_\xi$. If $\lambda\! =\!\mbox{sup}_{p\in\omega}~\lambda_p$ is a limit ordinal, then Lemma \ref{compa} applied to $X$, 
$G_{\lambda_p}\! :=\! G_f\cap X_{\lambda_p}^2$ and $(X_{\lambda_p})_{p\in\omega}$ implies, setting $X_\lambda\! :=\!\bigcap_{p\in\omega}~X_{\lambda_p}$, that $X_\lambda$ is $f$-invariant and $(X_\lambda ,G_f\cap X_\lambda^2)\!\in\!\mathcal{H}_\kappa$. As 
$X_\lambda\!\subsetneqq\! X_{\lambda_p}$ for each $p\!\in\!\omega$, we are done.
\hfill{$\square$}

\vfill\eject\medskip

 We now study subshifts. We have to find another solution when fixed points cannot be avoided since $(\mathbb{X}_1,G_{f_1})$ is not conjugate to the shift of a two-sided subshift. If $x\!\in\! A^\mathbb{Z}$ and $j\!\leq\! k$ are integers, then we define $x_{[j,k]}\!\in\! A^{k-j+1}$ by $x_{[j,k]}\! :=\!\big( x(j),\cdots ,x(k)\big)$.
  
\begin{lem} \label{fixshif} Let $\Sigma\!\subseteq\! A^\mathbb{Z}$ be a two-sided subshift, $l\!\in\!\omega$, 
$a_0,\cdots\! ,a_l\!\in\! A$, and $(x_n)_{n\in\omega}$ be an injective sequence of points of $\Sigma$ converging to 
$(a_0\!\cdots\! a_l)^\mathbb{Z}$. Then we can find $s\!\in\! A^{l+1}\!\setminus\{ (a_0\!\cdots\! a_l)\}$ and 
$\gamma\!\in\! A^\omega$ with $(a_0\!\cdots\! a_l)^{-\infty}\!\cdot\! s\gamma\!\in\!\Sigma$ or 
$\gamma^{-1}s\!\cdot\! (a_0\!\cdots\! a_l)^\infty\!\in\!\Sigma$.\end{lem}

\noindent\bf Proof.\rm\ We may assume, for example, that 
${x_n}_{[-k_n(l+1),k_n(l+1)-1]}\! =\! (a_0\!\cdots\! a_l)^{2k_n}$, 
$${x_n}_{[k_n(l+1),(k_n+1)(l+1)-1]}$$ 
is a constant $s\!\not=\! (a_0\!\cdots\! a_l)$, and $k_n\!\rightarrow\!\infty$. By compactness, we may assume that the sequence 
$({x_n}_{[k_n(l+1),\infty )})_{n\in\omega}$ converges to some $s\gamma$ in $A^\omega$. Note that 
$(a_0\!\cdots\! a_l)^{-\infty}\!\cdot\! s\gamma\!\in\!\Sigma$.\hfill{$\square$}\medskip

 We also need a version of Lemma \ref{compa} for subshifts. Let $\mathcal{S}_\kappa$ be the class of graphs, induced by the shift of a two-sided subshift, having no continuous 
$\kappa$-coloring.

\begin{lem} \label{compa1shift} Let $(\Sigma_p)_{p\in\omega}$ be a decreasing sequence of two-sided subshifts such that, for each $p\!\in\!\omega$, $(\Sigma_p,G_\sigma )\!\in\!\mathcal{S}_{\aleph_0}$, and 
$\Sigma\! :=\!\bigcap_{p\in\omega}~\Sigma_p$. Then $(\Sigma ,G_\sigma )\!\in\!\mathcal{S}_{\aleph_0}$.\end{lem}

\noindent\bf Proof.\rm\ Assume that $\Sigma_0\!\subseteq\! A^\mathbb{Z}$. Note that $\sigma_{\vert\Sigma_0}$ has finitely many fixed points since these fixed points are of the form $a^\mathbb{Z}$ for $a\!\in\! A$ and $A$ is finite. As 
$(\Sigma_p,G_\sigma )\!\in\!\mathcal{S}_{\aleph_0}$, we can find $a\!\in\! A$ such that, for any $p\!\in\!\omega$, 
$a^\mathbb{Z}\!\in\!\Sigma_p$ is not isolated in $\Sigma_p$. Lemma \ref{fixshif} provides $b\!\in\! A\!\setminus\!\{ a\}$ and, for example and for each $p\!\in\!\omega$, $\gamma_p\!\in\! A^\omega$ such that $a^{-\infty}\!\cdot\! b\gamma_p$ is in 
$\Sigma_p$. Extracting a subequence if necessary, we may assume that 
$(a^{-\infty}\!\cdot\! b\gamma_p)_{p\in\omega}$ converges to $a^{-\infty}\!\cdot\! b\gamma\!\in\!\Sigma_0$, by compactness. Note that $a^{-\infty}\!\cdot\! b\gamma\!\in\!\Sigma$, so that $(\Sigma ,G_\sigma )\!\in\!\mathcal{S}_{\aleph_0}$.\hfill{$\square$}\medskip

 We are now ready to prove Theorem \ref{mainsubshift} (a).\medskip
 
\noindent\bf Proof of Theorem \ref{mainsubshift} (a).\rm\ As in the proof of Theorem \ref{mainclosed} (a), we may assume that 
$\kappa\!\geq\! 2$. We argue by contradiction, which gives $(\Sigma ,G_\sigma )\!\in\!\mathcal{S}_\kappa$ such that, for each 
$(\Sigma',G_\sigma )\!\in\!\mathcal{S}_\kappa$ with the property that $(\Sigma',G_\sigma )\preceq^i_c(\Sigma ,G_\sigma )$, there is $(\Sigma'',G_\sigma )\!\in\!\mathcal{S}_\kappa$ with $(\Sigma'',G_\sigma )\prec^i_c(\Sigma',G_\sigma )$.\medskip

\noindent\bf Case 1.\rm\ There is $(\Sigma',G_\sigma )\!\in\!\mathcal{S}_\kappa$ with 
$(\Sigma',G_\sigma )\preceq^i_c(\Sigma ,G_\sigma )$ such that $\sigma_{\vert\Sigma'}$ is fixed point free.\medskip
 
 We can copy the proof of Theorem \ref{mainhomeo} (a) to conclude.\medskip

\noindent\bf Case 2.\rm\ For each $(\Sigma',G_\sigma )\!\in\!\mathcal{S}_\kappa$ with 
$(\Sigma',G_\sigma )\preceq^i_c(\Sigma ,G_\sigma )$, $\sigma_{\vert\Sigma'}$ is not fixed point free.\medskip
 
\noindent\bf Claim 1.\it\ For each $(\Sigma',G_\sigma )\!\in\!\mathcal{S}_\kappa$ with 
$(\Sigma',G_\sigma )\preceq^i_c(\Sigma ,G_\sigma )$, there is $(\Sigma'',G_\sigma )\!\in\!\mathcal{S}_\kappa$ with 
$\Sigma''\!\subseteq\!\Sigma'$ and $\Sigma''$ has a dense infinite orbit, and $(\Sigma',G_\sigma )\!\in\!\mathcal{S}_{\aleph_0}$.\rm\medskip

 Indeed, assume that $\Sigma'\!\subseteq\! A^\mathbb{Z}$. Note that $\sigma_{\vert\Sigma'}$ has finitely many fixed points since these fixed points are of the form $a^\mathbb{Z}$ for $a\!\in\! A$ and $A$ is finite. Let 
$U\! :=\!\{ a^\mathbb{Z}\!\in\!\Sigma'\mid a^\mathbb{Z}\mbox{ is isolated in }\Sigma'\}$. Note that $U$ is a clopen 
$\sigma_{\vert\Sigma'}$-invariant subset of $\Sigma'$, and 
$G_{\sigma_{\vert\Sigma'\setminus U}}\! =\! G_{\sigma_{\vert\Sigma'}}$, so that 
$(\Sigma'\!\setminus\! U,G_\sigma )\!\in\!\mathcal{S}_\kappa$. As we are in Case 2, $\sigma_{\vert\Sigma'\setminus U}$ is not fixed point free, which gives $a^\mathbb{Z}\!\in\!\Sigma'$ and a sequence $(x_n)_{n\in\omega}$ of points of 
$\Sigma'\!\setminus\! F_\sigma$ converging to $a^\mathbb{Z}$. Lemma \ref{fixshif} applied to $\Sigma'$, $l\! :=\! 0$ and 
$a_0\! :=\! a$ provides $b\!\in\! A\!\setminus\!\{ a\}$ and $\gamma\!\in\! A^\omega$ such that, for example,  
$x\! :=\! a^{-\infty}\!\cdot\! b\gamma\!\in\!\Sigma'$. In particular, $a^\mathbb{Z}\!\in\!\overline{\mbox{Orb}_\sigma (x)}$. So we proved the existence of $x\!\in\!\Sigma'\!\setminus\! F_\sigma$ such that $a^\mathbb{Z}\!\in\!\overline{\mbox{Orb}_\sigma (x)}$. So $\Sigma''\! :=\!\overline{\mbox{Orb}_\sigma (x)}$ is as desired since $a^\mathbb{Z}$ is a witness for the fact that $(\Sigma'',G_\sigma )\!\in\!\mathcal{S}_{\aleph_0}\!\subseteq\!\mathcal{S}_\kappa$.\hfill{$\diamond$}

\vfill\eject\medskip
 
\noindent\bf Claim 2.\it\ There is $\Sigma'\!\subseteq\!\Sigma$ such that $(\Sigma',G_\sigma )\!\in\!\mathcal{S}_\kappa$, 
$\Sigma'$ contains a dense infinite orbit $O$ and, for each $(\Sigma'',G_\sigma )\!\in\!\mathcal{S}_\kappa$ with  
$\Sigma''\!\subseteq\!\Sigma'$, $\Sigma''\cap O$ is infinite.\rm\medskip

 Indeed, we argue by contradiction.  We inductively construct a strictly $\subseteq$-decreasing sequence 
$(\Sigma_\xi )_{\xi <\aleph_1}$ such that $\Sigma_0\! =\!\Sigma$ and $(\Sigma_\xi ,G_\sigma )\!\in\!\mathcal{S}_\kappa$, which will contradict the fact that $\Sigma$ is a 0DMC space. Assume that $\Sigma_\xi$ is constructed, which is the case for 
$\xi\! =\! 0$. Claim 1 gives $(\Sigma',G_\sigma )\!\in\!\mathcal{S}_\kappa$ with $\Sigma'\!\subseteq\!\Sigma_\xi$ and $\Sigma'$ has a dense infinite orbit $O$. Our assumption gives $(\Sigma_{\xi +1},G_\sigma )\!\in\!\mathcal{S}_\kappa$ with 
$\Sigma_{\xi +1}\!\subseteq\!\Sigma'$ and $\Sigma_{\xi +1}\cap O$ is finite. In particular, 
$\Sigma_{\xi +1}\!\subsetneqq\!\Sigma_\xi$. If $(\lambda_p)_{p\in\omega}$ is strictly increasing and 
$\lambda\! =\!\mbox{sup}_{p\in\omega}~\lambda_p$ is a limit ordinal, then we set 
$\Sigma_\lambda\! :=\!\bigcap_{p\in\omega}~\Sigma_{\lambda_p}$. By Lemma \ref{compa1shift}, 
$(\Sigma_\lambda ,G_\sigma )\!\in\!\mathcal{S}_\kappa$. As $\Sigma_\lambda\!\subsetneqq\!\Sigma_{\lambda_p}$ for each 
$p\!\in\!\omega$, we are done.\hfill{$\diamond$}\smallskip

 Let $(\Sigma_0,G_\sigma )\!\in\!\mathcal{S}_\kappa$ such that $(\Sigma_0,G_\sigma )\prec^i_c(\Sigma',G_\sigma )$. Claim 1 provides $(\Sigma'_0,G_\sigma )\!\in\!\mathcal{S}_\kappa$ with $\Sigma'_0\!\subseteq\!\Sigma_0$ and 
$\Sigma'_0$ has a dense infinite orbit $O_0$. Note that $(\Sigma'_0,G_\sigma )\preceq^i_c(\Sigma',G_\sigma )$, with $h$ as a witness, and $h[O_0]$ is an infinite orbit, by Lemma \ref{orbthree}. We set $\Sigma''\! :=\!\overline{h[O_0]}$. Then $\Sigma''\!\subseteq\!\Sigma'$, $h$ is a witness for the fact that 
$(\Sigma'_0,G_\sigma )\preceq^i_c(\Sigma'',G_\sigma )$, and thus $(\Sigma'',G_\sigma )\!\in\!\mathcal{S}_\kappa$. Let $O$ be the orbit given by Claim  2. By Claim 2, $\Sigma''\cap O$ is infinite. As $\Sigma''$ is contained in the closed set $h[\Sigma'_0]$, $h[\Sigma'_0]\cap O$ is infinite. As $\Sigma'_0$ has a dense infinite orbit, $F_2^{\Sigma'_0,\sigma}$ is nowhere dense in $\Sigma'_0$. As $\Sigma'_0$ has finitely many fixed points and a finite orbit of size at least two is sent onto an orbit of the same size by $h$ by Lemmas \ref{orbthree}, \ref{orbtwo}, there is $z_0\!\in\!\Sigma'_0$ with an infinite orbit sent into $O$. This implies that $h[\mbox{Orb}(z_0)]\! =\! O$ and there is $\eta\!\in\!\{ -1,1\}$ such that 
$h\!\circ\!\sigma\! =\!\sigma^\eta\!\circ\! h$ on $\mbox{Orb}(z_0)$, by Lemma \ref{orbthree}. In particular, the set $O$ is contained in the compact set $h[\Sigma'_0]$, showing that $h$ is onto, and thus a homeomorphism, by compactness. In particular, $\mbox{Orb}(z_0)$ is dense in $\Sigma'_0$. This implies that $h\!\circ\!\sigma\! =\!\sigma^\eta\!\circ\! h$ on 
$\Sigma'_0$. If $y\!\not=\!\sigma (y)\!\in\!\Sigma'$ and, for example, $\eta\! =\! -1$, then we set $z\! :=\!\sigma (y)$. Let 
$x\!\in\!\Sigma'_0$ with $z\! =\! h(x)$. Note that 
$\big( y,\sigma (y)\big)\! =\!\big(\sigma^{-1}(z),z\big)\! =\!\big(\sigma^{-1}\big( h(x)\big) ,h(x)\big)\! =\!
\big( h\big(\sigma (x)\big) ,h(x)\big)\!\in\!(h\!\times\! h)[G_\sigma ]$, showing that 
$(h\!\times\! h)[G_\sigma ]\! =\! G_\sigma$. Thus $(\Sigma',G_\sigma )\preceq^i_c(\Sigma'_0,G_\sigma )\preceq^i_c(\Sigma_0,G_\sigma )$, which is the desired contradiction concluding the proof.\hfill{$\square$}

\section{$\!\!\!\!\!\!$ Concrete countable basis}\indent

 We already checked the part (b) of Theorems \ref{mainclosed}, \ref{mainhomeo} and \ref{mainsubshift} when 
$\kappa\!\leq\! 1$.\medskip

\noindent\bf Proof of Theorem \ref{mainhomeo} (b) when $\kappa\!\geq\! 3$.\rm\ Let $X$ be a 0DMC space, and $f$ be a homeomorphism of $X$ such that $(X,G_f)$ has no continuous $\kappa$-coloring. If $F_f$ is an open subset of $X$, then there is a continuous $\aleph_0$-coloring of $G_f$, by Corollary \ref{corfp}. By Theorem 1.12 in [L], there is a continuous $3$-coloring of $G_f$, which contradicts the fact that $\kappa\!\geq\! 3$. Thus $F_f$ is not an open subset of $X$ and we can apply Corollary \ref{exfp}.\hfill{$\square$}\medskip

 In this section, it remains to study the part (a) of Theorems \ref{main} and \ref{mainshift}.

\subsection{$\!\!\!\!\!\!$ Some general facts about symmetric relations}\indent

 Let $X$ be a set, $R$ be a relation on $X$, $\kappa$ be a countable cardinal, and $\overline{x}\! :=\! (x_i)_{i<\kappa}$ be a sequence of elements of $X$. Recall that $\overline{x}$ is a $R$-{\bf walk} if $(x_i,x_{i+1})\!\in\! R$ whenever 
$i\! +\! 1\! <\!\kappa$. A $R$-{\bf path} is an injective $R$-walk. We say that $\overline{x}$ is a $R$-{\bf cycle} if 
$3\!\leq\!\kappa\! <\!\aleph_0$, $\overline{x}$ is a $R$-path and $(x_{\kappa -1},x_0)\!\in\! R$. A {\bf connected component} of 
$(X,R)$ is a subset $C$ of $X$ such that, for each $x\!\in\! C$,\medskip

\centerline{$C\! =\!\{ y\!\in\! X\mid\exists\overline{x}~R\mbox{-path with }1\!\leq\!\kappa\! <\!\aleph_0\mbox{, }x_0\! =\! x\mbox{ and }x_{\kappa -1}\! =\! y\} .$}\medskip

\noindent We say that $(X,R)$ is {\bf connected} if $X$ is a connected component of $(X,R)$.

\vfill\eject\medskip

 The following fact is very classical.
 
\begin{lemm} \label{classical} Let $X$ be a set, and $G$ be a symmetric relation on $X$. Exactly one of the following holds:\smallskip

(1) there is a 2-coloring of $(X,G)$,\smallskip

(2) we can find $m\!\in\!\omega$ and $(x_i)_{i\leq 2m}\!\in\! X^{2m+1}$ such that 
$(x_i,x_{i+1})\!\in\! G$ for each $i\!\leq\! 2m$ (with the convention $x_{2m+1}\! :=\! x_0$).\end{lemm}

 In particular, $(2q\! +\! 3,C_{2q+3})$ has no continuous 2-coloring.

\subsection{$\!\!\!\!\!\!$ Isolated finite 2-colorable connected components}\indent

 The next two results will allow us to remove the isolated finite 2-colorable connected components. 

\begin{lemm} \label{removeisolated} Let $X$ be a 0DMC space, $G$ be a closed graph on $X$, and\medskip

\centerline{$O\! :=\!\bigcup\{ C\!\subseteq\! X\!\setminus\! X'\mid C\mbox{ finite }(X,G)\mbox{-connected component and }
(C,G\cap C^2)\mbox{ has a 2-coloring}\} .$}\medskip
  
\noindent If $(X\!\setminus\! O,G\cap (X\!\setminus\! O)^2)$ has a continuous 2-coloring, then so does $(X,G)$.\end{lemm}

\noindent\bf Proof.\rm\ Let $c\! :\! X\!\setminus\! O\!\rightarrow\! 2$ be a continuous 2-coloring of 
$(X\!\setminus\! O,G\cap (X\!\setminus\! O)^2)$. As $O$ is an open subset of the 0DM space $X$, $X\!\setminus\! O$ is closed and there is a clopen partition $(P_\varepsilon )_{\varepsilon\in 2}$ of $X$ with 
$P_\varepsilon\!\setminus\! O\! =\! c^{-1}(\{\varepsilon\} )$ (see 22.16 in [K1]).\medskip

 Let us prove that if $x\!\notin\! O$, then there is an open neighbourhood $N_x$ of $x$ such that 
$N_x\!\subseteq\! P_{c(x)}$, and $y\!\in\! P_{1-c(x)}$ if there is $x'\!\in\! N_x$ with $(x',y)\!\in\! G$. We argue by contradiction. Let $(N_i)_{i\in\omega}$ be a decreasing basis of open neighbourhoods of $x$ contained in $P_{c(x)}$. Then for each $i$ there is $(x_i,y_i)\!\in\! G\cap (N_i\!\times\! P_{c(x)})$. Note that $(x_i)_{i\in\omega}$ converges to $x$, and we may assume that 
$(y_i)_{i\in\omega}$ converges to some $y\!\in\! P_{c(x)}$ by compactness of $X$. As 
$G$ is closed, $(x,y)\!\in\! G$. As $x\!\notin\! O$ and $O$ is a union of connected components, $y\!\notin\! O$, which contradicts the fact that $c$ is a coloring.\medskip

 We now set $N\! :=\!\bigcup\{ N_x\mid x\!\notin\! O\}$, and define $\overline{c}\! :\! N\!\rightarrow\! 2$ by 
$\overline{c}(y)\! :=\! c(x)$ if $x\!\notin\! O$ and $y\!\in\! N_x$. This definition is correct since $x,x'\!\notin\! O$ and 
$y\!\in\! N_x\cap N_{x'}$ imply that $c(x)\! =\! c(x')$. If $x,z\!\notin\! O$ and $(x',y)\!\in\! G\cap (N_x\!\times\! N_z)$, then 
$\overline{c}(x')\! =\! c(x)\!\not=\! c(z)\! =\!\overline{c}(y)$ since $y\!\in\! P_{1-c(x)}\cap P_{c(z)}$. Thus $\overline{c}$ is a 2-coloring of $(N,G\cap N^2)$. By definition, $\overline{c}$ is continuous.\medskip

 Now note that $X\!\setminus\! N$ is finite, since otherwise there is an injective sequence $(w_i)_{i\in\omega}$ of elements of $X\!\setminus\! N$, and we may assume that it converges to some $w\!\in\! X$ by compactness of $X$. As 
$X\!\setminus\! O\!\subseteq\! N$, $X\!\setminus\! N\!\subseteq\! O\!\subseteq\! X\!\setminus\! X'$, so that 
$w\!\in\! X'\!\subseteq\! N$ and $w_i\!\in\! N$ if $i$ is big enough, which is the desired contradiction.\medskip

 As $X\!\setminus\! N$ is finite, the set 
$I\! :=\!\bigcup\{ C\mid C\mbox{ appears in the definition of }O\mbox{ and }C\!\setminus\! N\!\not=\!\emptyset\}$ is finite and 
$G$-invariant. We restrict $\overline{c}$ to $N\!\setminus\! I$ and extend this restriction using any 2-coloring on each of the components of $I$ to conclude.\hfill{$\square$}

\begin{coro} \label{corremoveisolated} Let $X$ be a 0DMC space, $f$ be a fixed point free homeomorphism of $X$, and 
$$O\! :=\!\bigcup\{\mbox{Orb}(x)\mbox{ finite of even cardinality}\mid x\!\in\! X\!\setminus\! X'\} .$$  
If $(X\!\setminus\! O,G_f\cap (X\!\setminus\! O)^2)$ has a continuous 2-coloring, then so does $(X,G_f)$.\end{coro}

 Note that in this section there is no upper bound on the Cantor-Bendixson rank of $X$.

\subsection{$\!\!\!\!\!\!$ Finiteness results}

\bf Convention.\rm\ In this section, $X$ is a countable MC space, and $f$ is a homeomorphism of $X$.\medskip

 The following classical fact will be used a lot.
  
\begin{lemm} \label{lastfinite} Assume that $X$ has Cantor-Bendixson rank $\beta\! +\! 1$. Then $X^\beta$ is finite.\end{lemm}

\noindent\bf Proof.\rm\ $X^\beta$ is compact and discrete, and thus finite.\hfill{$\square$}
\medskip

 In practice, in our spaces of Cantor-Bendixson rank at most two, we will consider a partition of the Cantor-Bendixson derivative into finitely many (closed) invariant sets.

\begin{lemm} \label{clopenfinite} Let $C$ be a clopen $f$-invariant subset of $X'$, $O$ be an open subset of $X$ containing 
$C$, and $(x_n)_{n\in\omega}$ be a sequence of points of $X\!\setminus\! X'$ converging to a point of $C$. Then 
$$\{\mbox{Orb}(x_n)\mid n\!\in\!\omega\mbox{ and }\mbox{Orb}(x_n)\!\not\subseteq\! O\}$$ 
is finite.\end{lemm}

\noindent\bf Proof.\rm\ Note that $C,X'\!\setminus\! C$ are disjoint and clopen in $X'$, and thus closed in the compact space $X$. If $y\!\in\! C$, then $y,f(y),f^{-1}(y)$ are in the open set 
$V\! :=\! X\!\setminus\! (X'\!\setminus\! C)$, which gives an open neighbourhood $N_y$ of $y$ with 
$\overline{N_y}\!\subseteq\! O\cap V\cap f^{-1}(V)\cap f[V]$. The compactness of $C$ provides $F\!\subseteq\! C$ finite such that $C\!\subseteq\! N\! :=\!\bigcup_{y\in F}~N_y\!\subseteq\!\overline{N}\! =\!
\bigcup_{y\in F}~\overline{N_y}\!\subseteq\! V\cap f^{-1}(V)\cap f[V]$. In particular, 
$X'\!\setminus\! C\! =\! X\!\setminus\! V$ is contained in the open set 
$U\! :=\! X\!\setminus\!\big(\overline{N}\cup f[\overline{N}]\cup f^{-1}(\overline{N})\big)$. Note that we can find 
$n_0\!\in\!\omega$ such that $x_n\!\in\! N$ if $n\!\geq\! n_0$. We set 
$M\! :=\!\{ n\!\geq\! n_0\mid\exists m\!\in\!\omega ~~f^m(x_n)\!\in\! N\wedge f^{m+1}(x_n)\!\notin\! N\}$. If $n\!\in\! M$, then there is $m_n\!\in\!\omega$ with $f^{m_n+1}(x_n)\!\notin\! N\cup U$. As $X'\!\subseteq\! N\cup U$, 
$X\!\setminus\! (N\cup U)$ is finite. This shows that $\{\mbox{Orb}_f(x_n)\mid n\!\in\! M\}$ is finite. Moreover, 
$\{ f^m(x_n)\mid m\!\in\!\omega\}\!\subseteq\! N\!\subseteq\! O$ if $n_0\!\leq\! n\!\notin\! M$. We can argue similarly with $f^{-1}$ instead of $f$, so that $\{\mbox{Orb}_f(x_n)\mid n\!\in\!\omega\wedge\mbox{Orb}_f(x_n)\!\not\subseteq\! O\}$ is finite.
\hfill{$\square$}\medskip

\noindent\bf Notation.\rm\ If $b$ is a bijection of a set $S$, $x\!\in\! S$ and $\bf{d}\!\in\{ -,+\}$, then we set 
$$\mbox{Orb}^{\bf d}_b(x)\! :=\!\{ b^{{\bf d} i}(x)\mid i\!\in\!\omega\} .$$ 

\noindent\bf Convention.\rm\ In the rest of this section, we assume the existence of $\kappa\!\in\!\omega$ and a (finite) partition $(C_j)_{j\leq\kappa}$ of $X'$ into closed $f$-invariant sets.\medskip

 The next lemma controls the closures of the orbits and is a basic tool.

\begin{lemm} \label{twosides} Assume that $x\!\in\! X\!\setminus\! X'$ has an infinite orbit. Then we can find 
$j^-,j^+\!\leq\!\kappa$ such that, for each ${\bf d}\!\in\{ -,+\}$, 
$\overline{\mbox{Orb}^{\bf d} (x)}\!\subseteq\!\mbox{Orb}^{\bf d} (x)\cup C_{j^{\bf d}}$.\smallskip

 If moreover the $C_j$'s are orbits, then $X'$ is finite and we can find $y^-,y^+\!\in\! X'$ such that, for each 
${\bf d}\!\in\{ -,+\}$, $\overline{\mbox{Orb}^{\bf d} (x)}\! =\!\mbox{Orb}^{\bf d} (x)\cup\mbox{Orb}(y^{\bf d} )\mbox{ and }
\big( f^{{\bf d}q\vert\mbox{Orb}(y^{\bf d} )\vert}(x)\big)_{q\in\omega}\mbox{ converges to }y^{\bf d}$.\end{lemm}

\noindent\bf Proof.\rm\ We first prove the following.\medskip

\noindent\bf Claim.\it\ \label{prepatwoside} Let $S$ be a closed subset of $X'$, and 
$1\!\leq\!\kappa'\!\in\!\omega$ such that the limit points of $\{ f^{q\kappa'}(x)\mid q\!\in\!\omega\}$ are in $S$. Then it is not possible to find disjoint $f^{\kappa'}$-invariant subsets $S_0,S_1$ clopen in $S$ for which we can find, for each 
$\varepsilon\!\in\! 2$, $y_\varepsilon\!\in\! S_\varepsilon$ and $(q^\varepsilon_j)_{j\in\omega}$ such that 
$\big( f^{q^\varepsilon_j\kappa'}(x)\big)_{j\in\omega}$ converges to $y_\varepsilon$.\rm\medskip

 Indeed, we argue by contradiction. Note that $S_0,S\!\setminus\! S_0$ are disjoint and closed in $S$ and $X'$, and thus closed in the compact space $X$. If $y\!\in\! S_0$, then $y,f^{\kappa'}(y)$ are in the open set 
$O\! :=\! X\!\setminus\! (S\!\setminus\! S_0)$, which gives an open neighbourhood $N_y$ of $y$ with 
$\overline{N_y}\!\subseteq\! O\cap f^{-\kappa'}(O)$. The compactness of $S_0$ provides $F\!\subseteq\! S_0$ finite such that 
$S_0\!\subseteq\! N\! :=\!\bigcup_{y\in F}~N_y\!\subseteq\!\overline{N}\! =\!
\bigcup_{y\in F}~\overline{N_y}\!\subseteq\! O\cap f^{-\kappa'}(O)$. In particular, $S\!\setminus\! S_0\! =\! X\!\setminus\! O$ is contained in the open set $N'\! :=\! X\!\setminus\! (\overline{N}\cup f^{\kappa'}[\overline{N}])$. Note that the set 
$\{ q\!\in\!\omega\mid f^{q\kappa'}(x)\!\in\! N\wedge f^{(q+1)\kappa'}(x)\!\notin\! N\}$, and thus  
$\{ q\!\in\!\omega\mid f^{(q+1)\kappa'}(x)\!\notin\! N\cup N'\}$, are infinite. By compactness, a subsequence of these 
$f^{(q+1)\kappa'}(x)$'s has to converge to a point of $S\!\subseteq\! N\cup N'$, which is the desired contradiction.
\hfill{$\diamond$}\smallskip

 As $\mbox{Orb}(x)$ is infinite, the sequence $\big( f^n(x)\big)_{n\in\omega}$ is injective and contained in the discrete space $X\!\setminus\! X'$. The compactness of $X$ provides a strictly increasing sequence $(n_q)_{q\in\omega}$ and 
$y^+\!\in\! X'$ such that $\big( f^{n_q}(x)\big)_{q\in\omega}$ converges to $y^+$. Fix $j^+\!\leq\!\kappa$ with 
$y^+\!\in\! C_{j^+}$. The claim applied to $S\! :=\! X'\! =\!\bigcup_{j\leq\kappa}~C_j$ and $\kappa'\! :=\! 1$ shows that 
$\overline{\mbox{Orb}^+(x)}\!\subseteq\!\mbox{Orb}^+(x)\cup C_{j^+}$.\smallskip

 If the $C_j$'s are orbits, then $C_{j^+}\! =\!\mbox{Orb}(y^+)$, so that 
$\overline{\mbox{Orb}^+(x)}\!\subseteq\!\mbox{Orb}^+(x)\cup\mbox{Orb}(y^+)$. As $y^+$ is in 
$\overline{\mbox{Orb}^+(x)}$, we actually have equality since $f$ is a homeomorphism. As $X'$ is a nonempty countable MC space, there is a countable ordinal $\beta$ such that $X'$ has Cantor-Bendixson rank $\beta\! +\! 1$. Thus $(X')^\beta$ is nonempty finite by Lemma \ref{lastfinite}. If $\beta\! =\! 0$, then $X'$ is finite. If $\beta\!\geq\! 1$ and 
$z\!\in\! X'\!\setminus\! (X')^\beta$ has an infinite orbit, then this orbit is not closed, which contradicts our assumptions on the $C_j$'s. This shows that $X'$ has finite orbits and 
$\beta\! =\! 0$. In particular, we may assume that $y^+$ is a limit point of 
$\{ f^{q\kappa'}(x)\mid q\!\in\!\omega\}$, where 
$\kappa'\! :=\!\vert\mbox{Orb}(y^+)\vert$. This gives $(q_l)_{l\in\omega}$ such that 
$\big( f^{q_l\kappa'}(x)\big)_{l\in\omega}$ converges to $y^+$. The claim applied to $S\! :=\!\mbox{Orb}(y^+)$ and 
$\kappa'$ implies that $\big( f^{q\kappa'}(x)\big)_{q\in\omega}$ converges to $y^+$.\medskip

 We argue similarly with $\mbox{Orb}^-(x)$ instead of $\mbox{Orb}^+(x)$.\hfill{$\square$}\medskip

 There is no upper bound on the Cantor-Bendixson rank of $X$ in Lemma \ref{clopenfinite}, the first part of Lemma \ref{twosides}, and Lemmas \ref{twosidesfinite}, \ref{orbsquare} to come. The next two results complete Lemma \ref{twosides}.
 
\begin{lemm} \label{twosidesfinite} Let $j^-\!\not=\! j^+\!\leq\!\kappa$, $O^-,O^+$ be disjoint open subsets of $X$ such that $C_{j^-}\!\subseteq\! O^-$ and $\bigcup_{j^-\not= j\leq\kappa}~C_j\!\subseteq\! O^+$, and  
$(x_n)_{n\in\omega}$ be a sequence of points of $X\!\setminus\! X'$ such that $\mbox{Orb}(x_n)$ is infinite and 
$\mbox{Orb}^+(x_n)$ meets $O^-$ and $O^+$ for each $n$. Then the set 
$\{\mbox{Orb}(x_n)\mid n\!\in\!\omega\}$ is finite. In particular, the set 
$\{\mbox{Orb}(x)\mbox{ infinite}\mid x\!\in\! X\!\setminus\! X'\mbox{ and }\forall {\bf d}\!\in\!\{ -,+\} ~~
\overline{\mbox{Orb}^{\bf d} (x)}\!\subseteq\!\mbox{Orb}^{\bf d} (x)\cup C_{j^{\bf d}}\}$ is finite.\end{lemm}

\noindent\bf Proof.\rm\ Note that the sets $C_{j^-}$, $\bigcup_{j^-\not= j\leq\kappa}~C_j$ are disjoint and clopen in $X'$, and thus closed in the compact space $X$. We argue by contradiction, so that we may assume that 
$\big(\mbox{Orb}(x_n)\big)_{n\in\omega}$ is injective. If $n\!\in\!\omega$, then we can find, for each 
${\bf d}\!\in\!\{ -,+\}$, $(m^{n,{\bf d}}_q)_{q\in\omega}$ strictly increasing and $y^{\bf d}_n\!\in\! O^{\bf d}$ such that 
$\big( f^{m^{n,{\bf d}}_q}(x_n)\big)_{q\in\omega}$ converges to $y^{\bf d}_n$. In particular, we may assume that 
$f^{m^{n,{\bf d}}_q}(x_n)\!\in\! O^{\bf d}$. This gives $(p_n)_{n\in\omega}$ such that 
$f^{p_n}(x_n)\!\in\! X\!\setminus\! O^+$ and $f^{p_n+1}(x_n)\!\in\! O^+\!\subseteq\! X\!\setminus\! O^-$. The compactness of $X$ provides $y\!\in\! X'$ such that $\big( f^{p_n}(x_n)\big)_{n\in\omega}$ converges to 
$y\!\in\! X'\!\setminus\! O^+\! =\! C_{j^-}$. As $C_{j^-}$ is $f$-invariant, $f(y)\!\in\! C_{j^-}\!\subseteq\! O^-$. On the other hand, $f(y)\! =\!\mbox{lim}_{n\rightarrow\infty}~f^{p_n+1}(x_n)\!\notin\! O^-$, which is the desired contradiction.\medskip

 For the last assertion, assume that $\big(\mbox{Orb}(x_n)\big)_{n\in\omega}$ is a sequence of elements of our set. The compactness provides, for each $n$ and each ${\bf d}\!\in\!\{ -,+\}$, $(l^{n,{\bf d}}_q)_{q\in\omega}$ strictly increasing and 
$y^{\bf d}_n$ in $C_{j^{\bf d}}$ such that $\big( f^{{\bf d}l^{n,{\bf d}}_q}(x_n)\big)_{q\in\omega}$ converges to 
$y^{\bf d}_n$. In particular, we may assume that $f^{{\bf d}l^{n,{\bf d}}_q}(x_n)\!\in\! O^{\bf d}$, so that 
$\mbox{Orb}\big( f^{-l^{n,-}_q}(x_n)\big)\! =\!\mbox{Orb}(x_n)$ and $\mbox{Orb}^+\big( f^{-l^{n,-}_q}(x_n)\big)$ meets $O^-$ and $O^+$.\hfill{$\square$}

\vfill\eject\medskip

\noindent\bf Convention.\rm\ In the rest of this section, we assume the existence of a continuous 2-coloring $\overline{c}$ of 
$(X',G_f\cap (X')^2)$.\medskip

\noindent\bf Notation.\rm\ In order to simplify the notation, we will sometimes identify $\kappa$ with 
$\mathbb{Z}/{\kappa\mathbb{Z}}$ when $\kappa\!\geq\! 2$ is finite. For example, the parity $\mbox{par}(n)$ of 
$n\!\in\!\mathbb{Z}$ will often be viewed as an element of $\mathbb{Z}/{2\mathbb{Z}}$. We set, for $j\!\leq\!\kappa$ and 
$\varepsilon\!\in\! 2$, $C^\varepsilon_j\! :=\! C_j\cap\overline{c}^{-1}(\{\varepsilon\} )$.

\begin{lemm} \label{orbsquare} Let ${\bf d}\!\in\!\{ -,+\}$, $j\!\leq\!\kappa$, $\varepsilon\!\in\! 2$, and 
$x\!\in\! X\!\setminus\! X'$ for which there is a sequence $(m_q)_{q\in\omega}$ of natural numbers of constant parity such that $\big( f^{{\bf d} m_q}(x)\big)_{q\in\omega}$ converges to a point of $C^\varepsilon_j$. Then 
$$\left\{\!\!\!\!\!\!
\begin{array}{ll}
& \overline{\mbox{Orb}_{f^2}^{\bf d}\big( f^{\text{par}(m_0)}(x)\big)}
\!\subseteq\!\mbox{Orb}_{f^2}^{\bf d}\big( f^{\text{par}(m_0)}(x)\big)\cup C^\varepsilon_j\mbox{,}\cr\cr 
& \overline{\mbox{Orb}_{f^2}^{\bf d}\big( f^{1-\text{par}(m_0)}(x)\big)}
\!\subseteq\!\mbox{Orb}_{f^2}^{\bf d}\big( f^{1-\text{par}(m_0)}(x)\big)\cup C^{1-\varepsilon}_j.
\end{array}
\right.$$
\end{lemm}

\noindent\bf Proof.\rm\ Note that $\mbox{Orb}(x)$ is infinite since $x\!\in\! X\!\setminus\! X'$ and 
$C_j\!\subseteq\! X'$. It remains to apply Lemma \ref{twosides} to $f^2$ and 
$(C^\varepsilon_j)_{j\leq\kappa ,\varepsilon\in 2}$.\hfill{$\square$}

\subsection{$\!\!\!\!\!\!$ Some general facts about homeomorphisms}\indent 

 The next lemma provides a sufficient condition for minimality.\medskip

\noindent\bf Notation.\rm\ Let $\mathfrak{G}$ be the class of graphs induced by a homeomorphism of a MC space having no continuous 2-coloring.

\begin{lemm} \label{suffmini} Let $Y$ be a 0DMC space, $h$ be a fixed point free homeomorphism of $Y$ such that 
$(Y,G_h)$ has no continuous 2-coloring, and $S$ be a dense subset of $Y$ with the property that for any 
$V\!\subseteq\! Y$, for any graph $H$ on $V$ contained in $G_h$ such that $(V,H)$ has no continuous 2-coloring, and for any $y\!\in\! S$, $\big( y,h(y)\big)\!\in\! H$ holds. Then $(Y,G_h)$ is $\preceq^i_c$-minimal in $\mathfrak{G}$ and in the class of closed graphs on a MC space having no continuous 2-coloring.\end{lemm}

\noindent\bf Proof.\rm\ Assume that $(K,G)\!\in\!\mathfrak{G}$ and 
$(K,G)\preceq^i_c(Y,G_h)$ with $\varphi$ as a witness. Corollary 2.3 in [Kr-St] shows that $(Y,G_h)$ has a continuous 3-coloring, which implies that $(K,G)$ too. By compactness, $K$ is homeomorphic to a subspace of $Y$, so that $K$ is 0D. As 
$(K,G)\!\in\!\mathfrak{G}$, there is a homeomorphism 
$g\! :\! K\!\rightarrow\! K$ with $G\! =\! G_g$. In particular, the set $F_g$ of fixed points of $g$ is a clopen subset of $K$, and 
$(K\!\setminus\! F_g,G_g\cap (K\!\setminus\! F_g)^2)$ has no continuous 2-coloring, by Corollary \ref{corfp}.(b). This implies that we may assume that $g$ is fixed point free, so that $G$ is compact. We set ${V\! :=\!\varphi [K]}$ and 
$H\! :=\! (\varphi\!\times\!\varphi )[G]$, so that $V\!\subseteq\! Y$, 
$H\!\subseteq\! G_h$ is a compact graph on $V$, ${(K,G)\preceq^i_c(V,H)}$ with 
$\varphi$ as a witness, and 
$(V,H)\preceq^i_c(K,G)$ with $\varphi^{-1}$ as a witness by compactness. Note that 
$\big( y,h(y)\big)\!\in\! H$ if $y\!\in\! S$, by our assumption. The density of $S$ in $Y$ and the compactness of $H$ then imply that 
$\mbox{Graph}(h)\!\subseteq\! H$. As $H$ is a graph, we get $H\! =\! G_h$ and therefore $V\! =\! Y$. Thus 
$(Y,G_h)\preceq^i_c(K,G)$ and $(Y,G_h)$ is $\preceq^i_c$-minimal in $\mathfrak{G}$, and also in the class of closed graphs on a MC space having no continuous 2-coloring.\hfill{$\square$}

\begin{coro} \label{oddcyclesmini} Let $q$ be a natural number. Then $(2q\! +\! 3,C_{2q+3})$ is $\preceq^i_c$-minimal in 
$\mathfrak{G}$ and in the class of closed graphs on a MC space having no continuous 2-coloring.\end{coro}

\noindent\bf Proof.\rm\ Note that $Y\! :=\! 2q\! +\! 3$, equipped with the discrete topology, is a 0DMC space. The formula 
$h(i)\! :=\! (i\! +\! 1)\mbox{ mod }(2q\! +\! 3)$ defines a fixed point free homeomorphism $h$ of $Y$, and $C_{2q+3}\! =\! G_h$. Lemma \ref{classical} implies that $(Y,G_h)$ has no continuous 2-coloring. Any dense subset of $Y$ is equal to $Y$. If 
$V\!\subseteq\! Y$, $H$ is a graph on $V$ contained in $G_h$ such that $(V,H)$ has no continuous 2-coloring, and $y\!\in\! Y$, then $\big( y,h(y)\big)\!\in\! H$. Indeed, we argue by contradiction, and we may assume that $y\! =\! 2q\! +\! 2$, the other cases being similar. Then the formula $c(x)\! :=\!\mbox{par}(x)$ defines a continuous 2-coloring of $(V,H)$, which cannot be. It remains to apply Lemma \ref{suffmini}.\hfill{$\square$}\medskip

\noindent\bf Remark.\rm\ $(2q\! +\! 3,C_{2q+3})$ is in fact $\preceq^i_c$-minimal in the class of graphs on a Hausdorff topological space having no continuous 2-coloring.\medskip

 The following result is also without upper bound on the rank.\medskip
 
\noindent\bf Convention.\rm\ In the rest of this section, we assume that $X$ is a countable 0DMC space, $f$ is a homeomorphism of $X$ such that $(X',G_f\cap (X')^2)$ has a continuous 2-coloring $\overline{c}$, 
$\kappa\!\in\!\omega$, and $(C_j)_{j\leq\kappa}$ is a partition of $X'$ into closed $f$-invariant sets.\medskip
 
\noindent\bf Notation.\rm\ Let $F\! :=\!\{ C^\varepsilon_j\mid\varepsilon\!\in\! 2\mbox{ and }j\!\leq\!\kappa\}$. We define relations $D',E'$ on $F$ by\medskip
 
\leftline{$(C^\varepsilon_j,C^{\varepsilon'}_{j'})\!\in\! D'\Leftrightarrow (\varepsilon\!\not=\!\varepsilon'\mbox{ and }j\! =\! j')\mbox{ or }$}\smallskip

\rightline{$\exists x\!\in\! X\!\setminus\! X'~~
\exists (m_q)_{q\in\omega},(n_q)_{q\in\omega}\!\in\!\omega^\omega\mbox{ with constant parity such that}$}\smallskip

\rightline{$m_0\! +\! n_0\mbox{ is odd and }
\mbox{lim}_{q\rightarrow\infty}~f^{-m_q}(x)\!\in\! C^\varepsilon_j\mbox{ and }
\mbox{lim}_{q\rightarrow\infty}~f^{n_q}(x)\!\in\! C^{\varepsilon'}_{j'}\mbox{,}$}\medskip

\leftline{$(C^\varepsilon_j,C^{\varepsilon'}_{j'})\!\in\! E'\Leftrightarrow\exists x\!\in\! X\!\setminus\! X'~~
\exists (m_q)_{q\in\omega},(n_q)_{q\in\omega}\!\in\!\omega^\omega\mbox{ with constant parity such that}$}\smallskip

\rightline{$m_0\! +\! n_0\mbox{ is even and }
\mbox{lim}_{q\rightarrow\infty}~f^{-m_q}(x)\!\in\! C^\varepsilon_j\mbox{ and }
\mbox{lim}_{q\rightarrow\infty}~f^{n_q}(x)\!\in\! C^{\varepsilon'}_{j'}\mbox{,}$}\medskip

\noindent and set $D\! :=\! s(D')$ and $E\! :=\! s(E')$.

\begin{lemm} \label{genextcol} Assume that $f$ is fixed point free and $X\!\setminus\! X'$ contains only infinite orbits. Then, with the notation just above, if there is $g\! :\! F\!\rightarrow\! 2$ satisfying
$$\left\{\!\!\!\!\!\!\!
\begin{array}{ll}
& \forall (C^\varepsilon_j,C^{\varepsilon'}_{j'})\!\in\! D~~g(C^\varepsilon_j)\!\not=\! g(C^{\varepsilon'}_{j'})\mbox{,}\cr
& \forall (C^\varepsilon_j,C^{\varepsilon'}_{j'})\!\in\! E~~g(C^\varepsilon_j)\! =\! g(C^{\varepsilon'}_{j'})\mbox{,}
\end{array}
\right.$$
then $(X,G_f)$ has a continuous 2-coloring.\end{lemm}

\noindent\bf Proof.\rm\ We define $c\! :\! X\!\rightarrow\! 2$ as follows. If $y\!\in\! X'$, then there is a unique 
$C^\varepsilon_j\!\in\! F$ with $y\!\in\! C^\varepsilon_j$. We put $c(y)\! :=\! g(C^\varepsilon_j)$. If 
$x\!\in\! X\!\setminus\! X'$, then $\mbox{Orb}(x)$ is infinite, and Lemma \ref{twosides} provides a unique 
$j^-\!\leq\!\kappa$ such that $\overline{\mbox{Orb}^-(x)}\!\subseteq\!\mbox{Orb}^-(x)\cup C_{j^-}$. As $X\!\setminus\! X'$ is discrete and $\mbox{Orb}(x)$ is infinite, there is $(m_q)_{q\in\omega}\!\in\!\omega^\omega$ with constant parity strictly increasing such that $\big( f^{-m_q}(x)\big)_{q\in\omega}$ converges to a point of $C_{j^-}$. Replacing $m_q$ with 
$m_q\! +\! 1$ if necessary, we may assume that this limit is in $C_{j^-}^0$. Lemma \ref{orbsquare} applied to ${\bf d}\! :=\! -$, 
$j\! :=\! j^-$ and $\varepsilon\! :=\! 0$ implies that 
$\overline{\mbox{Orb}_{f^2}^-\big( f^{\text{par}(m_0)}(x)}\big)
\!\subseteq\!\mbox{Orb}_{f^2}^-\big( f^{\text{par}(m_0)}(x)\big)\cup C^0_{j^-}$ and 
$\overline{\mbox{Orb}_{f^2}^-\big( f^{1-\text{par}(m_0)}(x)}\big)
\!\subseteq\!\mbox{Orb}_{f^2}^-\big( f^{1-\text{par}(m_0)}(x)\big)\cup C^1_{j^-}$. Thus, if 
$\big( f^{-m'_q}(x)\big)_{q\in\omega}$ converges to a point of $C_{j^-}^0$ and the parity of $m'_q$ is constant, then 
$\mbox{par}(m_0)\! =\!\mbox{par}(m'_0)$. This allows us to put 
$c(x)\! :=\! g(C^0_{j^-})\! +\!\mbox{par}(m_0)$.\medskip

 If $y\!\in\! C^\varepsilon_j$, then $f(y)\!\in\! C^{1-\varepsilon}_j$ since $f$ is fixed point free, so that 
$c(y)\! =\! g(C^\varepsilon_j)\!\not=\! g(C^{1-\varepsilon}_j)\! =\! c\big( f(y)\big)$ since 
$(C^\varepsilon_j,C^{1-\varepsilon}_j)\!\in\! D$. If $x\!\in\! X\!\setminus\! X'$, then $c(x)\! =\! 
g(C^0_{j^-})\! +\!\mbox{par}(m_0)\!\not=\! g(C^0_{j^-})\! +\!\mbox{par}(m_0\! +\! 1)\! =\! c\big( f(x)\big)$. Thus $c$ is a 2-coloring of $(X,G_f)$.\medskip

 Assume that $(x_n)_{n\in\omega}\!\in\! (X\!\setminus\! X')^\omega$ converges to $y\!\in\! C^\varepsilon_j$, so that 
$c(y)\! =\! g(C^\varepsilon_j)$. Let us prove that $c(x_n)\! =\! c(y)$ if $n$ is big enough. As $\mbox{Orb}(x_n)$ is infinite, Lemma \ref{twosides} provides a unique $j_n^-\!\leq\!\kappa$ such that 
$\overline{\mbox{Orb}^-(x_n)}\!\subseteq\!\mbox{Orb}^-(x_n)\cup C_{j_n^-}$. Splitting the sequence $(x_n)_{n\in\omega}$ into finitely many subsequences if necessary, we may assume that the sequence $(j_n^-)_{n\in\omega}$ is a constant $j^-$. Replacing $x_0$ with $f(x_0)$ if necessary, we may assume that $c(x_0)\! =\! g(C^0_{j^-})$. As above, there is 
${(m^n_q)_{q\in\omega}\!\in\!\omega^\omega}$ with constant parity strictly increasing such that 
$\big( f^{-m^n_q}(x_n)\big)_{q\in\omega}$ converges to a point $y_n$ of $C_{j^-}$. If $n\! >\! 0$, then, replacing 
$m^n_q$ with $m^n_q\! +\! 1$ if necessary, we may assume that $y_n\!\in\! C_{j^-}^0$, so that 
$c(x_n)\! =\! g(C^0_{j^-})\! +\!\mbox{par}(m^n_0)$. If $n\! =\! 0$, then we choose $m^0_0$ even, and 
$y_0\!\in\! C_{j^-}^0$ since $c(x_0)\! =\! g(C^0_{j^-})$.\medskip

 Let us prove that we can find disjoint clopen subsets $O^0,O^1$ of $X$ satisfying, for each $\eta\!\in\! 2$,\medskip
 
\noindent - $C^\eta_j\!\subseteq\! O^\eta\!\subseteq\! X\!\setminus\! (X'\!\setminus\! C_j)$,\smallskip

\noindent - $O^\eta\cap f^{-1}(O^\eta )\! =\!\emptyset$.\medskip

 Note that $C^0_j,C^1_j$ are disjoint and clopen in $X'$, and thus closed in the 0DMC space $X$. By 22.16 in [K1], there is a clopen subset $C$ of $X$ with $C^0_j\!\subseteq\! C\!\subseteq\! X\!\setminus\! C^1_j$. Similarly, we can find clopen subsets $C^0,C^1$ of $X$ with the properties that $C^0_j\!\subseteq\! C^0\!\subseteq\! C\!\setminus\! (X'\!\setminus\! C_j)$ and 
$C^1_j\!\subseteq\! C^1\!\subseteq\! X\!\setminus\!\big( C\cup (X'\!\setminus\! C_j)\big)$. Note that 
$C^\eta_j\!\subseteq\! f^{-1}(C^{1-\eta})$ since $f$ is fixed point free, so that we can set 
$O^\eta\! :=\! C^\eta\cap f^{-1}(C^{1-\eta})$.\medskip

 We then put $O\! :=\! O^0\cup O^1$, so that, by Lemma \ref{clopenfinite} applied to $C_j$ and $O$, the set 
$$\{\mbox{Orb}(x_n)\mid n\!\in\!\omega\mbox{ and }\mbox{Orb}(x_n)\!\not\subseteq\! O\}$$ 
is finite. We set $I\! :=\!\{ n\!\in\!\omega\mid\mbox{Orb}(x_n)\!\not\subseteq\! O\}$. Note that we can find 
$n_0\!\in\!\omega$ such that $x_n\!\in\! O^\varepsilon$ if $n\!\geq\! n_0$.\medskip

 We first prove that $c(x_n)\! =\! c(y)$ if $n\!\notin\! I$ is big enough. If $n\!\notin\! I$, then 
$\mbox{Orb}(x_n)\!\subseteq\! O$, $f^{-m^n_q}(x_n)$ is in the clopen set $O$ and thus 
$y_n\!\in\! C^0_{j^-}\cap O\!\subseteq\! X'\!\setminus\! (X'\!\setminus\! C_j)\! =\! C_j$. This implies that $j^-\! =\! j$ and 
$c(x_n)\! =\! g(C^0_j)\! +\!\mbox{par}(m^n_0)$. If $q\!\geq\! q_0$ is big enough, then 
$f^{-m^n_q}(x_n)$ is in $O^0$. As $\mbox{Orb}(x_n)$ is contained in $O$ and $O^0\cap f^{-1}(O^0)\! =\!\emptyset$, 
$f^{1-m^n_q}(x_n)\!\in\! O\!\setminus\! O^0\! =\! O^1$. As ${\mbox{Orb}(x_n)\!\subseteq\! O}$ and 
${O^1\cap f^{-1}(O^1)\! =\!\emptyset}$, $f^{2-m^n_q}(x_n)\!\in\! O\!\setminus\! O^1\! =\! O^0$. Inductively, as 
$\big(\mbox{par}(m^n_q)\big)_{q\in\omega}$ is constant, ${f^{-m^n_0}(x_n)\!\in\! O^0}$. Similarly, 
$x_n\!\in\! O^{\text{par}(m^n_0)}$. This implies that $\mbox{par}(m^n_0)\! =\!\varepsilon$ if $n\!\geq\! n_0$. As 
$g(C^\varepsilon_j)$ is different from $g(C^{1-\varepsilon}_j)$, 
$c(x_n)\! =\! g(C^0_j)\! +\!\mbox{par}(m^n_0)\! =\! g(C^0_j)\! +\!\varepsilon\! =\! g(C^\varepsilon_j)\! =\! c(y)$ if $n\!\geq\! n_0$, as desired.\medskip

 As $\{\mbox{Orb}(x_n)\mid n\!\in\!\omega\mbox{ and }\mbox{Orb}(x_n)\!\not\subseteq\! O\}$ is finite, it remains to see that if 
$x_n\!\in\!\mbox{Orb}(x_0)$ and $n$ is big enough, then $c(x_n)\! =\! c(y)$. As $x_n\!\in\!\mbox{Orb}(x_0)$, we may assume that either there is $(p_n)_{n\in\omega}\!\in\!\omega^\omega$ with constant parity such that $x_n\! =\! f^{p_n}(x_0)$ for each $n$, or there is $(r_n)_{n\in\omega}\!\in\!\omega^\omega$ with constant parity such that $x_n\! =\! f^{-r_n}(x_0)$ for each $n$.\medskip

\noindent\bf Case 1.\rm\ $x_n\! =\! f^{p_n}(x_0)$ for each $n$.\medskip

\noindent\bf Subcase 1.1.\rm\ $g(C^0_{j^-})\! =\! g(C^\varepsilon_j)$.\medskip

 Note that $\big( f^{-m^0_q}(x_0)\big)_{q\in\omega}\! =\!\big( f^{-(m^0_q+p_n)}(x_n)\big)_{q\in\omega}$ converges to 
$y_0\!\in\! C^0_{j^-}$. Thus, by definition of $c$, $c(x_n)\! =\! g(C^0_{j^-})\! +\!\mbox{par}(p_n)$. Note that 
$\big( f^{p_q}(x_0)\big)_{q\in\omega}$ converges to $y\!\in\! C^\varepsilon_j$. As $g(C^0_{j^-})\! =\! g(C^\varepsilon_j)$, 
$(C^0_{j^-},C^\varepsilon_j)\!\notin\! D$, and $m^0_q\! +\! p_q$ is even, like $p_n$. Thus 
$c(x_n)\! =\! g(C^0_{j^-})\! =\! g(C^\varepsilon_j)\! =\! c(y)$, as desired.\medskip
 
\noindent\bf Subcase 1.2.\rm\ $g(C^0_{j^-})\!\not=\! g(C^\varepsilon_j)$.\medskip

 Arguing as in Subcase 1.1, as $g(C^0_{j^-})\!\not=\! g(C^\varepsilon_j)$, $(C^0_{j^-},C^\varepsilon_j)\!\notin\! E$, and thus $m^0_q\! +\! p_q$ is odd, just like $p_n$. This implies that $c(x_n)\! =\! g(C^0_{j^-})\! +\! 1\! =\! g(C^\varepsilon_j)\! =\! c(y)$, as desired.\medskip

\noindent\bf Case 2.\rm\ $x_n\! =\! f^{-r_n}(x_0)$ for each $n$.\medskip

 Note that $r_0\! =\! 0$ since $\mbox{Orb}(x_0)$ is infinite. Assume that $\varepsilon\! =\! 1$, the other case being similar. Note that $\overline{\mbox{Orb}^-(x_0)}\!\subseteq\!\mbox{Orb}^-(x_0)\cup C_{j^-}$, so that 
$y\! =\!\mbox{lim}_{n\rightarrow\infty}~x_n\!\in\! C_j\cap C_{j^-}$ and $j\! =\! j^-$. In particular, 
$c(x_n)\! =\! g(C^0_j)\! +\! 1\! -\!\mbox{par}(r_n)\! =\! g(C^0_j)\! +\! 1\! =\! g(C^1_j)\! =\! c(y)$, as desired.\medskip

 This proves the continuity of $c$.\hfill{$\square$}\medskip
  
 A key consequence of Lemma \ref{genextcol} is the following.
  
\begin{lemm} \label{both23} Assume that $f$ is fixed point free, $X\!\setminus\! X'$ contains only infinite orbits, and $(X,G_f)$ has no continuous 2-coloring. Then we can find $l\!\in\!\omega$, $C^{\varepsilon_0}_{j_0},\cdots\! ,C^{\varepsilon_l}_{j_l}$ with 
$(j_i)_{i\leq l}$ injective, a sequence $(z'_i)_{i\leq l}$ of elements of $X\!\setminus\! X'$, and, for $i\!\leq\! l$, a sequence 
$(m^i_q)_{q\in\omega}$ of even natural numbers and a sequence $(n^i_q)_{q\in\omega}$ of natural numbers with constant parity satisfying the following:\medskip

\noindent (a) if $i\! <\! l$, then $n^i_0$ is even and one of the following holds:\medskip

$(\alpha )_i~y_i^-\! :=\!\mbox{lim}_{q\rightarrow\infty}~f^{-m^i_q}(z'_i)\!\in\! C^{\varepsilon_i}_{j_i}\mbox{ and }
y_i^+\! :=\!\mbox{lim}_{q\rightarrow\infty}~f^{n^i_q}(z'_i)\!\in\! C^{\varepsilon_{i+1}}_{j_{i+1}}\mbox{,}$\smallskip

$(\beta )_i~y_i^-\! :=\!\mbox{lim}_{q\rightarrow\infty}~f^{-m^i_q}(z'_i)\!\in\! C^{\varepsilon_{i+1}}_{j_{i+1}}\mbox{ and }
y_i^+\! :=\!\mbox{lim}_{q\rightarrow\infty}~f^{n^i_q}(z'_i)\!\in\! C^{\varepsilon_i}_{j_i}\mbox{,}$\medskip

\noindent (b) $n^l_0$ is odd and one of the following holds:\medskip

$(\alpha )_l~y_l^-\! :=\!\mbox{lim}_{q\rightarrow\infty}~f^{-m^l_q}(z'_l)\!\in\! C^{\varepsilon_l}_{j_l}\mbox{ and }
y_l^+\! :=\!\mbox{lim}_{q\rightarrow\infty}~f^{n^l_q}(z'_l)\!\in\! C^{\varepsilon_0}_{j_0}\mbox{,}$\smallskip

$(\beta )_l~y_l^-\! :=\!\mbox{lim}_{q\rightarrow\infty}~f^{-m^l_q}(z'_l)\!\in\! C^{\varepsilon_0}_{j_0}\mbox{ and }
y_l^+\! :=\!\mbox{lim}_{q\rightarrow\infty}~f^{n^l_q}(z'_l)\!\in\! C^{\varepsilon_i}_{j_l}.$
\end{lemm}

\noindent\bf Proof.\rm\ By Lemma \ref{genextcol} and the notation just above its statement, it is not possible to find 
$g\! :\! F\!\rightarrow\! 2$ satisfying
$$\left\{\!\!\!\!\!\!\!
\begin{array}{ll}
& \forall (C^\varepsilon_j,C^{\varepsilon'}_{j'})\!\in\! D~~g(C^\varepsilon_j)\!\not=\! g(C^{\varepsilon'}_{j'})\mbox{,}\cr
& \forall (C^\varepsilon_j,C^{\varepsilon'}_{j'})\!\in\! E~~g(C^\varepsilon_j)\! =\! g(C^{\varepsilon'}_{j'}).
\end{array}
\right.$$
Note that if $(C^\varepsilon_j,C^{\varepsilon'}_{j'})\!\in\! E$, then 
$(C^\varepsilon_j,C^{1-\varepsilon}_j), (C^{1-\varepsilon}_j,C^{\varepsilon'}_{j'})\!\in\! D$. In particular, the first condition (on $D$) just above implies the second one (on $E$). So  there is no $g$ satisfying the first condition (on $D$). Lemma \ref{classical} provides 
$m\!\in\!\omega$ and $(\varepsilon_i,j_i)_{i\leq 2m}\!\in\! (2\!\times\!\kappa )^{2m+1}$ such that $(C^{\varepsilon_i}_{j_i},C^{\varepsilon_{i+1}}_{j_{i+1}})\!\in\! D$ for each 
$i\!\leq\! 2m$. We may assume that the sequence $\big( (\varepsilon_i,j_i)\big)_{i<n}$ is injective. We set $l\! :=\! 2m$. If $i\!\leq\! l$, then 
$(C^{\varepsilon_i}_{j_i},C^{\varepsilon_{i+1}}_{j_{i+1}})\!\in\! D\cup E$, so that 
$(\varepsilon_i\!\not=\!\varepsilon_{i+1}\mbox{ and }j_i\! =\! j_{i+1})$, or we can find $z'_i\!\in\! X\!\setminus\! X'$ and 
$(m^i_q)_{q\in\omega}$, $(n^i_q)_{q\in\omega}\!\in\!\omega^\omega$ with constant parity such that one of the following holds:\medskip

\centerline{$(\alpha )_i~y_i^-\! :=\!\mbox{lim}_{q\rightarrow\infty}~f^{-m^i_q}(z'_i)\!\in\! C^{\varepsilon_i}_{j_i}\mbox{ and }
y_i^+\! :=\!\mbox{lim}_{q\rightarrow\infty}~f^{n^i_q}(z'_i)\!\in\! C^{\varepsilon_{i+1}}_{j_{i+1}}\mbox{,}$}\medskip

\centerline{$(\beta )_i~y_i^-\! :=\!\mbox{lim}_{q\rightarrow\infty}~f^{-m^i_q}(z'_i)\!\in\! C^{\varepsilon_{i+1}}_{j_{i+1}}\mbox{ and }
y_i^+\! :=\!\mbox{lim}_{q\rightarrow\infty}~f^{n^i_q}(z'_i)\!\in\! C^{\varepsilon_i}_{j_i}.$}\medskip

\noindent Note that, changing $\varepsilon_{p+1}$ enough times if necessary, we may assume that 
$(C^{\varepsilon_i}_{j_i},C^{\varepsilon_{i+1}}_{j_{i+1}})\!\in\! E$ if $i\! <\! l$, so that 
$(C^{\varepsilon_l}_{j_l},C^{\varepsilon_0}_{j_0})\!\in\! D$. Note then that, canceling $C^{\varepsilon_l}_{j_l}$ if necessary, we may assume that the case when $(\varepsilon_i\!\not=\!\varepsilon_{i+1}\mbox{ and }j_i\! =\! j_{i+1})$ never holds. Also, replacing $(z'_i,m^i_q,n^i_q)$ with $\big( f(z'_i),m^i_q\! +\! 1,n^i_q\! -\! 1)$ if necessary, we may assume that $m^i_0$ is even if $i\!\leq\! l$.\hfill{$\square$}

\subsection{$\!\!\!\!\!\!$ General homeomorphisms}\indent

 We first study the $\Sigma_{\bf p}$'s. Recall that the space $\Sigma_{\bf p}$ is MC with Cantor-Bendixson rank two, and 
$\sigma_{\vert\Sigma_{\bf p}}$ is a homeomorphism of the $\sigma$-invariant space $\Sigma_{\bf p}$. It will be convenient to set, for 
$${\bf p}\! =\! (l,\lambda_0,\cdots ,\lambda_l,m,\varepsilon_0,\cdots ,\varepsilon_{l-1})\!\in\!\omega^{l+3}\!\times\! 2^l\mbox{,}$$ 
$y_i\! :=\! w_i^\mathbb{Z}\mbox{ if }i\!\leq\! l\mbox{, }
z_i\! :=\! w_{i+\varepsilon_i}^{-\infty}\!\cdot\! w_{i+1-\varepsilon_i}^\infty\mbox{ if }i\! <\! l\mbox{, }
z_l\! :=\! w_l^{-\infty}\!\cdot\! b_0\!\cdots\! b_{m-1}(w_0^\infty )$. It is important to note that 
$\mbox{lim}_{q\rightarrow\infty}~~\sigma^{-q\lambda_{i+\varepsilon_i}}(z_i)\! =\! y_{i+\varepsilon_i}$ and 
$$\mbox{lim}_{q\rightarrow\infty}~~\sigma^{q\lambda_{i+1-\varepsilon_i}}(z_i)\! =\! y_{i+1-\varepsilon_i}$$ 
if $i\! <\! l$. Similarly, 
$\mbox{lim}_{q\rightarrow\infty}~~\sigma^{-q\lambda_l}(z_l)\! =\! y_l
\mbox{ and lim}_{q\rightarrow\infty}~~\sigma^{q\lambda_0+m}(z_l)\! =\! y_0$.

\begin{them} \label{bas} Let $X$ be a countable MC space with Cantor-Bendixson rank at most two, and $f$ be a fixed point free homeomorphism of $X$ such that $(X,G_f)$ contains no odd cycle and has no continuous 2-coloring. Then there is 
${\bf p}\!\in\!\mathcal{P}$ such that $(\Sigma_{\bf p},G_\sigma )\preceq^i_c(X,G_f)$.\end{them}

\noindent\bf Proof.\rm\ Note that $X$ is countable, so that $X$ is 0D by 7.12 in [K1]. By Corollary \ref{corremoveisolated}, we may assume that $X\!\setminus\! X'$ contains only infinite orbits. Note that $X'$ is finite by Lemma \ref{lastfinite}, which gives 
$\kappa\!\in\!\omega$ and a partition $(C_j)_{j\leq\kappa}$ of $X'$ into orbits, which are closed and $f$-invariant sets. Note that the $C_j$'s have even cardinality, which gives a (continuous) 2-coloring $\overline{c}$ of $(X',G_f\cap (X')^2)$. Lemma \ref{both23} provides $l$, $C^{\varepsilon_0}_{j_0},\cdots\! ,C^{\varepsilon_l}_{j_l}$, $(z'_i)_{i\leq l}$, and, for $i\!\leq\! l$, 
$(m^i_q)_{q\in\omega}$ and $(n^i_q)_{q\in\omega}$. We set, for $i\!\leq\! l$, $\lambda_i\! :=\!\vert C_{j_i}\vert$, so that 
$\lambda_i\! >\! 0$ is even. Note that $f^{-1}$ is a homeomorphism and $G_{f^{-1}}\! =\! G_f$. So, replacing $f$ and $z'_l$ with $f^{-1}$ and $f(z'_l)$ respectively if necessary, we may assume that $(\alpha )_l$ holds. We set, for $i\! <\! l$, 
$$\varepsilon_i\! :=\!\left\{\!\!\!\!\!\!\!
\begin{array}{ll}
& 0\mbox{ if }(\alpha )_i\mbox{ holds,}\cr
& 1\mbox{ if }(\beta )_i\mbox{ holds.}
\end{array}
\right.$$
We also define, for $i\!\leq\! l$, $\overline{i\! +\! 1}\! :=\! i\! +\! 1\mbox{ mod }(l\! +\! 1)$ and ${\bf d}_i\!\in\!\{ -,+\}$ by 
$${\bf d}_i\! :=\!\left\{\!\!\!\!\!\!\!
\begin{array}{ll}
& -\mbox{ if }(\alpha )_i\mbox{ holds,}\cr
& +\mbox{ if }(\beta )_i\mbox{ holds,}
\end{array}
\right.$$
and we will use the conventions $--\! =\! +$ and $-+\! =\! -$. We will now show that we may assume that $m^i_q$ or $n^i_q$ is equal to $q\lambda_i$ if the limit coming from Lemma \ref{both23} is in $C_{j_i}$, except $n^l_q$ that will be $q\lambda_0\! +\! m$ with $m\! <\!\lambda_0$ odd. We will also ensure that $y_i^{-{\bf d}_i}\! =\! y_{\overline{i+1}}^{{\bf d}_{\overline{i+1}}}$. If $x$ is in $C^0_{j_i}$, then 
$\mbox{Orb}_{f^2}(x)\!\subseteq\! C^0_{j_i}$ and $\mbox{Orb}_{f^2}\big( f(x)\big)\!\subseteq\! C^1_{j_i}$, showing that $C^{\varepsilon_i}_{j_i}$ is a $f^2$-orbit. The key fact is as follows.\medskip

\noindent\bf Claim.\it\ Let $x\!\in\! X\!\setminus\! X'$, ${\bf d}\!\in\!\{ -,+\}$, $j\!\leq\!\kappa$ and $\varepsilon\!\in\! 2$ for which there is a sequence $(m_q)_{q\in\omega}$ of natural numbers of constant parity such that 
$\big( f^{{\bf d} m_q}(x)\big)_{q\in\omega}$ converges to a point of $C^\varepsilon_j$. Then there is 
$y\!\in\! C^{\varepsilon +\text{par}(m_0)}_j$ such that $\big( f^{{\bf d}q\vert C_j\vert}(x)\big)_{q\in\omega}$ converges to $y$.\rm\medskip

 Indeed, Lemma \ref{twosides} provides $y\!\in\! C_j$ such that $\big( f^{{\bf d}q\vert C_j\vert}(x)\big)_{q\in\omega}$ converges to $y$. Lemma \ref{orbsquare} implies that $y\!\in\! C^{\varepsilon +\text{par}(m_0)}_j$.\hfill{$\diamond$}\medskip

 Assume first that $l\! =\! 0$. The claim provides $y^-\!\in\! C_{j_0}^{\varepsilon_0},y^+\!\in\! C_{j_0}^{\varepsilon_0+1}$ such that $\big( f^{-q\lambda_0}(z'_0)\big)_{q\in\omega}$ converges to $y^-$ and $\big( f^{q\lambda_0}(z'_0)\big)_{q\in\omega}$ converges to $y^+$. As $C^{\varepsilon_0}_{j_0}$ is an $f^2$-orbit, there is $p\! <\!\vert C^{\varepsilon_0}_{j_0}\vert$ such that $y^-\! =\! f^{2p+1}(y^+)$. As $\big( f^{q\lambda_0+2p+1}(z'_0)\big)_{q\in\omega}$ converges to $y^-$, we are done.\medskip

 Assume now that $l\!\geq\! 1$ and $i\! =\! 0$. The claim provides\medskip
 
\noindent - $y'_0\!\in\! C_{j_0}^{\varepsilon_0}$ such that $\big( f^{{\bf d}_0q\lambda_0}(z'_0)\big)_{q\in\omega}$ converges to $y'_0$,\smallskip

\noindent - $y'_1\!\in\! C_{j_1}^{\varepsilon_1}$ such that $\big( f^{-{\bf d}_0q\lambda_1}(z'_0)\big)_{q\in\omega}$ converges to $y'_1$,\smallskip

\noindent - $y''_1\!\in\! C_{j_1}^{\varepsilon_1}$ such that $\big( f^{{\bf d}_1q\lambda_1}(z'_1)\big)_{q\in\omega}$ converges to $y''_1$.\medskip
 
 As $C^{\varepsilon_1}_{j_1}$ is a $f^2$-orbit, there is $p'\! <\!\vert C^{\varepsilon_1}_{j_1}\vert$ such that 
$y'_1\! =\! f^{2p'}(y''_1)$. As $\Big( f^{{\bf d}_1q\lambda_1}\big( f^{2p'}(z'_1)\big)\Big)_{q\in\omega}$ converges to $y'_1$, 
we are done if we replace $z'_1$ with $f^{2p'}(z'_1)$, which does not affect our convergence and parity properties. Iterating this process if necessary and arguing as in the case $l\! =\! 0$, we complete our construction. In other words, possibly changing the $z'_i$'s, we ensured the existence of $(y'_i)_{i\leq l}$ and $m\! <\!\lambda_0$ odd satisfying the following.\medskip

\noindent (a) if $i\! <\! l$, then one of the following holds:\medskip

$(\alpha )_i~y'_i\! :=\!\mbox{lim}_{q\rightarrow\infty}~f^{-q\lambda_i}(z'_i)\!\in\! C^{\varepsilon_i}_{j_i}\mbox{ and }
y'_{i+1}\! :=\!\mbox{lim}_{q\rightarrow\infty}~f^{q\lambda_{i+1}}(z'_i)\!\in\! C^{\varepsilon_{i+1}}_{j_{i+1}}\mbox{,}$\smallskip

$(\beta )_i~y'_{i+1}\! :=\!\mbox{lim}_{q\rightarrow\infty}~f^{-q\lambda_{i+1}}(z'_i)\!\in\! C^{\varepsilon_{i+1}}_{j_{i+1}}\mbox{ and }
y'_i\! :=\!\mbox{lim}_{q\rightarrow\infty}~f^{q\lambda_i}(z'_i)\!\in\! C^{\varepsilon_i}_{j_i}\mbox{,}$\medskip

\noindent (b) $y'_l\! :=\!\mbox{lim}_{q\rightarrow\infty}~f^{-q\lambda_l}(z'_l)\!\in\! C^{\varepsilon_l}_{j_l}\mbox{ and }
y'_0\! :=\!\mbox{lim}_{q\rightarrow\infty}~f^{q\lambda_0+m}(z'_l)\!\in\! C^{\varepsilon_0}_{j_0}.$\medskip

 We now completely defined 
${\bf p}\! :=\! (l,\lambda_0,\cdots\! ,\lambda_l,m,\varepsilon_0,\cdots ,\varepsilon_{l-1})\!\in\!\mathcal{P}$, and are ready to define $h\! :\!\Sigma_{\bf p}\!\rightarrow\! X$. We set $h\big(\sigma^j(y_i)\big)\! :=\! f^j(y'_i)$ if $i\!\leq\! l$ and 
$j\! <\!\lambda_i$, and $h\big(\sigma^k(z_q)\big)\! :=\! f^k(z'_q)$ if $q\!\leq\! l$ and $k\!\in\!\mathbb{Z}$ so that 
$h$ is an injective homomorphism from $(\Sigma_{\bf p},G_\sigma )$ into $(X,G_f)$. Assume that $i,q\!\leq\! l$ and 
$(\zeta_n)_{n\in\omega}$ is a sequence of points of $\mbox{Orb}(z_q)$ converging to a point $y$ of $\mbox{Orb}(y_i)$. Let $k_n\!\in\!\mathbb{Z}$ with $\zeta_n\! =\!\sigma^{k_n}(z_q)$, and $j\! <\!\lambda_i$ with $y\! =\!\sigma^j(y_i)$. Note that $(k_n)_{n\in\omega}$ tends to $\infty$ or $-\infty$. Let $i_n\!\in\!\mathbb{Z}$ and $0\!\leq\! r_n\! <\!\lambda_i$ with 
$k_n\! =\! i_n\lambda_i\! +\! r_n$. We may assume that $r_n$ is a constant $r$. Then 
$h(\zeta_n)\! =\! h\big(\sigma^{k_n}(z_q)\big)\! =\! f^{k_n}(z'_q)$, and $h(y)\! =\! h\big(\sigma^j(y_i)\big)\! =\! f^j(y'_i)$. If $q\! <\! l$ or $k_n$ tends to $-\infty$, then $\sigma^{k_n}(z_q)$ tends to $\sigma^r(y_i)$, so that $j\! =\! r$ and $h(\zeta_n)$ tends to $f^r(y'_i)\! =\! h(y)$ as desired. If $q\! =\! l$ and $k_n$ tends to $\infty$, then $\sigma^{k_n}(z_q)$ tends to $\sigma^{r-m}(y_i)$, so that $j\! =\! r\! -\! m$ and $h(\zeta_n)$ tends to $f^{r-m}(y'_i)\! =\! h(y)$ as desired. Thus $h$ is continuous.
\hfill{$\square$}

\begin{lemm} \label{minSigmap} Let ${\bf p}\!\in\!\mathcal{P}$. Then $(\Sigma_{\bf p},G_\sigma )$ is $\preceq^i_c$-minimal in $\mathfrak{G}$ and in the class of closed graphs on a MC space having no continuous 2-coloring.\end{lemm}

\noindent\bf Proof.\rm\ Note first that $\sigma_{\vert\Sigma_{\bf p}}$ is fixed point free, so that $G_\sigma$ is closed.\medskip

\noindent\bf Claim 1.\it\ $(\Sigma_{\bf p},G_\sigma )$ has no continuous 2-coloring.\rm\medskip

 Indeed, we argue by contradiction, which gives $c\! :\!\Sigma_{\bf p}\!\rightarrow\! 2$. Assume, for example, that 
$c(z_l)\! =\! 0$. As $\lambda_l$ is even, $c\big(\sigma^{-q\lambda_l}(z_l)\big)\! =\! 0$, so that $c(y_l)\! =\! 0$ by continuity. As 
$\lambda_0$ is even and $m$ is odd, $c\big(\sigma^{q\lambda_0+m}(z_l)\big)\! =\! 1$, so that $c(y_0)\! =\! 1$ by continuity. On the other hand, if $i\! <\! l$ and 
$\varepsilon_i\! =\! 0$, as $\lambda_{i+\varepsilon_i}$ is even, 
$c\big(\sigma^{-q\lambda_{i+\varepsilon_i}}(z_i)\big)\! =\! c(z_i)$, so that 
$c(y_i)\! =\! c(z_i)$ by continuity. Similarly, as $\lambda_{{i+1-\varepsilon_i}}$ is even, 
$c\big(\sigma^{q\lambda_{{i+1-\varepsilon_i}}}(z_i)\big)\! =\! c(z_i)$, so that 
$c(y_{i+1})\! =\! c(z_i)$ by continuity. This shows that $c(y_i)\! =\! c(y_{i+1})$ if $i\! <\! l$ (even if $\varepsilon_i\! =\! 1$). Thus $c(y_0)\! =\! c(y_l)$, which is the desired contradiction.\hfill{$\diamond$}\medskip

 This proves that $(\Sigma_{\bf p},G_\sigma )$ is in our classes. We now set 
$S\! :=\!\bigcup_{i\leq l}~~\mbox{Orb}_\sigma (z_i)$, so that $S$ is a dense open subset of $\Sigma_{\bf p}$, and 
$F^{\Sigma_{\bf p},\sigma}_2$ is nowhere dense in $\Sigma_{\bf p}$.\medskip
 
\noindent\bf Claim 2.\it\ Let $V\!\subseteq\! \Sigma_{\bf p}$, $H$ be a graph on $V$ contained in $G_\sigma$ such that $(V,H)$ has no continuous 2-coloring, and $x\!\in\! S$. Then $\big( x,\sigma (x)\big)\!\in\! H$.\rm\medskip

 Indeed, we argue by contradiction. Recall that the sets of the form 
$$[w]_q\! :=\!\{ y\!\in\! A_{\bf p}^\mathbb{Z}\mid\forall j\! <\!\vert w\vert ~~w(j)\! =\! y(q\! +\! j)\}\mbox{,}$$ 
where $w\!\in\! A_{\bf p}^{<\omega}$ and $q\!\in\!\mathbb{Z}$, form a basis made of clopen subsets of the space 
$A_{\bf p}^\mathbb{Z}$. Assume first that $x\! =\! z_l$. We set 
$C\! :=\! (\bigcup_{i\leq l,j<\lambda_i\text{ even}}~[a^i_j]_0\cup[b_0]_0\cup\bigcup_{j<m\text{ odd}}~[b_j]_0)\cap V$,  
so that $C$ is a clopen subset of $V$ and $H\cap\big( C^2\cup (V\!\setminus\! C)^2\big)\! =\!\emptyset$. Thus $C$ defines a continuous 2-coloring of $(V,H)$, which is the desired contradiction. If there is $k\!\in\!\mathbb{Z}$ with 
$x\! =\!\sigma^k(z_l)$, then we just replace the basic clopen sets of the form $[w]_0$ in the definition of $C$ with 
$[w]_{-k}$, the rest of the argument is the same. Similarly, if $i\! <\! l$ and $x\!\in\!\mbox{Orb}(z_i)$, then we may assume that $x\! =\! z_i$. If $\varepsilon_i\! =\! 0$, then we set 
$$C\! :=\! (\bigcup_{n\leq i,j<\lambda_n\text{ even}}~[a^n_j]_0\cup [a^i_{\lambda_i-1}a^{i+1}_0]_{-1}\cup
\bigcup_{i<n\leq l,j<\lambda_n\text{ odd}}~[a^n_j]_0\cup\bigcup_{j<m\text{ odd}}~[b_j]_0)\cap V$$ 
and conclude similarly. If $\varepsilon_i\! =\! 1$, then we set 
$$C\! :=\! (\bigcup_{n\leq i,j<\lambda_n\text{ odd}}~[a^n_j]_0\cup [a^{i+1}_{\lambda_{i+1}-1}a^i_0]_{-1}\cup
\bigcup_{i<n\leq l,j<\lambda_n\text{ even}}~[a^n_j]_0\cup\bigcup_{j<m\text{ even}}~[b_j]_0)\cap V$$ 
and conclude similarly.\hfill{$\diamond$}\medskip

 It remains to apply our claims and Lemma \ref{suffmini}.\hfill{$\square$}\medskip

\noindent\bf Proof of Theorem \ref{main} (a).\rm\ Let $X$ be a countable MC space with Cantor-Bendixson rank at most two, and $f$ be a homeomorphism of $X$ such that $(X,G_f)$ has no continuous 2-coloring. As $X$ is countable, $X$ is 0D by 7.12 in [K1]. By Corollary \ref{oddcyclesmini}, we may assume that $(X,G_f)$ contains no odd cycle. If $F_f$ is not open in $X$, then $(\mathbb{X}_1,G_{f_1})\preceq^i_c(X,G_f)$ by Corollary \ref{exfp}. So we may assume that $F_f$ is open in $X$. This implies that $(X\!\setminus\! F_f,G_f\cap (X\!\setminus\! F_f)^2)$ has no continuous 2-coloring, by Corollary \ref{corfp}. Thus we may assume that $f$ is fixed point free. It remains to apply Theorem \ref{bas} to get the basis result.\medskip

 Assume that ${\bf p},{\bf p}'\!\in\!\mathcal{P}$ and 
$(\Sigma_{\bf p},G_\sigma )\preceq^i_c(\Sigma_{{\bf p}'},G_\sigma )$ with $h$ as a witness. As $\sigma_{\vert\Sigma_{\bf p}}$ is fixed point free and 
$F^{\Sigma_{\bf p},\sigma}_2$ is nowhere dense in $\Sigma_{\bf p}$, Lemmas \ref{orbthree} and \ref{orbtwo} show that $h$ sends any orbit onto an orbit of the same size. This shows that the number of orbits of 
$\Sigma_{\bf p}$ is at most the finite number of orbits of $\Sigma_{{\bf p}'}$, by injectivity of $h$. By Lemma \ref{minSigmap}, 
$(\Sigma_{{\bf p}'},G_\sigma )\preceq^i_c(\Sigma_{\bf p},G_\sigma )$, so that $\Sigma_{\bf p}$ and 
$\Sigma_{{\bf p}'}$ have the same number of orbits, in bijection by $h$. In particular, $\Sigma_{\bf p}$ and 
$\Sigma_{{\bf p}'}$ have the same number of finite orbits, i.e., $l\! =\! l'$. This also shows that if 
$\Lambda\! :=\!\mbox{max}_{i\leq l}~\lambda_i$, then $\Lambda\! =\!\mbox{max}_{i\leq l}~\lambda'_i$. Note also that 
$m\! <\!\lambda_0\!\leq\!\Lambda$ and, similarly, $m'\! <\!\Lambda$. This shows that 
${\bf p}'\!\in\!\{ l\}\!\times\! (\Lambda\! +\! 1)^{l+1}\!\times\!\Lambda\!\times\! 2^l$, so that 
$\mathcal{F}_{\bf p}$ is finite, and 
${\bf m}_{\bf p}\! :=\!\mbox{min}_{\text{lex}}~\mathcal{F}_{\bf p}$ is defined.\medskip

\noindent\bf Claim.\it\ There is $\mathcal{P}_a\!\subseteq\! \mathcal{P}$, obtained by choosing $\mbox{min}_{\text{lex}}~\mathcal{F}_{\bf p}$ for each ${\bf p}\!\in\!\mathcal{P}$, with the properties that $\{ (\Sigma_{\bf m},G_\sigma )\mid {\bf m}\!\in\!\mathcal{P}_a\}$ is a basis for $\{ (\Sigma_{\bf p},G_\sigma )\mid {\bf p}\!\in\!\mathcal{P}\}$ and 
$$\{ (2q\! +\! 3,C_{2q+3})\mid q\!\in\!\omega\}\cup\{ (\mathbb{X}_1,G_{f_1})\}\cup
\{ (\Sigma_{\bf m},G_\sigma )\mid {\bf m}\!\in\!\mathcal{P}_a\}$$ 
is an antichain.\rm

\vfill\eject\medskip

 Indeed, let 
$\mathcal{P}_a\! :=\!\{ {\bf m}_{\bf p}\mid {\bf p}\!\in\!\mathcal{P}\}$, so that 
$\{ (\Sigma_{\bf m},G_\sigma )\mid {\bf m}\!\in\!\mathcal{P}_a\}$ is an antichain basis for 
$\{ (\Sigma_{\bf p},G_\sigma )\mid {\bf p}\!\in\!\mathcal{P}\}$. As $\mathbb{X}_1$ and the $\Sigma_{\bf p}$'s are infinite, the odd cycles are not above the other graphs. As these graphs contain no odd cycle, they are not below them. As a homomorphism sends an odd cycle of cardinality $l$ into an odd cycle of cardinality at most $l$ and by injectivity, the odd cycles form antichain. We saw that $\sigma_{\vert\Sigma_{\bf p}}$ is fixed point free, so that $(\Sigma_{\bf p},G_\sigma )$ has a continuous 
$\aleph_0$-coloring, by Corollary \ref{corfp}. Corollaries \ref{corfp} and \ref{exfp} imply that 
$(\mathbb{X}_1,G_{f_1})$ has no continuous $\aleph_0$-coloring. Thus the $(\Sigma_{\bf p},G_\sigma )$'s are not above 
$(\mathbb{X}_1,G_{f_1})$. As $(\Sigma_{\bf p},G_\sigma )$ has an infinite orbit, the orbits of $(\mathbb{X}_1,G_{f_1})$ have size at most two and an orbit has to be sent into an orbit, $(\mathbb{X}_1,G_{f_1})$ is not above the 
$(\Sigma_{\bf p},G_\sigma )$'s, by injectivity.\hfill{$\diamond$}\medskip

 By the claim, $\{ (2q\! +\! 3,C_{2q+3})\mid q\!\in\!\omega\}\cup\{ (\mathbb{X}_1,G_{f_1})\}\cup
\{ (\Sigma_{\bf m},G_\sigma )\mid {\bf m}\!\in\!\mathcal{P}_a\}$ is our antichain basis.\hfill{$\square$}

\subsection{$\!\!\!\!\!\!$ Subshifts}\indent

 We first prove a lemma useful to prove Theorem \ref{mainshift} (a). Note that the fixed points of the shift are the constant sequences, of the form $a^\mathbb{Z}\!\in\!\Sigma$ with $a\!\in\! A$.

\begin{lemm} \label{evcst} Let $\Sigma\!\subseteq\! A^\mathbb{Z}$ be a countable two-sided subshift with Cantor-Bendixson rank at most two, $l\!\in\!\omega$, $a_0,\cdots\! ,a_l\!\in\! A$, and $(x_n)_{n\in\omega}$ be an injective sequence of points of 
$\Sigma$ converging to $(a_0\!\cdots\! a_l)^\mathbb{Z}$. Then the sequence 
$\big(\mbox{Orb}(x_n)\big)_{n\in\omega}$ is eventually constant, and we can find 
$s\!\in\! A^{l+1}\!\setminus\{ (a_0\!\cdots\! a_l)\}$ and $\gamma\!\in\! A^\omega$ with 
$(a_0\!\cdots\! a_l)^{-\infty}\!\cdot\! s\gamma\!\in\!\Sigma$ or $\gamma^{-1}s\!\cdot\! (a_0\!\cdots\! a_l)^\infty\!\in\!\Sigma$.\end{lemm}

\noindent\bf Proof.\rm\ For the last assertion, we can apply Lemma \ref{fixshif}. We argue by contradiction, so that we may assume that the sequence $\big(\mbox{Orb}(x_n)\big)_{n\in\omega}$ is injective. We may also assume, for example, that 
$${x_n}_{[-k_n(l+1),k_n(l+1)-1]}\! =\! (a_0\!\cdots\! a_l)^{2k_n}\mbox{,}$$ 
${x_n}_{[k_n(l+1),(k_n+1)(l+1)-1]}$ is a constant $s\!\not=\! (a_0\!\cdots\! a_l)$, and $k_n\!\rightarrow\!\infty$. By compactness, we may assume that the sequence $({x_n}_{[k_n(l+1),\infty )})_{n\in\omega}$ converges to some $s\gamma\!\in\! A^\omega$. We put $x\! :=\! (a_0\!\cdots\! a_l)^{-\infty}\!\cdot\! s\gamma$, so that $x\!\in\!\Sigma$. Note that $x$ is the limit of 
$\big(\sigma^{k_n(l+1)}(x_n)\big)_{n\in\omega}$. As the sequence $\big(\mbox{Orb}(x_n)\big)_{n\in\omega}$ is injective, this sequence $\big(\sigma^{k_n(l+1)}(x_n)\big)_{n\in\omega}$ is also injective, so that $x\!\in\!\Sigma'$. Thus $\mbox{Orb}(x)$ is finite of cardinality $\kappa$, and $x\! =\!\sigma^{-\kappa (l+1)}(x)$, contradicting $s\!\not=\! (a_0\!\cdots\! a_l)$.\hfill{$\square$}\medskip

 We now provide an antichain basis when fixed points exist.

\begin{lemm} \label{basisfpshift} Let $\Sigma$ be a countable two-sided subshift with Cantor-Bendixson rank at most two such that $F_{\sigma_{\vert\Sigma}}$ is not open in $\Sigma$. Then there is $n\!\in\!\omega$ with the property that 
$({}_n\Sigma ,G_\sigma )\preceq^i_c(\Sigma ,G_\sigma )$. Moreover, 
$\{ ({}_n\Sigma ,G_\sigma )\mid n\!\in\!\omega\}$ is a $\preceq^i_c$-antichain.\end{lemm}

\noindent\bf Proof.\rm\ Assume that 
$\Sigma\!\subseteq\! A^\mathbb{Z}$ and $(x_n)_{n\in\omega}$ is an injective sequence of points of 
$\Sigma\!\setminus\!\{ a^\mathbb{Z}\mid a\!\in\! A\}$ converging to $a^\mathbb{Z}\!\in\!\Sigma$. By Lemma \ref{evcst}, we may assume that $\mbox{Orb}(x_n)\! =\!\mbox{Orb}(x_0)$ for each $n$, and that $x_0$ is of the form 
$a^{-\infty}\!\cdot\! b\gamma$ or $\gamma^{-1}b\!\cdot\! a^\infty$ with $b\!\in\! A\!\setminus\!\{ a\}$ and $\gamma\!\in\! A^\omega$, so that $x_0\!\in\!\Sigma\!\setminus\!\Sigma'$. Lemma \ref{twosides} provides $y^-,y^+\!\in\!\Sigma'$, and at least one of them is $a^\mathbb{Z}$. The other one is of the form $(a_0\!\cdots\! a_l)^\mathbb{Z}$, where $l\!\in\!\omega$ and $a_0,\cdots ,a_l\!\in\! A$. If $l\!\geq\! 2$ is even, then the map 
$(0\!\cdots\! l)^\mathbb{Z}\!\mapsto\! (a_0\!\cdots\! a_l)^\mathbb{Z}$ is a witness for the fact that 
$({}_{l+1}\Sigma ,G_\sigma )\preceq^i_c(\Sigma ,G_\sigma )$.\medskip

 If $l\! =\! 0$ and $a_0\! =\! a$, then the map defined by $0^\mathbb{Z}\!\mapsto\! a^\mathbb{Z}$ and 
$\sigma^k(0^{-\infty}\!\cdot\! 10^\infty )\!\mapsto\!\sigma^k(x_0)$ when $k\!\in\!\mathbb{Z}$ is a witness for the fact that 
$({}_0\Sigma ,G_\sigma )\preceq^i_c(\Sigma ,G_\sigma )$.\medskip

 If $l\! =\! 0$ and $a_0\!\not=\! a$, then the map defined by\medskip
 
\noindent - $0^\mathbb{Z}\!\mapsto\! a^\mathbb{Z}$, $1^\mathbb{Z}\!\mapsto\! a_0^\mathbb{Z}$ and 
$\sigma^k(0^{-\infty}\!\cdot\! 1^\infty )\!\mapsto\!\sigma^k(x_0)$ when $k\!\in\!\mathbb{Z}$ and 
$x_0\! =\! a^{-\infty}\!\cdot\! b\gamma$,\smallskip

\noindent - $0^\mathbb{Z}\!\mapsto\! a_0^\mathbb{Z}$, $1^\mathbb{Z}\!\mapsto\! a^\mathbb{Z}$ and 
$\sigma^k(0^{-\infty}\!\cdot\! 1^\infty )\!\mapsto\!\sigma^k(x_0)$ when $k\!\in\!\mathbb{Z}$ and 
$x_0\! =\!\gamma^{-1}b\!\cdot\! a^\infty$,\medskip

\noindent is a witness for the fact that $({}_1\Sigma ,G_\sigma )\preceq^i_c(\Sigma ,G_\sigma )$.\medskip

 If $l$ is odd, then the map defined by\medskip
 
\noindent - $0^\mathbb{Z}\!\mapsto\! a^\mathbb{Z}$, $\sigma^i\big( (1,\cdots ,l\! +\! 1)^\mathbb{Z}\big)\!\mapsto\!\sigma^i\big( (a_0\!\cdots\! a_l)^\mathbb{Z}\big)$ when $i\!\leq\! l$ and 
$\sigma^k(0^{-\infty}\!\cdot\! (1,\cdots ,l\! +\! 1)^\infty )\!\mapsto\!\sigma^k(x_0)$ when $k\!\in\!\mathbb{Z}$ and 
$x_0\! =\! a^{-\infty}\!\cdot\! b\gamma$,\smallskip

\noindent - $0^\mathbb{Z}\!\mapsto\! a^\mathbb{Z}$, 
$\sigma^i\big( (1,\cdots ,l\! +\! 1)^\mathbb{Z}\big)\!\mapsto\!\sigma^{-i}\big( (a_0\!\cdots\! a_l)^\mathbb{Z}\big)$ when $i\!\leq\! l$ and $\sigma^k(0^{-\infty}\!\cdot\! (1,\cdots ,l\! +\! 1)^\infty )\!\mapsto\!\sigma^{-k}(x_0)$ when $k\!\in\!\mathbb{Z}$ and 
$x_0\! =\!\gamma^{-1}b\!\cdot\! a^\infty$,\medskip

\noindent is a witness for the fact that $({}_{l+1}\Sigma ,G_\sigma )\preceq^i_c(\Sigma ,G_\sigma )$.\medskip

 By Theorem \ref{main}, the odd cycles $({}_{2q+3}\Sigma ,G_\sigma )$ form antichain. The other $({}_n\Sigma ,G_\sigma )$'s are infinite and contain no odd cycle, so they are incomparable with the odd cycles. Assume, towards a contradiction, that 
$m\!\not=\! n$ and $({}_m\Sigma ,G_\sigma )\preceq^i_c({}_n\Sigma ,G_\sigma )$ with 
$h$ as a witness. The previous discussion shows that we may assume that $m,n$ are not of the form $2q\! +\! 3$. This implies that the two subshifts have a unique infinite orbit, which is dense. Note that an orbit has to be sent into an orbit. In particular, the infinite orbit 
$\mbox{Orb}(x_m)$ of $({}_m\Sigma ,G_\sigma )$ has to be sent into that $\mbox{Orb}(x_n)$ of $({}_n\Sigma ,G_\sigma )$, by injectivity.  Note that $h[\mbox{Orb}(x_m)]\! =\!\mbox{Orb}(x_n)$, and $h\!\circ\!\sigma\! =\!\sigma\!\circ\! h$ on 
$\mbox{Orb}(x_m)$ or $h\!\circ\!\sigma\! =\!\sigma^{-1}\!\circ\! h$ on $\mbox{Orb}(x_m)$, by Lemma \ref{orbthree}. In particular, $h[{}_m\Sigma\!\setminus\!\mbox{Orb}(x_m)]\!\subseteq\! {}_n\Sigma\!\setminus\!\mbox{Orb}(x_n)$ by injectivity, and 
$h\!\circ\!\sigma\! =\!\sigma\!\circ\! h$ or $h\!\circ\!\sigma\! =\!\sigma^{-1}\!\circ\! h$, by density. Thus 
$m\! <\! n$, by injectivity. If $m\! =\! 0$, then $\big(\sigma^{{\bf d} i}(x_m)\big)_{i\in\omega}$ converges to 
$0^\mathbb{Z}$ for each ${\bf d}\!\in\!\{ -,+\}$, so that $\Big(\sigma^{-i}\big( h(x_m)\big)\Big)_{i\in\omega}$, 
$\Big(\sigma^i\big( h(x_m)\big)\Big)_{i\in\omega}$ have the same limit, which cannot be. If $m\! =\! 1$, then 
$\big(\sigma^{{\bf d} i}(x_m)\big)_{i\in\omega}$ converges for each ${\bf d}\!\in\!\{ -,+\}$, so that 
$\Big(\sigma^{{\bf d} i}\big( h(x_m)\big)\Big)_{i\in\omega}$ also converges for each 
${\bf d}\!\in\!\{ -,+\}$, which cannot be. If $m\! =\! 2q\! +\! 2$, then 
$h\big[\mbox{Orb}\big( (1,\cdots ,m)^\mathbb{Z}\big)\big]\! =\!
\mbox{Orb}\big( (1,\cdots ,n)^\mathbb{Z}\big)$, which cannot be.\hfill{$\square$}\smallskip

\noindent\bf Proof of Theorem \ref{mainshift} (a).\rm\ Let $\Sigma$ be a countable two-sided subshift with Cantor-Bendixson rank at most two such that $(\Sigma ,G_\sigma )$ has no continuous 2-coloring. If $F_{\sigma_{\vert\Sigma}}$ is not open in 
$\Sigma$, then Lemma \ref{basisfpshift} provides $n\!\in\!\omega$ such that 
$({}_n\Sigma ,G_\sigma )\preceq^i_c(\Sigma ,G_\sigma )$. If $F_{\sigma_{\vert\Sigma}}$ is open in $\Sigma$, then Corollary \ref{corfp} implies that $(\Sigma\!\setminus\! F_{\sigma_{\vert\Sigma}},G_\sigma )$ has no continuous 2-coloring. Thus we may assume that $\sigma_{\vert\Sigma}$ is fixed point free. We may also assume that 
$(\Sigma ,G_\sigma )$ contains no odd cycle. It remains to apply Theorem \ref{bas} to get the basis result.\smallskip

 The set $\mathcal{P}_a\!\subseteq\! \mathcal{P}$ provided by the claim in the proof of Theorem \ref{main} has the properties that $\{ (\Sigma_{\bf m},G_\sigma )\mid {\bf m}\!\in\!\mathcal{P}_a\}$ is a basis for 
$\{ (\Sigma_{\bf p},G_\sigma )\mid {\bf p}\!\in\!\mathcal{P}\}$ and\medskip
 
\centerline{$\{ ({}_{2q+3}\Sigma ,G_\sigma )\mid q\!\in\!\omega\}\cup\{ (\Sigma_{\bf m},G_\sigma )\mid {\bf m}\!\in\!\mathcal{P}_a\}$}\medskip
 
\noindent is an antichain. By Lemma \ref{basisfpshift}, 
$\{ ({}_n\Sigma ,G_\sigma )\mid n\!\in\!\omega\}$ is also an antichain. We saw that 
$\sigma_{\vert\Sigma_{\bf p}}$ is fixed point free, so that $(\Sigma_{\bf p},G_\sigma )$ has a continuous $\aleph_0$-coloring, by Corollary \ref{corfp}. We also saw after Corollary \ref{corfp} that $({}_n\Sigma ,G_\sigma )$ has no continuous $\aleph_0$-coloring if $n$ is 1 or even. Thus the $(\Sigma_{\bf p},G_\sigma )$'s are not above 
$({}_n\Sigma ,G_\sigma )$ if $n$ is 1 or even.

\vfill\eject

 Let us check that $({}_n\Sigma ,G_\sigma )$ is not above the 
$(\Sigma_{\bf p},G_\sigma )$'s if $n$ is 1 or even. We argue by contradiction, which provides 
$h\! :\!\Sigma_{\bf p}\!\rightarrow\! {}_n\Sigma$. Note that $h[\Sigma_{\bf p}]$ is compact and does not contain a fixed point of 
$\sigma_{\vert {}_n\Sigma}$. This provides a continuous 2-coloring of 
$(h[\Sigma_{\bf p}],(h\!\times\! h)[G_\sigma ])$, and thus one of $(\Sigma_{\bf p},G_\sigma )$, which is not possible by Claim 1 in the proof of Lemma \ref{minSigmap}.  This shows that 
$\{ ({}_n\Sigma ,G_\sigma )\mid n\!\in\!\omega\}\cup\{ (\Sigma_{\bf m},G_\sigma )\mid {\bf m}\!\in\!\mathcal{P}_a\}$ is our antichain basis.\hfill{$\square$}

\section{$\!\!\!\!\!\!$ Basis of size continuum}\indent

 Recall the definition of $\mathcal{C}_\kappa$, $\mathcal{H}_\kappa$, and 
$\mathcal{S}_\kappa$, before the proof of Theorem \ref{mainclosed} (a), Lemma \ref{orbtwo}, and Lemma \ref{compa1shift}, respectively. Theorem 1.17 (b) in [L] provides a $\preceq^i_c$-antichain $\big( (\Sigma_\nu ,G_\sigma )\big)_{\nu\in 2^\omega}$ made of countable two-sided subshifts with Cantor-Bendixson rank three such that 
$\sigma_{\vert\Sigma_\nu}$ is fixed point free (and thus $G_{\sigma_{\vert\Sigma_\nu}}$ is closed), $(\Sigma_\nu ,G_\sigma)$ has a continuous $3$-coloring and is $\preceq^i_c$-minimal in $\mathcal{C}_2$ and in $\mathcal{H}_2$. This proves Theorem \ref{mainclosed} (b) for $\kappa\! =\! 2$, finishes the proof of Theorem \ref{mainhomeo} (b), proves Theorem \ref{mainsubshift} (b) for 
$\kappa\! =\! 2$, and proves the second part of Theorems \ref{main} (b) (the first part comes from Theorem \ref{mainhomeo}), and \ref{mainshift} (b) (the first part comes from Theorem \ref{mainsubshift}). The proof of Theorem \ref{mainsubshift} (b) for 
$\kappa\!\geq\! 3$ is partly similar, so we recall the construction of Theorem 1.17 (b) in [L].\smallskip

\noindent\bf Notation.\rm\ Let $\alpha_0\! :=\! (01)^{-\infty}\!\cdot\! (01)^\infty$, 
$\alpha_1\! :=\! (01)^{-\infty}\!\cdot\! 1^2(01)^\infty$, $Q\! :=\! (q_j)_{j\in\omega}\!\in\!\omega^\omega$ converging to infinity, and 
$\beta_Q\! :=\! (01)^{-\infty}\!\cdot\! 1{^\frown}_{j\in\omega}~\big((01)^{q_j}1^2 \big)$. This defines $\Sigma_Q\! =\!
\bigcup_{m\leq 1}~\mbox{Orb}_\sigma (\alpha_m)\cup\mbox{Orb}_\sigma (\beta_Q)$.\medskip

 Note that $\Sigma_Q$ is a countable two-sided subshift. By Claim 8 in the proof of Theorem 1.17 (b) in [L], $\Sigma_Q$ has Cantor-Bendixson rank three, and the remark after this claim shows that $(\Sigma_Q,G_\sigma)$ has a continuous $3$-coloring. The proof of Theorem 1.17 (b) in [L] also shows the minimality of $(\Sigma_Q,G_\sigma)$. In order to get the antichain, we consider the sequence $(p_n)_{n\in\omega}$ of prime numbers. We set, for $\nu\!\in\! 2^\omega$ and 
$n\!\in\!\omega$, $q^\nu_0\! :=\! 0$ and $q^\nu_{n+1}\! :=\! p_0^{\nu (0)+2}\!\cdots\! p_n^{\nu (n)+2}\! -\! 1$, which defines 
$Q^\nu\!\in\!\omega^\omega$ converging to infinity. Then $\Sigma_\nu\! :=\!\Sigma_{Q_\nu}$.\medskip

\noindent\bf Proof of Theorem \ref{mainsubshift} (b) for $\kappa\!\geq\! 3$.\rm\ Let $L\! :=\! (l_j)_{j\in\omega}\!\in\!\omega^\omega$ converging to infinity, $\gamma_0\! :=\! 0^\mathbb{Z}$, $\gamma_1\! :=\! (01)^\mathbb{Z}$, 
$\gamma_2\! :=\! (01)^{-\infty}\!\cdot\! 1^2(01)^\infty$, and 
$\delta_L\! :=\! 0^{-\infty}\!\cdot\! {^\frown}_{j\in\omega}~\big( (01)^{l_j}1^2\big)$. This defines as above 
$\Sigma_L\! =\!\bigcup_{m\leq 2}~\mbox{Orb}_\sigma (\gamma_m)\cup\mbox{Orb}_\sigma (\delta_L)$. Note that $\Sigma_L$ is a countable two-sided subshift. As $\gamma_0\!\in\! F_{\sigma_{\vert\Sigma_L}}$ is the limit of 
$\big(\sigma^{-n}(\delta_L)\big)_{n\in\omega}\!\in\! (\Sigma_L\!\setminus\! F_{\sigma_{\vert\Sigma_L}})^\omega$, 
$F_{\sigma_{\vert\Sigma_L}}$ is not an open subset of $\Sigma_L$. By Corollary \ref{corfp}, there is no continuous 
$\aleph_0$-coloring of $G_{\sigma_{\vert\Sigma_L}}$, so that $(\Sigma_L,G_\sigma )\!\in\!\mathcal{S}_\kappa$.\medskip

\noindent\emph{Claim.}\it\ $(\Sigma_L,G_\sigma )$ is $\preceq^i_c$-minimal in $\mathcal{S}_\kappa$.\rm\medskip

 Indeed, let $(\Sigma ,G_\sigma)\!\in\!\mathcal{S}_\kappa$ such that 
$(\Sigma ,G_\sigma )\preceq^i_c(\Sigma_L,G_\sigma )$ with $h$ as a witness. We first prove that there is $(\Sigma',G_\sigma)\!\in\!\mathcal{S}_\kappa$ with 
$(\Sigma',G_\sigma )\preceq^i_c(\Sigma ,G_\sigma )$ and 
$F_2^{\Sigma',\sigma_{\vert\Sigma'}}$ is nowhere dense in 
$\Sigma'$. We argue as in the proof of the claim in the proof of Theorem \ref{mainhomeo} (a), by contradiction. We inductively construct a strictly $\subseteq$-decreasing sequence 
$(\Sigma_\xi )_{\xi <\aleph_1}$ such that $\Sigma_0\! =\!\Sigma$, $\Sigma_\xi$ is $\sigma$-invariant and 
$(\Sigma_\xi ,G_\sigma\cap\Sigma_\xi^2)\!\in\!\mathcal{S}_\kappa$, which will contradict the fact that $\Sigma$ is a 0DMC space. Assume that $\Sigma_\xi$ is constructed. Note that $F_{\sigma_{\vert\Sigma_\xi}}$ is finite. Let 
$I\! :=\!\{ x\!\in\! F_{\sigma_{\vert\Sigma_\xi}}\mid x\mbox{ is isolated in }\Sigma_\xi\}$. Note that $I$ is a finite $\sigma$-invariant  clopen subset of $\Sigma_\xi$, $G_{\sigma_{\vert\Sigma_\xi}}\! =\! G_{\sigma_{\vert\Sigma_\xi\setminus I}}$, 
$(\Sigma_\xi\!\setminus\! I,G_\sigma )\!\in\!\mathcal{S}_\kappa$, and 
$(\Sigma_\xi\!\setminus\! I,G_\sigma )\preceq^i_c(\Sigma_\xi ,G_\sigma )$. So, restricting $\Sigma_\xi$ to 
$\Sigma_\xi\!\setminus\! I$ if necessary, we may assume that $F_{\sigma_{\vert\Sigma_\xi}}$ is nowhere dense in $\Sigma_\xi$. Note that $F^{\Sigma_\xi ,\sigma_{\vert\Sigma_\xi}}_2$ is closed and not nowhere dense in $\Sigma_\xi$.

\vfill\eject

 This gives a nonempty clopen subset $C$ of $\Sigma_\xi$ with the property that 
$C\!\subseteq\! F^{\Sigma_\xi,\sigma_{\vert\Sigma_\xi}}_2\!\setminus\! F_{\sigma_{\vert\Sigma_\xi}}$. Note that 
$U\! :=\! C\cup\sigma [C]$ is a nonempty clopen $\sigma$-invariant subset of $\Sigma_\xi$ contained in 
$F^{\Sigma_\xi ,\sigma_{\vert\Sigma_\xi}}_2\!\setminus\! F_{\sigma_{\vert\Sigma_\xi}}$. In particular, $U$ is a ODM separable space and $\sigma_{\vert U}$ is a fixed point free continuous involution. Proposition 7.5 in [L] provides a continuous 2-coloring of $(U,G_{\sigma_{\vert U}})$. All this implies that 
$\Sigma_{\xi +1}\! :=\!\Sigma_\xi\!\setminus\! U\!\subsetneqq\!\Sigma_\xi$, $\Sigma_{\xi +1}$ is $\sigma$-invariant and 
$(\Sigma_{\xi +1},G_\sigma\cap\Sigma_{\xi +1}^2)\!\in\!\mathcal{S}_\kappa$. If $(\lambda_p)_{p\in\omega}$ is strictly increasing and $\lambda\! =\!\mbox{sup}_{p\in\omega}~\lambda_p$ is a limit ordinal, then we set 
$\Sigma_\lambda\! :=\!\bigcap_{p\in\omega}~\Sigma_{\lambda_p}$. As $\kappa\!\geq\! 3$, the 
$(\Sigma_{\lambda_p},G_\sigma )$'s are in $\mathcal{S}_{\aleph_0}$, by Theorem 1.12 in [L]. By Lemma \ref{compa1shift}, 
$(\Sigma_\lambda ,G_\sigma )\!\in\!\mathcal{S}_\kappa$. As $\Sigma_\lambda\!\subsetneqq\!\Sigma_{\lambda_p}$ for each 
$p\!\in\!\omega$, we are done. In other words, we may assume that $F_2^{\Sigma,\sigma_{\vert\Sigma}}$ is nowhere dense in $\Sigma$. By Lemmas \ref{orbthree}, \ref{orbtwo}, $h$ sends an orbit of size at least two onto an orbit of the same size. We set $P\! :=\! h^{-1}\big(\mbox{Orb}_\sigma (\delta_L)\big)$. If $P$ is contained in $F_{\sigma_{\vert\Sigma}}$, then $P$ is finite since a two-sided subshift has only finitely many fixed points. Moreover, these points are isolated in $\Sigma$ since so are the elements of $\mbox{Orb}_\sigma (\delta_L)$ in $\Sigma_L$. This shows that $P$ is a finite $\sigma$-invariant clopen subset of $\Sigma$. In particular, $\Sigma\!\setminus\! P$ is a $\sigma$-invariant clopen subset of $\Sigma$ and 
$(\Sigma\!\setminus\! P,G_\sigma )\!\in\!\mathcal{S}_\kappa$. On the other hand, 
$(\Sigma\!\setminus\! P,G_\sigma )\preceq^i_c(\bigcup_{m\leq 2}~\mbox{Orb}_\sigma (\gamma_m),G_\sigma )$, which has a continuous $2$-coloring, which cannot be. This shows that $P$ contains an element of 
$\Sigma\!\setminus\! F_{\sigma_{\vert\Sigma}}$. Thus the dense set $\mbox{Orb}_\sigma (\delta_L)$ is contained in the compact set $h[\Sigma ]$ by the previous size argument, proving that $h$ is onto, and thus a homeomorphism by compactness. In particular, $P$ is a dense orbit $\mbox{Orb}_\sigma (x)$. By Lemma \ref{orbthree}, there is 
$\theta\!\in\!\{ -1,1\}$ such that $h\!\circ\!\sigma_{\vert\Sigma}\! =\!\sigma_{\vert\Sigma_L}^\theta\!\circ\! h$ on 
$\mbox{Orb}_\sigma (x)$, and thus on $\Sigma$. In particular, 
$(h\!\times\! h)[G_{\sigma_{\vert\Sigma}}]\! =\! G_{\sigma_{\vert\Sigma_L}}$ and 
$(\Sigma_L,G_\sigma )\preceq^i_c(\Sigma ,G_\sigma )$ with $h^{-1}$ as a witness.
\hfill{$\diamond$}\medskip

 We now define $L^\nu$ as we defined $Q^\nu$ just before this proof. It remains to check that the family 
$\big( (\Sigma_{L^\nu},G_{\sigma_{\vert\Sigma_{L^\nu}}})\big)_{\nu\in 2^\omega}$ is a 
$\preceq^i_c$-antichain. Assume, towards a contradiction, that $\nu\!\not=\!\nu'$ and 
$$(\Sigma_{L^\nu},G_{\sigma_{\vert\Sigma_{L^\nu}}})\preceq^i_c(\Sigma_{L^{\nu'}},G_{\sigma_{\vert\Sigma_{L^{\nu'}}}})$$ 
with $h$ as a witness. Let $m_0$ be minimal with the property that 
$\nu (m_0)\!\not=\!\nu'(m_0)$. By minimality of 
$(\Sigma_{L^{\nu'}},G_{\sigma_{\vert\Sigma_{L^{\nu'}}}})$, we may assume that $\nu (m_0)$ is smaller than $\nu'(m_0)$. By Lemma \ref{orbthree}, $h[\mbox{Orb}_\sigma (\gamma_2)]$, $h[\mbox{Orb}_\sigma (\delta_{L^{\nu}})]$ are disjoint infinite orbits in 
$\Sigma_{L^{\nu'}}$, so they are $\mbox{Orb}_\sigma (\gamma_2)$, $\mbox{Orb}_\sigma (\delta_{L^{\nu'}})$. As 
$\mbox{Orb}_\sigma (\delta_{L^{\nu'}})$ is dense in $\Sigma_{L^{\nu'}}$, the compact set $h[\Sigma_{L^{\nu}}]$ is 
$\Sigma_{L^{\nu'}}$, so that $h$ is a homeomorphism from $\Sigma_{L^{\nu}}$ onto $\Sigma_{L^{\nu'}}$. Moreover, $h$ is a witness for the fact that $\sigma_{\vert\Sigma_{L^{\nu}}}$ and $\sigma_{\vert\Sigma_{L^{\nu'}}}$ are flip-conjugate, by density of $\mbox{Orb}_\sigma (\delta_{L^{\nu}})$ in $\Sigma_{L^{\nu}}$ and Lemma \ref{orbthree}. In particular, 
$h[\Sigma'_{L^{\nu}}]\! =\!\Sigma'_{L^{\nu'}}$ and $h[\Sigma''_{L^{\nu}}]\! =\!\Sigma''_{L^{\nu'}}$, so that 
$$h[\mbox{Orb}_\sigma (\gamma_1)]\! =\!\mbox{Orb}_\sigma (\gamma_1)\mbox{,}$$ 
$h[\mbox{Orb}_\sigma (\gamma_2)]\! =\!\mbox{Orb}_\sigma (\gamma_2)$, 
$h[\mbox{Orb}_\sigma (\delta_{L^{\nu}})]$ is $\mbox{Orb}_\sigma (\delta_{L^{\nu'}})$ and $h(\gamma_0)\! =\!\gamma_0$.  
This gives $n_0,n_1\!\in\!\mathbb{Z}$ with $h(\gamma_2)\! =\!\sigma^{n_1}(\gamma_2)$ and 
$h(\delta_{L^{\nu}})\! =\!\sigma^{n_0}(\delta_{L^{\nu'}})$. We then set, for $r\!\in\!\omega$, 
${K^\nu_r\! :=\!\big(\Sigma_{j<r}~(2l_j\! +\! 2)\big)\! +\! 2l_r}$.   
Note that the sequence $\big(\sigma^{K^\nu_r}(\delta_{L^{\nu}})\big)_{r\in\omega}$ converges to $\gamma_2$, so that 
$\Big( h\big(\sigma^{K^\nu_r}(\delta_{L^{\nu}})\big)\Big)_{r\in\omega}$ converges to 
$h(\gamma_2)\! =\!\sigma^{n_1}(\gamma_2)$. As 
$h\big(\sigma^{K^\nu_r}(\delta_{L^{\nu}})\big)\! =\!\sigma^{n_0\pm K^\nu_r}(\delta_{L^{\nu'}})$, 
this implies that $\big(\sigma^{n_0-n_1\pm K^\nu_r}(\delta_{L^{\nu'}})\big)_{r\in\omega}$ converges to $\gamma_2$. As 
$(K^\nu_r)_{r\in\omega}$ is strictly increasing, this implies that $\sigma_{\vert\Sigma_{L^{\nu}}}$ and 
$\sigma_{\vert\Sigma_{L^{\nu'}}}$ are conjugate and $\big(\sigma^{n_0-n_1+K^\nu_r}(\delta_{L^{\nu'}})\big)_{r\in\omega}$ converges to $\gamma_2$. In particular, 
$${\sigma^{n_0-n_1+K^\nu_r}(\delta_{L^{\nu'}})}_{[-2,2]}\! =\!
\big(\delta_{L^{\nu'}}(n_0\! -\! n_1\! +\! K^\nu_r\! -\! 2),\cdots ,\delta_{L^{\nu'}}(n_0\! -\! n_1\! +\! K^\nu_r\! +\! 2)\big)
\! =\! {\gamma_2}_{[-2,2]}\! =\! 01^30$$ 
if $r$ is large enough. Using similar notation, this implies that 
$n_0\! -\! n_1\! +\! K^\nu_r\!\in\!\{ K^{\nu'}_m\mid m\!\in\!\omega\}$ if $r$ is large enough.

\vfill\eject

 In particular, this gives, for $r$ large enough, $m\! <\! M\!\in\!\omega$ with $n_0\! -\! n_1\! +\! K^\nu_r\! =\! K^{\nu'}_m$ and 
$$n_0\! -\! n_1\! +\! K^\nu_{r+1}\! =\! K^{\nu'}_M.$$ 
Thus $K^\nu_{r+1}\! -\! K^\nu_r\! =\! 2l^\nu_{r+1}\! +\! 2\! =\!\Sigma_{m<j\leq M}~(2l^{\nu'}_j\! +\! 2)$ and 
$$p_0^{\nu (0)+2}\!\cdots\! p_r^{\nu (r)+2}\! =\! l^\nu_{r+1}\! +\! 1\! =\!\Sigma_{m\leq n<M}~(l^{\nu'}_{n+1}\! +\! 1)
\! =\!\Sigma_{m\leq n<M}~(p_0^{\nu'(0)+2}\!\cdots\! p_n^{\nu'(n)+2}).$$ 
We may assume that $r$ is large enough to ensure that $r,m\!\geq\! m_0$, which implies that $p_{m_0}^{\nu'(m_0)+2}$ divides $p_0^{\nu (0)+2}\!\cdots\! p_r^{\nu (r)+2}$, which cannot be since $\nu (m_0)\! <\!\nu'(m_0)$.\hfill{$\square$}\medskip

 It remains to prove Theorem \ref{mainclosed} (b) for $\kappa\!\geq\! 3$.\medskip

\noindent\bf Proof of Theorem \ref{mainclosed} (b) for $\kappa\!\geq\! 2$.\rm\ We set, for $\varepsilon\!\in\!\kappa$, 
$\varepsilon^+\! :=\! (\varepsilon\! +\! 1)\mbox{ mod }\kappa$. We then set, for $\nu\!\in\! (\omega\!\setminus\!\{ 0\} )^\omega$ and $j,k\!\in\!\omega$, $\beta_0^{\nu,k,j}\! :=\! 0^{2+j+\Sigma_{i<k}\nu (i)}1^\infty$ if $j\! <\!\nu (k)$, 
$\beta_\varepsilon^{\nu,2k+1,j}\! :=\!\varepsilon^{2+j+\Sigma_{i<2k+1}\nu (i)}(\varepsilon^+)^\infty$ if 
$0\! <\!\varepsilon\! <\!\kappa$ and $j\! <\!\nu (2k\! +\! 1)$, and 
$\beta_\varepsilon^{\nu,2k,j}\! :=\!\varepsilon^{2+j+\Sigma_{i<2k}\nu (i)}(\varepsilon^+)^\infty$ if 
$\kappa\!\leq\!\varepsilon\! <\!2\kappa\! -\! 1$ and $j\! <\!\nu (2k)$. This allows us to define the (countable) set of vertices\medskip
 
\leftline{$X_\nu\! :=\!\{\varepsilon^\infty\mid\varepsilon\!\in\! 2\kappa\! -\! 1\}\cup\{ 01^\infty\}\cup
\{\beta^{\nu ,k,j}_0\mid k\!\in\!\omega\wedge j\! <\!\nu (k)\} ~\cup$}\smallskip

\rightline{$\{\beta^{\nu ,2k+1,j}_\varepsilon\mid k\!\in\!\omega\wedge 0\! <\!\varepsilon\! <\!\kappa\wedge j\! <\!\nu (2k\! +\! 1)\}\cup\{\beta^{\nu ,2k,j}_\varepsilon\mid k\!\in\!\omega\wedge \kappa\!\leq\!\varepsilon\! <\! 2\kappa\! -\! 1\wedge j\! <\!\nu (2k)\}$}\medskip

\noindent and the set of edges

\medskip
 
\leftline{$G_\nu\! :=\! s(\{ (\varepsilon^\infty ,\eta^\infty )\mid\varepsilon\!\not=\!\eta\!\in\!\kappa\}\cup
\{ (\varepsilon^\infty ,\eta^\infty )\mid\varepsilon\!\not=\!\eta\!\in\!\{ 0\}\cup [\kappa ,2\kappa\! -\! 1)\}\cup\{ (0^\infty ,01^\infty )\} 
~\cup$}\smallskip

\rightline{$\{ (01^\infty ,\beta^{\nu ,0,0}_\varepsilon )\mid\kappa\!\leq\!\varepsilon\! <\! 2\kappa\! -\! 1\} ~\cup$}\smallskip

\rightline{$\{ (\beta^{\nu ,2k,j}_\varepsilon ,\beta^{\nu ,2k,j}_\eta )\mid k\!\in\!\omega\wedge
\varepsilon\!\not=\!\eta\!\in\!\{ 0\}\cup [\kappa ,2\kappa\! -\! 1)\wedge j\! <\!\nu (2k)\} ~\cup$}\smallskip

\rightline{$\{ (\beta^{\nu ,2k,\nu (2k)-1+j}_0,\beta^{\nu ,2k+1,j}_\varepsilon )\mid k\!\in\!\omega\wedge
0\!\not=\!\varepsilon\!\in\!\kappa\wedge j\! <\!\nu (2k\! +\! 1)\} ~\cup$}\smallskip

\rightline{$\{ (\beta^{\nu ,2k+1,j}_\varepsilon ,\beta^{\nu ,2k+1,j}_\eta )\mid k\!\in\!\omega\wedge
\varepsilon\!\not=\!\eta\!\in\!\kappa\wedge j\! <\!\nu (2k\! +\! 1)\} ~\cup$}\smallskip

\rightline{$\{ (\beta^{\nu ,2k+1,\nu (2k+1)-1+j}_0,\beta^{\nu ,2k+2,j}_\varepsilon )\mid k\!\in\!\omega\wedge 
\varepsilon\!\in\! [\kappa ,2\kappa\! -\! 1)\wedge j\! <\!\nu (2k\! +\! 2)\} ).$}\medskip

\noindent Note that $X_\nu$ is a closed subspace of $\kappa^\omega$, and thus a 0DMC space, with Cantor-Bendixson rank two. Also, the set $G_\nu$ is closed graph on $X_\nu$. If $c\! :\! X_\nu\!\rightarrow\!\kappa$ is a continuous coloring of $(X_\nu ,G_\nu )$, then, for example, $c(0^\infty )\! =\! 0$. This implies that $c(\beta_0^{\nu ,k,j})\! =\! 0$ if $k$ is big enough. Inductively, this implies that $c(\beta_0^{\nu ,k,j})\! =\! 0$ for each $k\!\in\!\omega$ since $0^n1^\infty$ and $0^{n+1}1^\infty$ have 
 $\kappa\! -\! 1$ common $G_\nu$-neighbors which are all pairwise $G_\nu$-related. Thus $c(01^\infty )\! =\! 0$. In particular, 
$c(01^\infty )\! =\! c(0^\infty )\! =\! 0$, contradicting $(0^\infty ,01^\infty )\!\in\! G_\nu$. This shows that 
$(X_\nu ,G_\nu )\!\in\!\mathcal{C}_\kappa$. For the minimality, it is enough to see that 
$(X_\nu ,G_\nu )\preceq^i_c(X,G)$ if $X\!\subseteq\! X_\nu$ and $G\!\subseteq\! G_\nu$ is a graph on $X$ such that 
$(X,G)\!\in\!\mathcal{C}_\kappa$, by compactness. The previous discussion shows that 
$$G\!\supseteq\! G_\nu\!\setminus\! s(\{ (\varepsilon^\infty ,\eta^\infty )\mid\varepsilon\!\not=\!\eta\!\in\!\kappa\}\cup
\{ (\varepsilon^\infty ,\eta^\infty )\mid\varepsilon\!\not=\!\eta\!\in\!\{ 0\}\cup [\kappa ,2\kappa\! -\! 1)\} ).$$ 
Indeed, if one edge $e$ in the difference is not in $G$, then we may assume that $e\!\notin\! s(\{ (0^\infty ,01^\infty )\} )$, and we can give the same color to the two vertices ($\varepsilon^{n+2}(\varepsilon^+)^\infty$ with $\varepsilon\!\not=\! 0$ for one of them) of $e$ and ensure that $c(0^{n+1}1^\infty )\!\not=\! 0$. Thus $G\! =\! G_\nu$ since $G$ is closed. This implies that 
$X\! =\! X_\nu$, and $(X_\nu ,G_\nu )\preceq^i_c(X,G)$.\medskip

\vfill\eject

 For instance, for $\kappa\! =\! 2,3$ and $\nu\!\in\! N_{132}$, this gives the following pictures.\bigskip

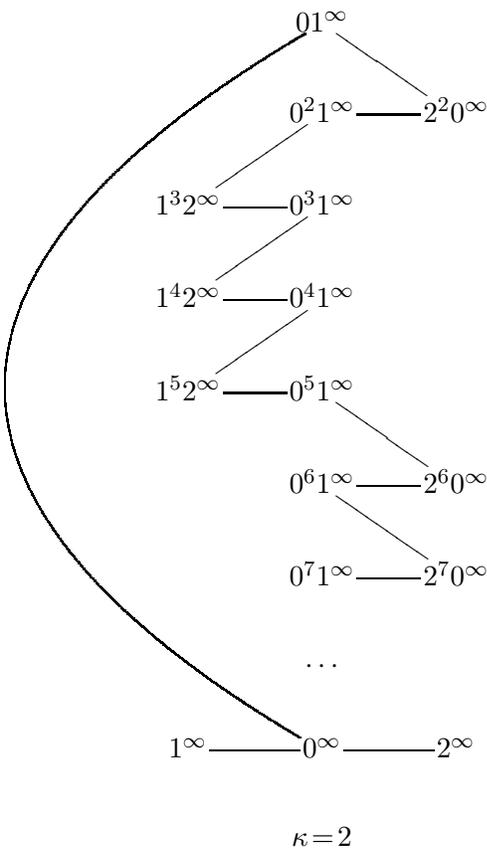
\centerline{~~~~~~~~~~~~~~~~~~~~~~~~~~~\scalebox{1.0}{\xymatrix@1{ 
& 01^\infty\ar@{-}@/_10pc/[dddddddd]\ar@{-}[dr] & \\ 
& 0^21^\infty\ar@{-}[dl]\ar@{-}[r] & 2^20^\infty \\
1^32^\infty\ar@{-}[r] & 0^31^\infty\ar@{-}[dl] & \\  
1^42^\infty\ar@{-}[r] & 0^41^\infty\ar@{-}[dl] & \\  
1^52^\infty\ar@{-}[r] & 0^51^\infty\ar@{-}[dr] & \\  
& 0^61^\infty\ar@{-}[r]\ar@{-}[dr] & 2^60^\infty \\  
& 0^71^\infty\ar@{-}[r] & 2^70^\infty & \\    
& \cdots & \\ 
1^\infty\ar@{-}[r] & 0^\infty\ar@{-}[r] & 2^\infty \\
& \kappa\! =\! 2 &  }}~~~~~~~~~~~~~~~~~~~~~~~~~~~~~~~~~~~~
\scalebox{0.71}
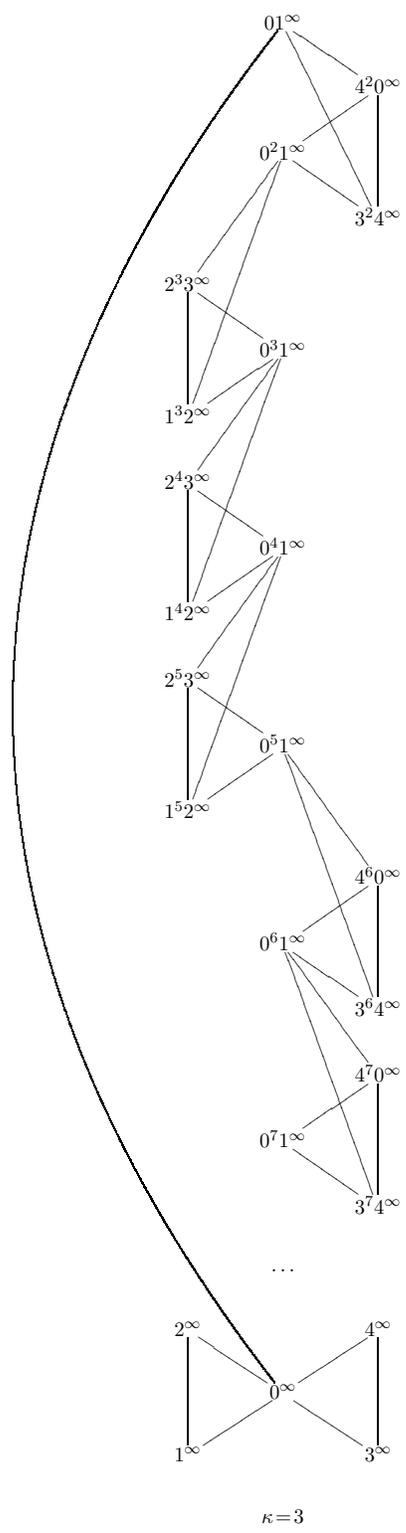
{\xymatrix@1{ 
& 01^\infty\ar@{-}@/_12pc/[ddddddddddddddddddddd]\ar@{-}[dr]\ar@{-}[dddr] & \\ 
& & 4^20^\infty\ar@{-}[dd] \\
& 0^21^\infty\ar@{-}[dr]\ar@{-}[ddl]\ar@{-}[ddddl]\ar@{-}[ur] & \\
& & 3^24^\infty \\
2^33^\infty\ar@{-}[dr]\ar@{-}[dd] & & \\  
& 0^31^\infty\ar@{-}[ddl]\ar@{-}[ddddl] & \\  
1^32^\infty\ar@{-}[ur] & & \\  
2^43^\infty\ar@{-}[dr]\ar@{-}[dd] & & \\  
& 0^41^\infty\ar@{-}[ddl]\ar@{-}[ddddl] & \\  
1^42^\infty\ar@{-}[ur] & & \\  
2^53^\infty\ar@{-}[dr]\ar@{-}[dd] & & \\  
& 0^51^\infty\ar@{-}[ddr]\ar@{-}[ddddr] & \\  
1^52^\infty\ar@{-}[ur] & & \\  
& & 4^60^\infty\ar@{-}[dd] \\  
& 0^61^\infty\ar@{-}[ur]\ar@{-}[dr]\ar@{-}[ddr]\ar@{-}[ddddr] & \\  
& & 3^64^\infty \\  
& & 4^70^\infty\ar@{-}[dd] & \\  
& 0^71^\infty\ar@{-}[ur]\ar@{-}[dr] & & \\    
& & 3^74^\infty & \\    
& \cdots & \\ 
2^\infty\ar@{-}[dr]\ar@{-}[dd] & & 4^\infty\ar@{-}[dd] \\ 
& 0^\infty\ar@{-}[ur]\ar@{-}[dr] & \\ 
1^\infty\ar@{-}[ur] & & 3^\infty \\
& \kappa\! =\! 3 & }}}

\vfill\eject

 Let us prove that $(X_\nu ,G_\nu )\not\preceq^i_c(X_{\nu'},G_{\nu'})$ if $(\nu ,\nu')\!\notin\!\mathbb{E}_t$, where 
$\mathbb{E}_t$ is the equivalence relation on $(\omega\!\setminus\!\{ 0\} )^\omega$ (introduced for example in [Do-J-K]) defined by 
$$(\nu ,\nu')\!\in\!\mathbb{E}_t\Leftrightarrow
\exists l,m\!\in\!\omega ~~\forall n\!\in\!\omega ~~\nu (l\! +\! n)\! =\!\nu (m\! +\! n).$$ 
We argue by contradiction, which gives $h\! :\! X_\nu\!\rightarrow\! X_{\nu'}$. Note that $0^\infty$ is the only limit point in $X_\nu$ with 
$\kappa\! +\! 1$ neighbors in $G_\nu$, so that $h(0^\infty )\! =\! 0^\infty$. The limit vertex $1^\infty$ is $G_\nu$-related to 
$\varepsilon^\infty$ if $1\!\not=\!\varepsilon\! <\!\kappa$ only, so that 
$$h[\{\varepsilon^\infty\mid\varepsilon\! <\!\kappa\} ],
h[\{\varepsilon^\infty\mid\varepsilon\!\in\!\{ 0\}\cup [\kappa ,2\kappa\! -\! 1)\} ]
\!\in\!\big\{\{\varepsilon^\infty\mid\varepsilon\! <\!\kappa\} ,
\{\varepsilon^\infty\mid\varepsilon\!\in\!\{ 0\}\cup [\kappa ,2\kappa\! -\! 1)\}\big\} .$$ 
Assume that $h(\varepsilon^\infty )\! =\!\eta_\varepsilon^\infty$. The continuity of $h$ implies that 
$$\left\{\!\!\!\!\!\!\!
\begin{array}{ll}
& h(\beta_0^{\nu ,k,j})\!\in\! N_0\mbox{,}\cr
& h(\beta_\varepsilon^{\nu ,2k+1,j})\!\in\! N_{\eta_\varepsilon}\mbox{ if }0\! <\!\varepsilon\! <\!\kappa\mbox{,}\cr
& h(\beta_\varepsilon^{\nu ,2k,j})\!\in\! N_{\eta_\varepsilon}\mbox{ if }\kappa\!\leq\!\varepsilon\! <\! 2\kappa\! -\! 1
\end{array}
\right.$$ 
if $k\!\geq\! k_0$ is big enough. Assume, for example, that $\eta_1\!\geq\!\kappa$, the case $\eta_1\! <\!\kappa$ being similar. Then $h(\beta_1^{\nu ,2k_0+1,0})$ is of the form $\beta_{\eta_1}^{\nu',2K_0,J_0}$. Thus 
$h(\beta_0^{\nu ,2k_0+1,0})\! =\!\beta_0^{\nu',2K_0,J_0}$. Then, for $0\! <\!\varepsilon\! <\!\kappa$ and inductively on 
$0\! <\! j\! <\!\nu (2k_0\! +\! 1)$, $h(\beta_\varepsilon^{\nu ,2k_0+1,j})\! =\!\beta_{\eta_\varepsilon}^{\nu',2K_0,J_0+j}$ and 
$h(\beta_0^{\nu ,2k_0+1,j})\! =\!\beta_0^{\nu',2K_0,J_0+j}$. Note then that, for $\kappa\!\leq\!\varepsilon\! <\! 2\kappa\! -\! 1$,  
$h(\beta_\varepsilon^{\nu ,2k_0+2,0})\! =\!\beta_{\eta_\varepsilon}^{\nu',2K_0+1,0}$ and 
$h(\beta_0^{\nu ,2k_0+2,0})\! =\!\beta_0^{\nu',2K_0+1,0}$. Then, for $\kappa\! <\!\varepsilon\! <\! 2\kappa\! -\! 1$ and inductively on $0\! <\! j\! <\!\nu (2k_0\! +\! 2)$, $h(\beta_\varepsilon^{\nu ,2k_0+2,j})\! =\!\beta_{\eta_\varepsilon}^{\nu',2K_0+1,j}$ and 
$h(\beta_0^{\nu ,2k_0+2,j})\! =\!\beta_0^{\nu',2K_0+1,j}$. Note then that, for $0\! <\!\varepsilon\! <\!\kappa$,  
$h(\beta_\varepsilon^{\nu ,2k_0+3,0})\! =\!\beta_{\eta_\varepsilon}^{\nu',2K_0+2,0}$ and 
$h(\beta_0^{\nu ,2k_0+3,0})\! =\!\beta_0^{\nu',2K_0+2,0}$. This implies that 
$$2\! +\!\nu (2k_0\! +\! 2)\! -\! 1\! +\!\big(\Sigma_{i<2K_0+1}~\nu'(i)\big)\! +\! 1\! =\! 2\! +\!\Sigma_{i<2K_0+2}~\nu'(i)\mbox{,}$$ 
so that $\nu (2k_0\! +\! 2)\! =\!\nu'(2K_0\! +\! 1)$. Then, for $0\! <\!\varepsilon\! <\!\kappa$ and inductively on 
$0\! <\! j\! <\!\nu (2k_0\! +\! 3)$, $h(\beta_\varepsilon^{\nu ,2k_0+3,j})\! =\!\beta_{\eta_\varepsilon}^{\nu',2K_0+2,j}$ and 
$h(\beta_0^{\nu ,2k_0+3,j})\! =\!\beta_0^{\nu',2K_0+2,j}$. Note then that, for $\kappa\!\leq\!\varepsilon\! <\! 2\kappa\! -\! 1$,  
$h(\beta_\varepsilon^{\nu ,2k_0+4,0})\! =\!\beta_{\eta_\varepsilon}^{\nu',2K_0+3,0}$ and 
$h(\beta_0^{\nu ,2k_0+4,0})\! =\!\beta_0^{\nu',2K_0+3,0}$. This implies that 
$$2\! +\!\nu (2k_0\! +\! 3)\! -\! 1\! +\!\big(\Sigma_{i<2K_0+2}~\nu'(i)\big)\! +\! 1\! =\! 2\! +\!\Sigma_{i<2K_0+3}~\nu'(i)\mbox{,}$$ 
so that $\nu (2k_0\! +\! 3)\! =\!\nu'(2K_0\! +\! 2)$. Inductively, we get $\nu (2k_0\! +\! 2\! +\! n)\! =\!\nu'(2K_0\! +\! 1\! +\! n)$ for each $n\!\in\!\omega$, so that $(\nu ,\nu')\!\in\!\mathbb{E}_t$.\medskip 

 Consider now the sequence $(p_n)_{n\in\omega}$ of prime numbers. We set, for $\alpha\!\in\! 2^\omega$, 
$$S_\alpha\! :=\!\{ p_0^{\alpha (0)+1}\!\cdots\! p_n^{\alpha (n)+1}\mid n\!\in\!\omega\}\mbox{,}$$ 
so that $S_\alpha\!\subseteq\!\omega\!\setminus\!\{ 0\}$ is infinite, and $S_\alpha\cap S_\beta$ is finite if $\alpha\!\not=\!\beta$. Let $\nu_\alpha\!\in\! (\omega\!\setminus\!\{ 0\} )^\omega$ be an injective enumeration of $S_\alpha$. Then 
$(\nu_\alpha ,\nu_\beta)\!\notin\!\mathbb{E}_t$ if $\alpha\!\not=\!\beta$, so that 
$(X_{\nu_\alpha},G_{\nu_\alpha})\not\preceq^i_c(X_{\nu_\beta},G_{\nu_\beta})$. So 
$\big( (X_{\nu_\alpha},G_{\nu_\alpha})\big)_{\alpha\in 2^\omega}$ is an antichain of size $2^{\aleph_0}$ made of minimal elements of $(\mathcal{C}_\kappa ,\preceq^i_c)$, proving that any basis for this class must have size $2^{\aleph_0}$.
\hfill{$\square$}\medskip

\noindent\bf Remark.\rm\ This proof improves the previous proof of Theorem \ref{mainclosed} (b) for $\kappa\! =\! 2$, in the  sense that the $\Sigma_\nu$'s mentioned at the beginning of this section had Cantor-Bendixson rank three, while the $X_\nu$'s here have Cantor-Bendixson rank two.

\section{$\!\!\!\!\!\!$ Some analytic complete sets}\indent

\noindent\bf Notation.\rm\ The space $\mathcal{K}(X)$ of compact subsets of a metrizable compact space $X$, equipped with the Vietoris topology, is a metrizable compact space (see Theorem 4.22 in [K1]). Let $\kappa\! <\!\aleph_0$ (see the remark just before Theorem \ref{mainclosed}). We denote by $\mathfrak{C}_\kappa$ the set of closed graphs on $2^\omega$ having no continuous $\kappa$-coloring, and code $\mathcal{C}_\kappa$ by 
$$C_\kappa\! :=\!\{ (K,G)\!\in\!\mathcal{K}(2^\omega )\!\times\!\mathcal{K}(2^\omega\!\times\! 2^\omega )\mid 
G\mbox{ is a graph on }K\mbox{ having no continuous }\kappa\mbox{-coloring}\} .$$ 
Note that $\mathfrak{C}_0$ is simply the set of closed graphs on $2^\omega$. We then set 
$$Q^\mathfrak{C}_\kappa\! :=\!\{ (G,H)\!\in\!\mathfrak{C}_\kappa^2\mid 
(2^\omega ,G)\preceq^i_c(2^\omega ,H)\}\mbox{,}$$
$E^\mathfrak{C}_\kappa\! :=\!\{ (G,H)\!\in\!\mathfrak{C}_\kappa^2\mid 
(2^\omega ,G)\equiv^i_c(2^\omega ,H)\}\! =\! i(Q^\mathfrak{C}_\kappa )$ (where 
$i(Q)\! :=\! Q\cap Q^{-1}$), 
$$Q^C_\kappa\! :=\!\big\{\big( (K,G),(L,H)\big)\!\in\! C_\kappa^2\mid 
(K,G)\preceq^i_c(L,H)\big\}$$ 
and $E^C_\kappa\! :=\! i(Q^C_\kappa )$. Note that $Q^\mathfrak{C}_\kappa$, 
$Q^C_\kappa$ are quasi-orders, while $E^\mathfrak{C}_\kappa$, 
$E^C_\kappa$ are equivalence relations.\medskip

 Now let $\kappa\!\leq\! 3$ (see Theorem 1.12.(b) in [L]), and $\mathcal{H}(2^\omega )$ be the set of homeomorphisms of $2^\omega$. We equip $\mathcal{H}(2^\omega )$ with the topology whose basic open sets are of the form 
$${O_{U_1,\ldots ,U_n,V_1,\ldots ,V_n}\! :=\!
\{ f\!\in\!\mathcal{H}(2^\omega )\mid\forall 1\!\leq\! i\!\leq\! n~~f[U_i]\! =\! V_i\}}\mbox{,}$$ where $n$ is a natural number and $U_i,V_i$ are clopen subsets of $2^\omega$. By Section 2 in [I-Me], this defines a structure of Polish group on $\mathcal{H}(2^\omega )$. A compatible complete distance is given by
$$d(f,g)\! :=\!\mbox{sup}_{\alpha\in 2^\omega}~
d_{2^\omega}\big( f(\alpha ),g(\alpha )\big)\! +\!
\mbox{sup}_{\alpha\in 2^\omega}~d_{2^\omega}\big( f^{-1}(\alpha ),g^{-1}(\alpha )\big) .$$
We denote by $\mathfrak{H}_\kappa$ the set of homeomorphisms of $2^\omega$ whose induced graph has no continuous $\kappa$-coloring and, as in the introduction of [L], code $\mathcal{H}_\kappa$ by 
$$H_\kappa\! :=\!\{ (K,f)\!\in\!
\mathcal{K}(2^\omega )\!\times\!\mathcal{H}(2^\omega )\mid f[K1]\! =\! K\wedge (K,G_{f_{\vert K}})\mbox{ has no continuous }\kappa\mbox{-coloring}\} .$$
Note that $\mathfrak{H}_0\! =\!\mathcal{H}(2^\omega )$. We then set 
$Q^\mathfrak{H}_\kappa\! :=\!\{ (f,g)\!\in\!\mathfrak{H}_\kappa^2\mid 
(2^\omega ,G_f)\preceq^i_c(2^\omega ,G_g)\}$, 
$E^\mathfrak{H}_\kappa\! :=\! i(Q^\mathfrak{H}_\kappa )$, 
$Q^H_\kappa\! :=\!\big\{\big( (K,f),(L,g)\big)\!\in\! H_\kappa^2\mid 
(K,G_f)\preceq^i_c(L,G_g)\big\}$ and $E^H_\kappa\! :=\! i(Q^H_\kappa )$.\medskip 
 
 Now let $\kappa\!\leq\!\aleph_0$ and $D$ be a countable dense subset of $2^\omega$. We identify 
$\mathcal{P}(D\!\times\! D)$ with $2^{D\times D}$ (equipped with the product topology of the discrete topology on $2$). We denote by $\mathfrak{D}_\kappa$ the set of graphs $G$ on $D$ such that $(2^\omega ,G)$ has no continuous $\kappa$-coloring and set
$$D_\kappa\! :=\!\{ (K,G)\!\in\!\mathcal{K}(2^\omega )\!\times\!\mathcal{P}(D\!\times\! D)\mid 
G\mbox{ is a graph on }K\mbox{ having no continuous }\kappa\mbox{-coloring}\} .$$
Note that $\mathfrak{D}_0$ is simply the set of graphs on $D$. We then set 
$$Q^\mathfrak{D}_\kappa\! :=\!\{ (G,H)\!\in\!\mathfrak{D}_\kappa^2\mid 
(2^\omega ,G)\preceq^i_c(2^\omega ,H)\}\mbox{,}$$ 
$E^\mathfrak{D}_\kappa\! :=\! i(Q^\mathfrak{D}_\kappa )$, 
$Q^D_\kappa\! :=\!\big\{\big( (K,G),(L,H)\big)\!\in\! D_\kappa^2\mid 
(K,G)\preceq^i_c(L,H)\big\}$ and $E^D_\kappa\! :=\! i(Q^D_\kappa )$.

\begin{thm} \label{Cor2} The spaces $\mathfrak{C}_\kappa$, $C_\kappa$,  
$\mathfrak{H}_\kappa$, $H_\kappa$, $\mathfrak{D}_\kappa$ and $D_\kappa$ are Polish, and $FCO$ is Borel reducible to the analytic relations $Q^\mathfrak{C}_\kappa$, 
$E^\mathfrak{C}_\kappa$, $Q^C_\kappa$, $E^C_\kappa$, $Q^\mathfrak{H}_\kappa$, 
$E^\mathfrak{H}_\kappa$, $Q^H_\kappa$, $E^H_\kappa$, $Q^\mathfrak{D}_\kappa$, 
$E^\mathfrak{D}_\kappa$, $Q^D_\kappa$ and $E^D_\kappa$. In particular, these relations are analytic complete as sets.\end{thm}

\noindent\bf Proof.\rm\ Let $G\!\in\!\mathcal{K}(2^\omega\!\times\! 2^\omega )$. Note that $G\!\in\!\mathfrak{C}_\kappa$ if and only if\medskip
 
\leftline{$G\! =\! G^{-1}\wedge G\cap\Delta (2^\omega )\! =\!\emptyset\wedge
\forall (C_i)_{i<\kappa}\!\in\!\big(\borone (2^\omega )\big)^\kappa$}\smallskip

\rightline{$(2^\omega\!\not\subseteq\!\bigcup_{i<\kappa}~C_i)\vee
(\exists i\!\not=\! j\! <\!\kappa ~~C_i\cap C_j\!\not=\!\emptyset )\vee 
(G\cap (\bigcup_{i<\kappa}~C_i^2)\!\not=\!\emptyset).$}\medskip

\noindent Note that, by continuity of the map $(x,y)\!\mapsto\! (y,x)$ and 4.29 in [K1], the condition ``$G\! =\! G^{-1}$" is closed. By definition of the Vietoris topology, the condition ``$G\cap\Delta (2^\omega )\! =\!\emptyset$" is open, while the last condition is closed. Thus $\mathfrak{C}_\kappa$ is a difference of two open sets, and thus $\bormtwo$ in 
$\mathcal{K}(2^\omega\!\times\! 2^\omega )$, and Polish. Note that 
$C_\kappa\! =\!
\{ (K,G)\!\in\!\mathcal{K}(2^\omega )\!\times\!\mathcal{K}(2^\omega\!\times\! 2^\omega )\mid G\!\in\!\mathfrak{C}_\kappa\wedge G\!\subseteq\! K^2\}$ by Theorem 2.2.1 in [E]. By 4.29 in [K1], $C_\kappa$ is also Polish. By Theorem 12.5 in [L], the set 
$\mathfrak{H}_\kappa$ is $\bormtwo$ in $\mathcal{H}(2^\omega )$, and Polish. By Theorem 1.14 in [L], the set $H_\kappa$ is $\bormtwo$ in 
$\mathcal{K}(2^\omega )\!\times\!\mathcal{H}(2^\omega )$, and Polish. Let 
$G\!\in\!\mathcal{P}(D\!\times\! D)$. Note that 
$G\!\in\!\mathfrak{D}_\kappa$ if and only if the formula above holds. We enumerate 
$D\! :=\!\{ d_n\mid n\!\in\!\omega\}$ injectively. The condition ``$G\! =\! G^{-1}$" can be written ``$\forall m,n\!\in\!\omega ~~(d_m,d_n)\!\notin\! G\vee (d_n,d_m)\!\in\! G$", which is a closed condition. The condition ``$G\cap\Delta (2^\omega )\! =\!\emptyset$" can be written ``$\forall n\!\in\!\omega ~~(d_n,d_n)\!\notin\! G$", which is a closed condition. For the last condition, note that $\borone (2^\omega )$ is countable. If $\kappa$ is finite, then the condition 
``$G\cap (\bigcup_{i<\kappa}~C_i^2)\!\not=\!\emptyset$" can be written 
``$\exists m,n\!\in\!\omega ~~\exists i\! <\!\kappa ~~(q_m,q_n)\!\in\! G\cap C_i^2$", so that $\mathfrak{D}_\kappa$ is 
$\bormtwo$ in $\mathcal{P}(D\!\times\! D)$, and Polish. If $\kappa\! =\!\aleph_0$, then by compactness the last condition can be written 
$$\forall n\!\in\!\omega ~~\forall (C_i)_{i<n}\!\in\!\big(\borone (2^\omega )\big)^n~~
(2^\omega\!\not\subseteq\!\bigcup_{i<n}~C_i)\vee
(\exists i\!\not=\! j\! <\! n~~C_i\cap C_j\!\not=\!\emptyset )\vee 
(G\cap (\bigcup_{i<n}~C_i^2)\!\not=\!\emptyset)\mbox{,}$$
so that $\mathfrak{D}_\kappa$ is Polish again. Note that 
$D_\kappa\! =\!\{ (K,G)\!\in\! \mathcal{K}(2^\omega )\!\times\!\mathcal{P}(D\!\times\! D)\mid 
G\!\in\!\mathfrak{D}_\kappa\wedge G\!\subseteq\! K^2\}$ by Theorem 2.2.1 in [E]. The condition ``$G\!\subseteq\! K^2$" is ``$\forall m,n\!\in\!\omega ~~(d_m,d_n)\!\notin\! G\vee d_m,d_n\!\in\! K$", which is a closed condition, so that 
$D_\kappa$ is Polish.\medskip

 Recall that 
$$(2^\omega ,G)\preceq^i_c(2^\omega ,H)\Leftrightarrow
\exists\varphi\! :\! 2^\omega\!\rightarrow\! 2^\omega\mbox{ injective continuous with }
G\!\subseteq\! (\varphi\!\times\!\varphi )^{-1}(H).$$
Note that $\varphi\! :\! 2^\omega\!\rightarrow\! 2^\omega$ is injective if and only if 
$\varphi [O\cap U]\! =\!\varphi[O]\cap\varphi [U]$ whenever $O,U$ are clopen subsets of $2^\omega$. By Lemma 12.4 in [L], and [K, 4.19, 4.29, 27.7], 
$$\{ (G,H)\!\in\!\mathcal{K}(2^\omega\!\times\! 2^\omega )^2\mid (2^\omega ,G)\preceq_c^i(2^\omega ,H)\}$$ 
is analytic, and thus $Q^\mathfrak{C}_\kappa$ and $E^\mathfrak{C}_\kappa$ are analytic.\medskip

 If now $(K,G)$, $(L,H)\!\in\!\mathcal{K}(2^\omega )\!\times\!
\mathcal{K}(2^\omega\!\times\! 2^\omega )$, then, since $K$ is a retract of $2^\omega$ by 2.8 in [K1], 
$$(K,G)\preceq^i_c(L,H)\Leftrightarrow
\exists\psi\!\in\!\mathcal{C}(2^\omega ,2^\omega )~~\psi [K1]\!\subseteq\! L\wedge
\psi_{\vert K}\mbox{ is injective}\wedge (\psi\!\times\!\psi )[G]\!\subseteq\! H.$$
By Lemma 12.4 in [L] and 4.29 in [K1], 
$$\big\{\big( (K,G),(L,H)\big)\!\in\!\big(\mathcal{K}(2^\omega )\!\times\!\mathcal{K}(2^\omega\!\times\! 2^\omega )\big)^2\mid (K,G)\preceq^i_c(L,H)\big\}\mbox{,}$$
$Q^C_\kappa$ and $E^C_\kappa$ are analytic.

\vfill\eject

 Let $f,g\!\in\!\mathcal{H}(2^\omega )$. Then\medskip
 
\leftline{$(2^\omega ,G_f)\preceq^i_c(2^\omega ,G_g)\Leftrightarrow\exists\varphi\! :\! 2^\omega\!\rightarrow\! 2^\omega\mbox{ injective continuous }$}\smallskip

\rightline{$\forall\alpha\!\in\! 2^\omega ~~f(\alpha )\! =\!\alpha\vee\varphi\big( f(\alpha )\big)\! =\! g\big(\varphi (\alpha )\big)\vee\varphi\big( f(\alpha )\big)\! =\! g^{-1}\big(\varphi (\alpha )\big) .$}\medskip

\noindent This shows that $\{ (f,g)\!\in\!\mathcal{H}(2^\omega )^2\mid 
(2^\omega ,G_f)\preceq^i_c(2^\omega ,G_g)\}$ and thus $Q^\mathfrak{H}_\kappa$ and $E^\mathfrak{H}_\kappa$ are analytic. The previous discussions show that $Q^H_\kappa$ and $E^H_\kappa$ are also analytic. If $G,H\!\in\!\mathcal{P}(D\!\times\! D)$, then the condition ``$G\!\subseteq\! (\varphi\!\times\!\varphi )^{-1}(H)$" can be written\medskip
 
\centerline{``$\forall m,n\!\in\!\omega ~~(d_m,d_n)\!\notin\! G\vee
\big(\varphi (d_m),\varphi (d_n)\big)\!\in\! H$",}\medskip

\noindent which is a closed condition, proving that 
$\{ (G,H)\!\in\!\big(\mathcal{P}(D\!\times\! D)\big)^2\mid (2^\omega ,G)\preceq^i_c(2^\omega ,H)\}$, 
$Q^\mathfrak{D}_\kappa$ and $E^\mathfrak{D}_\kappa$ are analytic. The previous discussions show that 
$Q^D_\kappa$ and $E^D_\kappa$ are also analytic.\medskip

 We define a map 
$\mathfrak{g}\! :\!\mathbb{M}\!\rightarrow\!\mathcal{K}(2^\omega\!\times\! 2^\omega )$ by 
$\mathfrak{g}(f)\! :=\! G_f$. Let $O$ be an open subset of $2^\omega\!\times\! 2^\omega$, and $(C^0_n)_{n\in\omega}$, $(C^1_n)_{n\in\omega}$ be sequences of clopen subsets of 
$2^\omega$ with the property that $O\! =\!\bigcup_{n\in\omega}~(C^0_n\!\times\! C^1_n)$. If $f\!\in\!\mathbb{M}$ and $G_f\!\subseteq\! O$, then there is a finite subset $F$ of 
$\omega$ with $G_f\! =\! s\big(\textup{Graph}(f)\big)\!\subseteq\!
\bigcup_{n\in F}~(C^0_n\!\times\! C^1_n)$. Note then that 
$$\bigcup_{n\in F}~(C^0_n\!\times\! C^1_n)\! =\!\bigcup_{S\subseteq F}~
\big( (\bigcap_{n\in S}~C^0_n\cap\bigcap_{n\in F\setminus S}~2^\omega\!\setminus\! C^0_n)\!\times\! (\bigcup_{n\in S}~C^1_n)\big) .$$ 
Thus 
$$\begin{array}{ll}
\textup{Graph}(f)\!\subseteq\!\bigcup_{n\in F}~(C^0_n\!\times\! C^1_n)\!\!\!\!
& \Leftrightarrow\forall S\!\subseteq\! F~
f[\bigcap_{n\in S}~C^0_n\cap\bigcap_{n\in F\setminus S}~
2^\omega\!\setminus\! C^0_n]\!\subseteq\!\bigcup_{n\in S}~C^1_n\cr 
& \Leftrightarrow\forall S\!\subseteq\! F~\exists R_n\!\in\!\borone (2^\omega)\cr
& \hfill{f[\bigcap_{n\in S}~C^0_n\cap
\bigcap_{n\in F\setminus S}~2^\omega\!\setminus\! C^0_n]\! =\! 
R_n\!\subseteq\!\bigcup_{n\in S}~C^1_n.}
\end{array}$$
This implies that $\{ f\!\in\!\mathbb{M}\mid G_f\!\subseteq\! O\}$ is an open subset of 
$\mathbb{M}$ since 
$$G_f\!\subseteq\! O\Leftrightarrow
\exists F\!\subseteq\!\omega\mbox{ finite with }
\textup{Graph}(f)\!\subseteq\!\bigcap_{\varepsilon\in 2}
\big(\bigcup_{n\in F}~(C^\varepsilon_n\!\times\! C^{1-\varepsilon}_n)\big) .$$ 
Now $G_f\cap O\!\not=\!\emptyset\Leftrightarrow
\exists n\!\in\!\omega ~~\exists\varepsilon\!\in\! 2~~
C^\varepsilon_n\cap f^{-1}(C^{1-\varepsilon}_n)\!\not=\!\emptyset\Leftrightarrow
\exists n\!\in\!\omega ~~\exists\varepsilon\!\in\! 2~~\exists\alpha\!\in\! C^\varepsilon_n~~
f(\alpha )\!\in\! C^{1-\varepsilon}_n$, so that 
$\{ f\!\in\!\mathbb{M}\mid G_f\cap O\!\not=\!\emptyset\}$ is an open subset of 
$\mathbb{M}$. Thus $\mathfrak{g}$ is continuous.\medskip
 
  By Lemma 7.11 in [L], if $f,g\!\in\!\mathbb{M}$, then $(f,g)\!\in\! FCO$ if and only if 
$(2^\omega ,G_f)\!\preceq^i_c(2^\omega ,G_g)$. As $FCO$ is symmetric, 
$(f,g)\!\in\! FCO$ if and only if $(2^\omega ,G_f)\!\equiv^i_c(2^\omega ,G_g)$. We define 
$\phi_1(f)\! :\! N_1\!\rightarrow\! N_1$ by $\phi_1(f)(1\alpha )\! :=\! 1f(\alpha )$, so that 
$\phi_1(f)$ is a homeomorphism with infinite orbits and $G_{\phi_1(f)}$ is a closed graph on $N_1$. We define a map 
$\mathfrak{g}^+\! :\!\mathbb{M}\!\rightarrow\!\mathfrak{C}_\kappa$ by 
$$\mathfrak{g}^+(f)\! :=\! G_{\phi_1(f)}\cup
\{ (0^{m+1}1^\infty ,0^{n+1}1^\infty )\mid m\!\not=\! n\!\in\!\kappa\! +\! 1\} .$$ 
Note that $G_{\phi_1(f)}\! =\!\{ (1\alpha ,1\beta )\mid (\alpha ,\beta )\!\in\! G(f)\}$, so that 
$\mathfrak{g}^+$ is continuous by 4.29 (iv and vi) in [K1]. In order to prove that $FCO$ is Borel reducible to 
$Q^\mathfrak{C}_\kappa$ and $E^\mathfrak{C}_\kappa$, it is enough to prove that if 
$f,g\!\in\!\mathbb{M}$, then $(2^\omega ,G_f)\!\preceq^i_c(2^\omega ,G_g)$ if and only if 
$\big( 2^\omega ,\mathfrak{g}^+(f)\big)\!\preceq^i_c
\big( 2^\omega ,\mathfrak{g}^+(g)\big)$. So let 
$\varphi\! :\! 2^\omega\!\rightarrow\! 2^\omega$ injective continuous with 
$G_f\!\subseteq\! (\varphi\!\times\!\varphi )^{-1}(G_g)$. We define 
$\Phi\! :\! 2^\omega\!\rightarrow\! 2^\omega$ by $\Phi (0\alpha )\! :=\! 0\alpha$ and 
$\Phi (1\alpha )\! :=\! 1\varphi (\alpha )$, so that $\Phi$ is injective continuous.

\vfill\eject

 Moreover, 
$\big(\Phi (0^{m+1}1^\infty ),\Phi (0^{n+1}1^\infty )\big)\! =\! (0^{m+1}1^\infty ,0^{n+1}1^\infty )\!\in\!\mathfrak{g}^+(g)$ and 
$$\Big(\Phi (1\alpha ),\Phi\big(\phi_1(f)(1\alpha )\big)\Big)\! =\!
\Big( 1\varphi (\alpha ),1\varphi\big( f(\alpha )\big)\Big)\! =\!
\big( 1\beta , 1g^{\pm 1}\big(\beta )\big)\!\in\!\mathfrak{g}^+(g)\mbox{,}$$ 
showing that $\Phi$ is a witness for the fact that 
$\big( 2^\omega ,\mathfrak{g}^+(f)\big)\preceq^i_c\big( 2^\omega ,\mathfrak{g}^+(g)\big)$.\medskip

 Conversely, let $\Phi\! :\! 2^\omega\!\rightarrow\! 2^\omega$ injective continuous with 
$\mathfrak{g}^+(f)\!\subseteq\! (\Phi\!\times\!\Phi )^{-1}\big(\mathfrak{g}^+(g)\big)$. Note that the $\mathfrak{g}^+(f)$-connected component of $N_0$ has size $\kappa\! +\! 1$, while each $\mathfrak{g}^+(f)$-connected component of $N_1$ is infinite. This implies that 
$\Phi (1\alpha )(0)\! =\! 1$. We now define $\varphi\! :\! 2^\omega\!\rightarrow\! 2^\omega$ by the formula 
$$\varphi (\alpha )\! :=\!\Phi (1\alpha )^-\! :=\!
\big(\Phi (1\alpha )(1),\Phi (1\alpha )(2),\cdots\big)\mbox{,}$$ 
so that $\varphi$ is injective continuous. Moreover, 
$$\Big( 1\varphi (\alpha ),1\varphi\big( f(\alpha )\big)\Big)\! =\!
\Big(\Phi (1\alpha ),\Phi\big( 1f(\alpha )\big)\Big)\! =\!
\Big(\Phi (1\alpha ),\Phi\big(\phi_1(f)(1\alpha )\big)\Big)\!\in\!\mathfrak{g}^+(g)\mbox{,}$$
so that $\varphi\big( f(\alpha )\big)\! =\! g^{\pm 1}\big(\varphi (\alpha )\big)$ and 
$\Big(\varphi (\alpha ),\varphi\big( f(\alpha )\big)\Big)\!\in\! G_g$, showing that $\varphi$ is a witness for the fact that $(2^\omega ,G_f)\preceq^i_c(2^\omega ,G_g)$. Thus $FCO$ is Borel (in fact continuously) reducible to $Q^\mathfrak{C}_\kappa$ and 
$E^\mathfrak{C}_\kappa$. Now note that the map 
$i_C\! :\!\mathfrak{C}_\kappa\!\rightarrow\! C_\kappa$ defined by 
$i_C(G)\! :=\! (2^\omega ,G)$ is continuous, 
$Q^\mathfrak{C}_\kappa\! =\! (i_C\!\times\! i_C)^{-1}(Q^C_\kappa )$ and 
$E^\mathfrak{C}_\kappa\! =\! (i_C\!\times\! i_C)^{-1}(E^C_\kappa )$, so that $FCO$ is also  Borel (in fact continuously) reducible to $Q^C_\kappa$ and $E^C_\kappa$.\medskip

 Note then that $\mathfrak{H}_3\!\subseteq\!\mathfrak{H}_2\!\subseteq\!\mathfrak{H}_1\!\subseteq\!\mathfrak{H}_0\! =\!\mathcal{H}(2^\omega )$. We define 
$h_3\! :\! 2^\omega\!\rightarrow\! 2^\omega$ by $h_3(0^\infty )\! :=\! 0^\infty$ and 
$h_3(0^{2n+\varepsilon}1\alpha )\! :=\! 0^{2n+1-\varepsilon}1\alpha$, so that $h_3\!\in\!\mathcal{H}(2^\omega )$ has orbits of size at most two and $F_{h_3}\! =\{ 0^\infty\}$ is not open in $2^\omega$. By Proposition \ref{fp}, there is no continuous $\aleph_0$-coloring of $G_{h_3}$, so that $h_3\!\in\!\mathfrak{H}_3$. We define, for 
$f\!\in\!\mathbb{M}$, $\phi (f)\!\in\!\mathcal{H}(2^\omega )$ by $\phi (f)(0\alpha )\! :=\! 0h_3(\alpha )$ and 
$\phi (f)(1\alpha )\! :=\! 1f(\alpha )$. Note that $\phi (f)\!\in\!\mathfrak{H}_3$ and  
$\phi\! :\!\mathbb{M}\!\rightarrow\!\mathfrak{H}_3$ is continuous (consider the distance 
$d$). In order to prove that $FCO$ is Borel reducible to $Q^\mathfrak{H}_3$ and 
$E^\mathfrak{H}_3$, it is enough to prove that if $f,g\!\in\!\mathbb{M}$, then 
$(2^\omega ,G_f)\preceq^i_c(2^\omega ,G_g)$ if and only if 
$(2^\omega ,G_{\phi (f)})\preceq^i_c(2^\omega ,G_{\phi (g)})$. We argue essentially as above, using the facts that the 
$\phi (f)$-orbit of $0\alpha$ has size at most two like the $h_3$-orbit of $\alpha$, while the $\phi (f)$-orbit of $1\alpha$ is infinite like the $f$-orbit of $\alpha$. Thus $FCO$ is Borel (in fact continuously) reducible to $Q^\mathfrak{H}_3$ and $E^\mathfrak{H}_3$, and in fact to $Q^\mathfrak{H}_\kappa$ and $E^\mathfrak{H}_\kappa$ because of the inclusions above. Now note that the map $i_H\! :\!\mathfrak{H}_\kappa\!\rightarrow\! H_\kappa$ defined by 
$i_H(f)\! :=\! (2^\omega ,f)$ is continuous, 
$Q^\mathfrak{H}_\kappa\! =\! (i_H\!\times\! i_H)^{-1}(Q^H_\kappa )$ and 
$E^\mathfrak{H}_\kappa\! =\! (i_H\!\times\! i_H)^{-1}(E^H_\kappa )$, so that $FCO$ is also  Borel (in fact continuously) reducible to $Q^H_\kappa$ and $E^H_\kappa$.\medskip 

 By Corollary 5.10 in [L], if $f,g\!\in\!\mathbb{M}$, then $(f,g)\!\in\! FCO$ if and only if 
$$(\overline{\mbox{proj}[\mathbb{G}_f]},\mathbb{G}_f)\!\preceq^i_c
(\overline{\mbox{proj}[\mathbb{G}_g]},\mathbb{G}_g)$$ 
(see Section 5 in [L] for the definition of the graph $\mathbb{G}_f$, whose vertices have degree at most one). This definition, as well as the notation before Theorem 13.2 in [L], show that $\mbox{proj}[\mathbb{G}_f]$ is contained in the copy 
$\mathcal{K}_{2^\infty}\! :=\! (2\cup\{ c,a,\overline{a}\} )^\omega$ of $2^\omega$. In fact, the definition of $\mathbb{G}_f$ shows that $\mbox{proj}[\mathbb{G}_f]$ is in fact contained in the closed nowhere dense subset 
$\{ x\!\in\!\mathcal{K}_{2^\infty}\mid
\forall m\!\in\!\omega ~~x(m)\! =\! c\vee\forall n\!\geq\! m~~x(n)\!\not=\! c\}$ of 
$\mathcal{K}_{2^\infty}$. In particular, $\overline{\mbox{proj}[\mathbb{G}_f]}$ is nowhere dense in $\mathcal{K}_{2^\infty}$.\medskip

\noindent\bf Claim 1.\it\ (Ryll-Nardzewski) Let $P,Q$ be closed nowhere dense subsets of $2^\omega$, and $\varphi\! :\! P\!\rightarrow\! Q$ be a continuous injection. Then there is a homeomorphism $\varphi^*$ of $2^\omega$ such that 
$\varphi^*(\alpha )\! =\!\varphi (\alpha )$ if $\alpha\!\in\! P$.\rm\medskip

 Indeed, the compact subset $R\! :=\!\varphi [P]$ of $Q$ is also closed and nowhere dense in $2^\omega$, and the map $\varphi'\! :\! P\!\rightarrow\! R$ defined by 
 $\varphi'(\alpha )\! :=\!\varphi (\alpha )$ is a homeomorphism. The Ryll-Nardzewski theorem (see Corollary 2 in [Kn-R]) provides a homeomorphism $\varphi^*$ of 
$2^\omega$ extending $\varphi'$, and thus having the desired property.\hfill{$\diamond$}

\vfill\eject\medskip

 Claim 1 implies that if $f,g\!\in\!\mathbb{M}$, then $(f,g)\!\in\! FCO$ if and only if 
$(\mathcal{K}_{2^\infty},\mathbb{G}_f)\!\preceq^i_c
(\mathcal{K}_{2^\infty},\mathbb{G}_g)$. Indeed, assume that $\varphi\! :\!\overline{\mbox{proj}[\mathbb{G}_f]}\!\rightarrow\!\overline{\mbox{proj}[\mathbb{G}_g]}$ is injective continuous and 
$\mathbb{G}_f\!\subseteq\! (\varphi\!\times\!\varphi )^{-1}(\mathbb{G}_g)$. Claim 1 provides $\varphi^*\! :\!\! :\!\mathcal{K}_{2^\infty}\!\rightarrow\!\mathcal{K}_{2^\infty}$ injective continuous coinciding with $\varphi$ on 
$\overline{\mbox{proj}[\mathbb{G}_f]}$, which is a witness for the fact that 
$(\mathcal{K}_{2^\infty},\mathbb{G}_f)\!\preceq^i_c
(\mathcal{K}_{2^\infty},\mathbb{G}_g)$. Conversely, if 
$\Phi\! :\!\! :\!\mathcal{K}_{2^\infty}\!\rightarrow\!\mathcal{K}_{2^\infty}$ is injective continuous and $\mathbb{G}_f$ is contained in $(\Phi\!\times\!\Phi )^{-1}(\mathbb{G}_g)$, then $\mbox{proj}[\mathbb{G}_f]$ and thus $\overline{\mbox{proj}[\mathbb{G}_f]}$ are contained in $\Phi^{-1}(\overline{\mbox{proj}[\mathbb{G}_g]})$, so that the map 
$\varphi\! :=\!\Phi_{\vert\overline{\mbox{proj}[\mathbb{G}_f]}}$ is a witness for the fact that $(\overline{\mbox{proj}[\mathbb{G}_f]},\mathbb{G}_f)\!\preceq^i_c
(\overline{\mbox{proj}[\mathbb{G}_g]},\mathbb{G}_g)$. As $FCO$ is symmetric, 
$(f,g)\!\in\! FCO$ if and only if 
$(\mathcal{K}_{2^\infty},\mathbb{G}_f)\!\equiv^i_c(\mathcal{K}_{2^\infty},\mathbb{G}_g)$.
\medskip

 Let $i\! :\!\mathcal{K}_{2^\infty}\!\rightarrow\! N_1\!\subseteq\! 2^\omega$ be a homeomorphism. The definition of $\mathbb{G}_f$ shows that it is contained in the countable dense subset $\mathbb{Q}\! :=\!\big\{ x\!\in\!\mathcal{K}_{2^\infty}\mid \exists l\!\in\!\omega ~\exists\varepsilon\!\in\!\{ a,\overline{a}\} ~\forall k\!\geq\! l~\,x(k)\! =\!\varepsilon\big\}$ of $\mathcal{K}_{2^\infty}$. In particular $Q\! :=\! i[\mathbb{Q}]$ is a countable dense subset of $N_1$, as well as $N_1\cap D$.\medskip

\noindent\bf Claim 2.\it\ (van Engelen) Let $Q,D$ be countable dense subsets of 
$2^\omega$. Then there is $h\!\in\!\mathcal{H}(2^\omega )$ such that $h[Q]\! =\! D$.\rm\medskip

 Indeed, we enumerate injectively $Q\! =\!\{ q_i\mid i\!\in\!\omega\}$ and 
$D\! =\!\{ d_i\mid i\!\in\!\omega\}$. Note that $\{ q_0\}$ is a zero-dimensional space homeomorphic to any of its nonempty clopen subsets, $\{ q_i\}$ (resp., $\{ d_i\}$) is closed nowhere dense in $Q$ (resp., $D$), and $\{ q_i\}\!\approx\!\{ q_0\}\!\approx\!\{ d_i\}$ for each $i$. Theorem 3.2.6 in [vE] provides the desired homeomorphism.\hfill{$\diamond$}\medskip

 Claim 2 provides $h\!\in\!\mathcal{H}(N_1)$ such that $h[Q]\! =\! N_1\cap D$. We set $H\! :=\! h\!\circ\! i$, so that 
$H\! :\!\mathcal{K}_{2^\infty}\!\rightarrow\! N_1$ is a homeomorphism. Recall that the {\bf chromatic number} of a graph $(X,G)$ is the smallest cardinal $\kappa$ for which there is a $\kappa$-coloring of $(X,G)$.\medskip

\noindent\bf Claim 3.\it\ There is a sequence $\big( (F_n,G_n)\big)_{n\in\omega}$ made of finite connected graphs which are pairwise $\preceq^i$-incomparable, have pairwise different chromatic numbers, and whose vertices have degree at least two.\rm\medskip

 Indeed, we use the {\bf Kneser graphs} $K(n,k)$. Recall that, if 
$n,k\!\in\!\omega\!\setminus\!\{ 0\}$, then $K(n,k)$ has set of vertices $[n]^k$, and 
$A,B\!\in\! [n]^k$ are $K(n,k)$-related if $A\cap B\! =\!\emptyset$. If $n\!\geq\! 3k$, then $K(n,k)$ is finite connected and its vertices have degree at least two. Note that 
$\mbox{Cardinality}([n]^k)\! =\!\binom{n}{k}$ and, by Theorem 6.29 in [H-N], $K(n,k)$ has chromatic number $n\! -\! 2k\! +\! 2$ if $n\!\geq\! 2k$. Moreover, by Proposition 6.27 in 
[H-N], $K(n,k)\not\preceq^i K(n',k')$ if $2\!\leq\!\frac{n'}{k'}\! <\!\frac{n}{k}$ (even without necessarily injectivity). All this implies that it is enough to construct a sequence 
$\big( (n_p,k_p)\big)_{p\in\omega}$ of pairs of positive natural numbers satisfying the following. 
$$\begin{array}{ll}
& (1)~3\!\leq\!\frac{n_{p+1}}{k_{p+1}}\! <\!\frac{n_p}{k_p}\cr
& (2)~\binom{n_p}{k_p}\! <\!\binom{n_{p+1}}{k_{p+1}}\cr
& (3)~(n_p\! -\! 2k_p\! +\! 2)_{p\in\omega}\mbox{ is injective}
\end{array}$$
We set $n_p\! :=\! 3\!\cdot\! 2^p\! +\! 1$ and $k_p\! :=\! 2^p$, so that (1) and (3) are satisfied. For (2), note that 
$$\binom{n_p}{k_p}\! =\!\frac{n_p!}{k_p!(n_p-k_p)!}\! =\!
\frac{(3\cdot 2^p+1)!}{(2^p)!(2\cdot 2^p+1)!}\mbox{,}$$ 
thus $\binom{n_p}{k_p}\! <\!\binom{n_{p+1}}{k_{p+1}}\Leftrightarrow
\frac{(3\cdot 2^p+1)!}{(2^p)!(2\cdot 2^p+1)!}\! <\!
\frac{(6\cdot 2^p+1)!}{(2\cdot 2^p)!(4\cdot 2^p+1)!}\Leftrightarrow 
(3\!\cdot\! 2^p\! +\! 1)!\!\cdot\! (4\!\cdot\! 2^p\! +\! 1)!\! <\! 
(6\!\cdot\! 2^p\! +\! 1)!\!\cdot\! (2^p)!\!\cdot\! (2\!\cdot\! 2^p\! +\! 1)$, 
$\binom{n_p}{k_p}\! <\!\binom{n_{p+1}}{k_{p+1}}\Leftrightarrow (3\!\cdot\! 2^p\! +\! 1)!\! <\! 
(6\!\cdot\! 2^p\! +\! 1)\!\cdots\! (4\!\cdot\! 2^p\! +\! 2)\!\cdot\! (2\!\cdot\! 2^p\! +\! 1)\!\cdot\! (2^p)!$. This holds since there are $2\!\cdot\! 2^p$ factors in 
$(6\!\cdot\! 2^p\! +\! 1)\!\cdots\! (4\!\cdot\! 2^p\! +\! 2)\! >\! (3\!\cdot\! 2^p\! +\! 1)\!\cdots\! (2^p\! +\! 2)$ and $2\!\cdot\! 2^p\! +\! 1\! >\! 2^p\! +\! 1$.\hfill{$\diamond$}\medskip

 We may assume that $F_n\!\subseteq\! N_{0^{n+1}1}\cap D$. Claim 3 allows us to define 
$\mathcal{G}\! :\!\mathbb{M}\!\rightarrow\!\mathfrak{D}_{\aleph_0}$, by 
$$\mathcal{G}(f)\! :=\! (H\!\times\! H)[\mathbb{G}_f]\cup\bigcup_{n\in\omega}~G_n.$$ 
This definition is correct since $\big( 2^\omega ,\mathcal{G}(f)\big)$ has no continuous $\aleph_0$-coloring. Indeed, we argue by contradiction to see that. By compactness of $2^\omega$, this graph would have a continuous $\kappa$-coloring for some $\kappa\!\in\!\omega$, which is not the case because the $G_n$'s have pairwise different chromatic numbers. The beginning of the proof of Theorem 13.2 in [L] shows that $\mathcal{G}$ is Borel. Let 
$f,g\!\in\!\mathbb{M}$. If $(f,g)\!\in\! FCO$, then 
$(\mathcal{K}_{2^\infty},\mathbb{G}_f)\!\preceq^i_c(\mathcal{K}_{2^\infty},\mathbb{G}_g)$, with $\varphi$ as a witness. We define ${\Phi\! :\! 2^\omega\!\rightarrow\! 2^\omega}$ by $\Phi (0\alpha )\! :=\! 0\alpha$ and 
$\Phi (1\alpha )\! :=\! H\Big(\varphi\big( H^{-1}(1\alpha )\big)\Big)$, so that $\Phi$ witnesses  
$\big( 2^\omega ,\mathcal{G}(f)\big)\preceq^i_c\big( 2^\omega ,\mathcal{G}(g)\big)$. Conversely, assume that 
${\Phi\! :\! 2^\omega\!\rightarrow\! 2^\omega}$ is injective continuous and 
$\mathcal{G}(f)\!\subseteq\! (\Phi\!\times\!\Phi )^{-1}(\mathcal{G}(g))$. The vertices in $N_1$ have degree at least two, while the other vertices of $\mathcal{G}(f)$ have degree at most one like those of $\mathbb{G}_f$. This implies that 
$\Phi\!\times\!\Phi$ sends $\bigcup_{n\in\omega}~G_n$ into itself, by injectivity. As the $G_n$'s are connected and pairwise $\preceq^i$-incomparable, $\Phi\!\times\!\Phi$ sends $G_n$ into itself, and onto itself by finiteness. Therefore 
$\Phi\!\times\!\Phi$ sends $(H\!\times\! H)[\mathbb{G}_f]$ into $(H\!\times\! H)[\mathbb{G}_g]$, by injectivity. Thus $\Phi$ sends $\mbox{proj}\big[ (H\!\times\! H)[\mathbb{G}_f]\big]$ into $\mbox{proj}\big[ (H\!\times\! H)[\mathbb{G}_g]\big]$, and 
$\overline{\mbox{proj}\big[ (H\!\times\! H)[\mathbb{G}_f]\big]}$ into 
$\overline{\mbox{proj}\big[ (H\!\times\! H)[\mathbb{G}_g]\big]}$. As $H$ is a homeomorphism, 
$$\overline{\mbox{proj}\big[ (H\!\times\! H)[\mathbb{G}_f]\big]}\! =\! 
H\big[\overline{\mbox{proj}[\mathbb{G}_f]}\big]\mbox{,}$$ 
and similarly with $g$. This allows us to define 
$\phi\! :\!\overline{\mbox{proj}[\mathbb{G}_f]}\!\rightarrow\!\overline{\mbox{proj}[\mathbb{G}_g]}$ by 
$\phi (\alpha )\! :=\! H^{-1}\Big(\Phi\big( H(\alpha )\big)\Big)$, and $\phi$ is a witness for the fact that 
$(\overline{\mbox{proj}[\mathbb{G}_f]},\mathbb{G}_f)\!\preceq^i_c
(\overline{\mbox{proj}[\mathbb{G}_g]},\mathbb{G}_g)$. Thus $(f,g)\!\in\! FCO$ and  
$\mathcal{G}$ Borel reduces $FCO$ to $Q^\mathfrak{D}_{\aleph_0}$ and $E^\mathfrak{D}_{\aleph_0}$. As 
$\mathfrak{D}_{\aleph_0}\!\subseteq\!\mathfrak{D}_\kappa$, this also holds for $Q^\mathfrak{D}_\kappa$ and 
$E^\mathfrak{D}_\kappa$ if $\kappa\! <\!\aleph_0$. Now note that the map 
$i_D\! :\!\mathfrak{D}_\kappa\!\rightarrow\! D_\kappa$ defined by $i_D(G)\! :=\! (2^\omega ,G)$ is continuous, 
$Q^\mathfrak{D}_\kappa\! =\! (i_D\!\times\! i_D)^{-1}(Q^D_\kappa )$ and 
$E^\mathfrak{D}_\kappa\! =\! (i_D\!\times\! i_D)^{-1}(E^D_\kappa )$, so that $FCO$ is also  Borel reducible to 
$Q^D_\kappa$ and $E^D_\kappa$.\medskip
 
 [De-GR-Ka-Kun-Kw] shows that $FCO$ is analytic complete as a set. Thus our sets are Borel analytic complete (using pre-images by Borel functions). By [K2], our sets are analytic complete.\hfill{$\square$}\medskip
 
\noindent\bf Question.\rm\ [De-GR-Ka-Kun-Kw] shows that $FCO$ is analytic complete as a set. On the other hand, Theorem 5 in [Ca-G] shows that the conjugacy relation on $\mathcal{H}(2^\omega )$ is Borel-bi-reducible with the most complicated of the orbit equivalence relations induced by a Borel action of the group of bijections of $\omega$. Also, in [Lo-R] it is proved that the bi-homomorphism relation between countable graphs is analytic complete as an equivalence relation. So we can ask about the position of the equivalence relations mentioned in Theorem \ref{Cor2} among analytic equivalence relations, in particular\medskip

\noindent (1) $E^\mathfrak{C}_0\! :=\!\{ (G,H)\!\in\!\mathcal{K}(2^\omega\!\times\! 2^\omega )^2\mid 
G,H\mbox{ are graphs}\wedge (2^\omega ,G)\equiv^i_c(2^\omega ,H)\}$,\smallskip

\noindent (2) $E^\mathfrak{H}_0\! :=\!\{ (f,g)\!\in\!\mathcal{H}(2^\omega )^2\mid 
(2^\omega ,G_f)\equiv^i_c(2^\omega ,G_g)\}$,\smallskip

\noindent (3) $E^\mathfrak{D}_0\! :=\!\{ (G,H)\!\in\!\mathcal{P}(D\!\times\! D)^2\mid 
G,H\mbox{ are graphs}\wedge (2^\omega ,G)\equiv^i_c(2^\omega ,H)\}$, where $D$ is a countable dense subset of 
$2^\omega$.

\vfill\eject

\section{$\!\!\!\!\!\!$ References}\indent

\noindent [Ca-G]\ \ R. Camerlo and S. Gao, The completeness of the isomorphism relation for countable Boolean algebras,\ \it Trans. Amer. Math. Soc.\rm\ 353, 2 (2001), 491-518

\noindent [C-D-To-W]\ \ D. Cenzer, A. Dashti, F. Toska and S. Wyman, Computability of countable subshifts in one dimension,\ \it Theory Comput. Syst.\rm\ 51 (2012), 352-371

\noindent [Co-M]\ \ C. T. Conley and B. D. Miller, An antibasis result for graphs of infinite Borel chromatic number,\ \it Proc. Amer. Math. Soc.\rm 142, 6 (2014) 2123-2133

\noindent [De-GR-Ka-Kun-Kw] K. Deka, F. Garc\'ia-Ramos, K. Kasprzak, P. Kunde and D. Kwietniak, The conjugacy and flip conjugacy problem for Cantor minimal systems,\ \it Preprint, in preparation\ \rm

\noindent [Do-J-K]\ \ R. Dougherty, S. Jackson and A. S. Kechris, The structure of hyperfinite Borel equivalence relations,\ \it Trans. Amer. Math. Soc.\rm 341, 1 (1994) 193-225

\noindent [E]\ \ R. L. Ellis, Extending continuous functions on zero-dimensional spaces,\ \it  Math. Ann.\ \rm 186 (1970), 114-122

\noindent [I-Me]\ \ T. Ibarluc\'ia and J. Melleray, Full groups of minimal homeomorphisms and Baire category methods,\ \it Ergodic Theory Dynam. Systems\ \rm 36, 2 (2016), 550-573

\noindent [G]\ \ S. Gao,~\it Invariant Descriptive Set Theory,~\rm Pure and Applied Mathematics, A Series of Monographs and Textbooks, 293, Taylor and Francis Group, 2009

\noindent [H-N]\ \ P. Hell and J. Nes\v et\v ril,~\it Graphs and homomorphisms,~\rm Oxford Lecture Series in Mathematicsand its Applications, 28, Oxford University Press, Oxford, 2004, xii+244 pp

\noindent [K1]\ \ A. S. Kechris,~\it Classical Descriptive Set Theory,~\rm Springer-Verlag, 1995

\noindent [K2]\ \ A. S. Kechris, On the concept of $\ca$-completeness,~\it Proc. Amer. Math. Soc.
\ \rm 125, 6 (1997), 1811-1814

\noindent [K-S-T]\ \ A. S. Kechris, S. Solecki, and S. Todor\v cevi\' c, Borel chromatic numbers,\ \it Adv. Math.\rm\ 141, 1 (1999), 1-44

\noindent [Kn-Re]\ \ B. Knaster and M. Reichbach, Notion d'homog\'en\'eit\'e et prolongements des hom\'eomorphies,\ \it Fund. Math.\rm\ 40 (1953), 180-193

\noindent [Kr-St]\ \ A. Krawczyk and J. Steprans, Continuous colorings of closed graphs,\ \it Topology Appl.~\rm 51 (1993), 13-26

\noindent [Ku]\ \ P. K\r{u}rka,~\it Topological and symbolic dynamics,~\rm Cours Sp\'ecialis\'es (Specialized Courses), 11, 
Soci\'et\'e Math\'ematique de France, Paris, 2003

\noindent [L-Z]\ \ D. Lecomte and M. Zelen\'y, Baire-class $\xi$ colorings: the first three levels,\ \it Trans. Amer. Math. Soc.\rm\ 366, 5 (2014), 2345-2373

\noindent [L]\ \ D. Lecomte, Continuous 2-colorings and topological dynamics,\ \it  Dissertationes Math.\rm\ 586 (2023), 1-92

\noindent [Lo-R]\ \ A. Louveau and C. Rosendal, Complete analytic equivalence relations,\ \it Trans. Amer. Math. Soc.
\ \rm 357, 12 (2005), 4839-4866

\noindent [Me]\ \ J. Melleray, Dynamical simplices and Borel complexity of orbit equivalence,\ \it Israel J. Math.\ \rm 236 (2020), 317-344

\noindent [P]\ \ Y. Pequignot, Finite versus infinite: an insufficient shift,\ \it Adv. Math.\rm\ 320, 7 (2017), 244-249

\noindent [T-V] S. Todor\v cevi\' c and Z. Vidny\'anszky, A complexity problem for Borel graphs,\ \it Invent. Math.\ \rm 226 (2021), 225-249 

\noindent [vE]\ \ A. J. M. van Engelen,~\it Homogeneous zero-dimensional absolute Borel sets,~\rm CWI Tract, 27. Stichting Mathematisch Centrum, Centrum voor Wiskunde en Informatica, Amsterdam, 1986. iv+133 pp

\end{document}